\documentclass[11pt,reqno]{amsart}

\usepackage[T1]{fontenc}
\usepackage{algorithm} 
\usepackage{algpseudocode} 
\usepackage{float}
\usepackage{enumerate}

\usepackage{subcaption}
\usepackage{caption} 
\usepackage{graphicx} 
\newtheorem{theorem}{Theorem}[section]
\newtheorem{lemma}[theorem]{Lemma}

\newtheorem{proposition}[theorem]{Proposition}

\theoremstyle{remark}
\newtheorem{remark}{Remark}[section]

\usepackage[svgnames]{xcolor}
\usepackage{color}
\usepackage{amssymb, amsmath, amsfonts, latexsym}

\usepackage{amsthm}

\usepackage{enumerate}
\usepackage{bm}
\usepackage{bbm} 
\usepackage{enumitem} 
\usepackage{hyperref}
\hypersetup{colorlinks=true, citecolor = teal, linkcolor=blue!80!black}

\usepackage[initials,alphabetic]{amsrefs}

\newcommand{\dd}{\mathrm{d}}

\newcommand{\supp}{\operatorname{supp}}
\newcommand{\argmin}{\operatorname{argmin}}

\newcommand{\bu}{\boldsymbol{u}}
\newcommand{\bx}{\boldsymbol{x}}
\newcommand{\bX}{\boldsymbol{X}}
\newcommand{\by}{\boldsymbol{y}}
\newcommand{\bY}{\boldsymbol{Y}}
\newcommand{\bz}{\boldsymbol{z}}
\newcommand{\pv}{\operatorname{p.v.}}

\newcommand{\PP}{\mathbb P}

\newcommand{\R}{\mathbb R}

\newcommand{\vol}{\operatorname{vol}}

\allowdisplaybreaks[4]

\newcommand{\8}{\infty}

\newcommand{\ax}{\mathcal{A}}

\newcommand{\ix}{\mathcal{I}}
\newcommand{\kx}{\mathcal{K}}

\newcommand{\mx}{\mathcal{M}}

\newcommand{\rz}{\mathbb{R}}
\newcommand{\cz}{\mathbb{C}}
\newcommand{\ez}{\mathbb{E}}

\newcommand{\pz}{\mathbb{P}}

\newcommand{\al}{\alpha}
\newcommand{\de}{\delta}
\newcommand{\si}{\sigma}

\newcommand{\la}{\lambda}
\newcommand{\bt}{\beta}

\newcommand{\eps}{\varepsilon}

\newcommand{\nc}{\normalcolor}

\numberwithin{equation}{section}
\numberwithin{figure}{section}

\begin{document}

\title[Schatten ball intersection and difference]{Volume asymptotics of Schatten ball intersection and difference}

\author{Siyu Wang}  
\address{Academy of Mathematics and Systems Science, Chinese Academy of Sciences, 100190 Beijing, China.}
\email{wang\_siyu@mail.bnu.edu.cn}
\author{Qiang Zeng}
\address{Academy of Mathematics and Systems Science, Chinese Academy of Sciences, 100190 Beijing, China.}
\email{qzeng@amss.ac.cn}


\date{November 25, 2026}


\begin{abstract}
We establish the first order high dimensional asymptotics for the volumes of intersection and difference of the Schatten $p$-ball and $q$-ball, quantifying the $0$-$1$ law obtained by Kabluchko, Prochno, and Th\"ale in 2020. The analogous question for classical $\ell_p^n$ balls was originally asked by Milman and was considered by Schechtman and Zinn in 1990. At the technical level, our approach to difference volumes can be applied to improve the large deviation principle of empirical measure moment of $\beta$-ensembles due to Augeri in 2018 from even integer monomial potential with convex perturbation to any continuous potential equivalent to power function at infinity. Our solution to the noncommutative problem relies primarily on logarithmic potential theory. 
\end{abstract}

\maketitle

\tableofcontents

\section{Introduction}

The volume of a (convex) body in Euclidean spaces is a classic object to study in asymptotic geometric analysis; see, e.g., \cite{Pis89}.  V. D. Milman raised the question of the volume left in the (unit volume) $\ell_p^n$ ball after removing a $t$-multiple of the $\ell_q^n$ ball for $p=1$ and $q=2$. Here we recall the quasinorm of a vector $x=(x_1,\dots,x_n)$ in $n$-dimensional Euclidean space $\|x\|_{\ell_p}=(\sum_{i=1}^n |x_i|^p)^{1/p}, p>0,$ and the $\ell_p^n$ ball of radius $r$ is the set of $x$ with $\|x\|_{\ell_p}\le r$.
Schechtman and Zinn \cite{SZ90} studied this problem and derived exponential tail bounds with the precise rate $n^{p/q}$ for $0<p<q\le \8$. 
Since the removed part is the intersection of the two balls, the question can be formulated in terms of the intersection or the difference volume. The works \cites{KPT19,Pr26} considered the large dimension limits of the volume.  

In this paper, we consider the noncommutative Milman problem. As the Schatten $p$-class is the noncommutative version of $\ell_p$ space, Schatten $p$-balls of finite dimensional matrices are natural analogs of the $\ell_p^n$ balls. Their asymptotic geometry has been studied from several viewpoints, including volume asymptotics and volume ratios \cites{SR84,KPT20a}, isotropic and concentration properties \cites{KMP98,GP07}, and limit theorems and large deviations \cites{KPT20c,JKP24}; see also \cites{Sim05,PX03} for general background on Schatten $p$-classes and noncommutative $L_p$ spaces. 

Let us introduce the notation. Fix $\beta\in\{1,2,4\}$. Following \cites{KPT20b,KPT20c}, let $\operatorname{Mat}_n(\mathbb F_\beta)$ denote the space of $n\times n$ matrices over $\mathbb F_\beta$, where $\mathbb F_1=\mathbb R$, $\mathbb F_2=\mathbb C$, and $\mathbb F_4=\mathbb H$ (quaternions). For $0<p<\infty$, the Schatten $p$-quasinorm on matrices is defined by $\|A\|_{S_p}=\bigl(\mathrm{ tr} (A^*A)^{p/2}\bigr)^{1/p},$ and $\|A\|_{S_\infty}=\|A\|_{\mathrm{op}}$, the operator norm on $\ell_2$. Here $\mathrm{tr}(A)=\sum_{i} A_{ii}$ is the usual matrix trace, and $A^*$ is the matrix transpose, complex conjugate transpose and quaternion conjugate transpose for $\bt=1,2$ and 4 respectively. When $1\le p\le \8$, $\|\cdot\|_{S_p}$ is a norm.   The Schatten $p$ unit ball is  $$\mathbb{SB}_{p,\beta}^n =\{A\in\operatorname{Mat}_n(\mathbb F_\beta):\|A\|_{S_p}\le1\}.$$
In the self-adjoint setting, let $\mathcal{H}_n\left(\mathbb{F}_\beta\right)=\left\{A \in \operatorname{Mat}_n(\mathbb F_\beta): A^*=A\right\}$. We consider the unit ball
$$\mathbb{B}_{p,\beta}^n =\{A\in\mathcal H_n(\mathbb F_\beta):\|A\|_{S_p}\le1\}.$$
Note that when $p<1$ these Schatten balls are not convex.
We regard both $\operatorname{Mat}_n(\mathbb F_\beta)$ and $\mathcal H_n(\mathbb F_\beta)$ as real Euclidean spaces with scalar product $\langle A,B\rangle=\Re \operatorname{tr}(AB^*).$
Let $\vol_{\beta,n}^{\mathrm{Mat}}$ and
$\vol_{\beta,n}^{\mathcal H}$ denote the corresponding Riemannian volumes, respectively. We write $\vol$ when the ambient space is clear.
With $d_{\beta,n}=\beta n(n-1)/2+n$, the volume normalized balls are
$$\mathbb D_{p,\beta}^n
=\vol(\mathbb B_{p,\beta}^n)^{-1/d_{\beta,n}}
\mathbb B_{p,\beta}^n,
\qquad
\mathbb{SD}_{p,\beta}^n
=\vol(\mathbb{SB}_{p,\beta}^n)^{-1/(\beta n^2)}
\mathbb{SB}_{p,\beta}^n.$$
For $0<p<\infty$, let $b_p=(p\sqrt{\pi} \Gamma(p/2)/\Gamma((p+1)/2))^{1/p}$ and define $b_\infty=\lim_{p\to\infty}b_p$.
For $0<p,q<\infty$, let
\begin{align}\label{eq:rpq}
r_{p,q}
=\frac{b_p}{b_q}
\left(\frac{2p}{p+q}\right)^{1/q}, \qquad
s_{p,q}
=\left(\frac{p+q}{2q}\right)^{1/p}
\frac{b_p}{b_q}=\frac{1}{r_{q,p}},
\end{align}
and
\begin{align}\label{eq:tcpq}
t_{p,q}^*=e^{\frac{1}{2p}-\frac{1}{2q}}
\Bigl(\frac{2p}{p+q}\Bigr)^{1/q}, \qquad c_{p,q}
=
e^{-\frac1{2p}+\frac1{2q}}\frac{b_p}{b_q}.
\end{align}
We extend these constants to $p=\infty$ or $q=\infty$ by limits; their explicit values are given in Section \ref{sec:ball}.
In particular, $c_{p,q}t_{p,q}^*=r_{p,q}$.

\begin{figure}[!ht]
 \centering
 \includegraphics[width=0.6\linewidth]{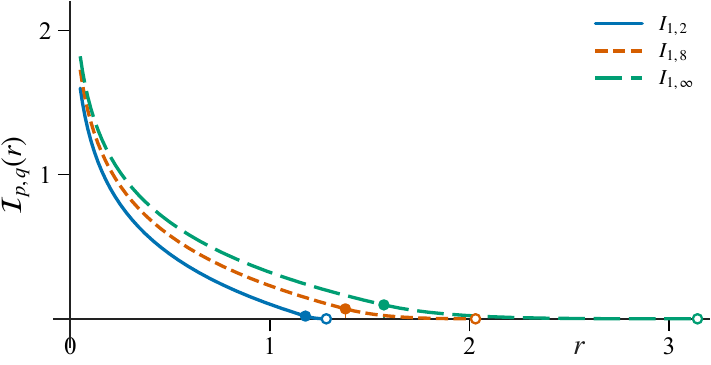}
 \caption{$\mathcal I_{p,q}(r)$ for several choices of $(p,q)$. Solid points indicate the value at $s_{p,q}$, while hollow points indicate $r_{p,q}$.}
 \label{fig:I1q}
\end{figure}

A basic question concerns the asymptotic intersection volume of two such balls. For $0<p,q<\infty$, Kabluchko, Prochno, and Th\"ale \cite{KPT20b}*{Theorem 5.1} established the following $0$-$1$ law: as $n\to\8$
$$
\vol\left(\mathbb D_{p,\beta}^n\cap
t\mathbb D_{q,\beta}^n\right)
 \to 
\begin{cases}
0,&t<t_{p,q}^*,\\
1,&t>t_{p,q}^*.
\end{cases}
$$
This law determines the limiting intersection volume with the convergence rates left open. The goal of this paper is to identify the exact exponential speeds and the corresponding rate functions for the intersection volume in the subcritical regime (i.e., $t<t^*_{p,q}$) and for the complementary difference volume in the supercritical regime (i.e., $t>t^*_{p,q}$), which are quantitative versions of the above 0-1 law.
Let us define the rate function $\mathcal I_{p,q}$ as follows. 
For $0<p,q<\infty$, 
\begin{align}\label{eq:ipqdef}
\mathcal I_{p,q}(r) =
\begin{cases}
\frac12 \log \frac{b_p}{r b_q} +\frac14\left(\frac1q-\frac1p\right),
& 0<r\le s_{p,q},\\
\frac12\log\frac{b_p}{a_{p,q}(r)}
+\frac{p+q}{4pq}
\Bigl(
\frac{b_p^p}{a_{p,q}(r)^p}-\frac{p+2q}{p+q}
\Bigr)^2
-\frac{1}{4p}+\frac{1}{4(p+q)},
& s_{p,q} <r \leq r_{p,q},
\end{cases}
\end{align}
where, for $p\ne q$ and $s_{p,q}<r\leq r_{p,q}$, $a_{p,q}(r)$ is the smaller positive solution to
\begin{align}\label{eq:a-def}
	\frac{b_q^q}{a^q}r^q-\frac{b_p^p}{a^p}=\frac{p-q}{p+q}.
\end{align}
For $0<p<\infty$, $q=\infty$, and $0<r\leq r_{p,\infty}=b_p$, define
\begin{equation}\label{eq:ipinf}
\mathcal I_{p,\infty}(r)
=
\frac12\log\frac{b_p}{r}-\frac{1}{4p}
+
\frac{1}{p}\Bigl(1-\frac{b_p^p}{2r^p}\Bigr)^2
\mathbf{1}_{\{2^{-1/p}b_p<r\leq b_p\}}.
\end{equation}
In Figure \ref{fig:I1q}, we compare these rate functions for several pairs $(p,q)$ and mark the thresholds $s_{p,q}$ and $r_{p,q}$. Further properties of $\mathcal I_{p,q}$ are collected in Lemma \ref{le:apq}. In Remarks \ref{rem:pintq} and \ref{rem:Ilim}, we explain the behaviors of its first branch and its limit as $q\to\infty$, respectively.

The following result establishes the exponential decay rate of the volume asymptotics of Schatten ball intersections for $0<t\leq t_{p,q}^*$. Our method works for arbitrary $p$ and $q$, including the nonconvex balls.

\begin{theorem}\label{th:ballint}
Let $\beta \in \{1,2,4\},$ $0<p<\infty$ and $0<q\leq \infty$. For $0<t \leq t_{p,q}^*$, the following asymptotics hold.
\begin{enumerate}
\item[(i)] \emph{Self-adjoint case:} 
$$ \lim_{n\to\infty}\frac{1}{\beta n^2}
\log \vol\left(\mathbb D_{p,\beta}^n\cap t\mathbb D_{q,\beta}^n\right)
=
-\mathcal I_{p,q}(c_{p,q}t).$$
\item[(ii)] \emph{General case:} 
$$ \lim_{n\to\infty}\frac{1}{\beta n^2}
\log \vol\left(\mathbb{SD}_{p,\beta}^n\cap t\mathbb{SD}_{q,\beta}^n\right)
=
-2\mathcal I_{p,q}(c_{p,q} t).$$
\end{enumerate}
Moreover, $\ix_{p,q}(r)\ge0$ is continuous and strictly decreasing for $0<r\le r_{p,q}$ with $\ix_{p,q}(r_{p,q})=0$.
\end{theorem}

\begin{figure}[htbp]
 \centering
 \begin{subfigure}[t]{.50\linewidth}
  \centering
  \includegraphics[width=\linewidth,viewport=0 0 260 183,clip]{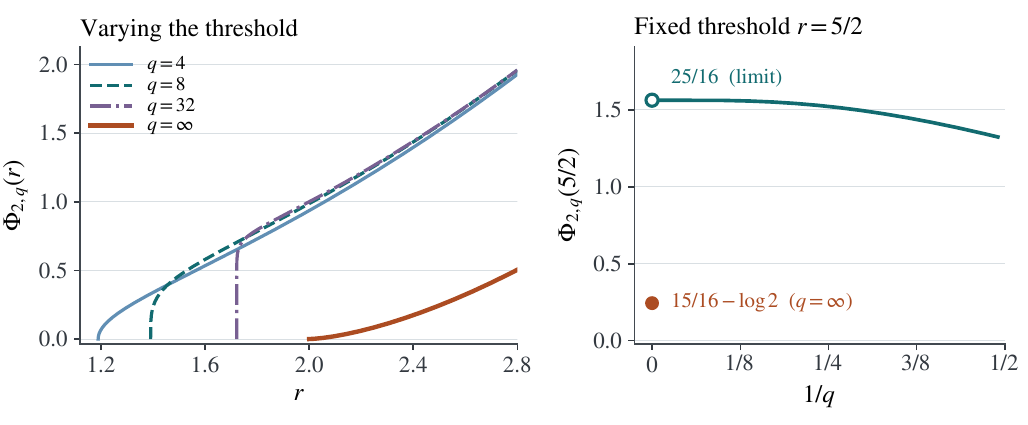}
  \caption{Varying the value $r$.}
 \end{subfigure}\hfill
 \begin{subfigure}[t]{.45\linewidth}
  \centering
  \includegraphics[width=\linewidth,viewport=260 0 493.2 183,clip]{figures/Phi2q.pdf}
  \caption{Varying $1/q$ at $r=5/2$.}
 \end{subfigure}
 \caption{Difference rate function for $p=2$, at speed $\beta n^{1+2/q}$ for finite
 $q$ and $\beta n$ for $q=\infty$. In (B) increasing $q$ moves
 from right to left. At $1/q=0$, the open point is the limit $25/16$ and the filled point is the $q=\8$ rate $15/16-\log 2$.}  \label{fig:Phi2q}
\end{figure}

We remark that the coefficient 2 in the rate function of general case originates from the roles of the symmetric eigenvalues and their squared quantities as singular values in maximizing the free entropy; see the discussion before \cite{HP00}*{Proposition 5.3.5} and the proofs in Section \ref{sec:ball} for more details.
We next consider differences of Schatten balls.
Let $0<p<q\le\infty$. For $r \geq r_{p,q}$, we define
\begin{align}\label{eq:phipq}
\Phi_{p,q}(r)
=\begin{cases}
	\frac{1}{2p}
\left(r^q-r_{p,q}^q\right)^{p/q}, & 0<p<q<\infty, \\
\frac{r^p}{2p}
- 
{}_2F_1\left(
\frac12,\frac p2;\frac p2+1;\frac{b_p^2}{r^2}
\right)\cdot \frac1p
-
\log\frac{r+\sqrt{r^2-b_p^2}}{b_p}, & 0<p<q=\infty,
\end{cases}	
\end{align}
where the hypergeometric function is given by
\begin{equation}\label{e:hypint}
{}_2F_1\left(\frac12,\frac p2;\frac p2+1;z\right)
=\frac p2\int_0^1\frac{u^{p/2-1}}{\sqrt{1-zu}}\dd u,
\qquad p>0,\quad 0\leq z\leq1.
\end{equation}
For fixed $0<r<b_p$,  Figure \ref{fig:I1q} illustrates the convergence $\mathcal I_{p,q}(r)\to\mathcal I_{p,\infty}(r)$ as $q\to \infty$. The difference rate function  $\Phi_{p,q}$ behaves differently. In Remark \ref{rem:Plim}, we explain why $\Phi_{p,q}(r)$ has a phase transition at $q=\infty$; see Figure \ref{fig:Phi2q} for an illustration.

Let $\mathcal M_1(D)$ denote the space of Borel probability measures on a closed set $D\subset \rz$ endowed with the topology of weak convergence. For $\mu\in\mathcal M_1(D)$, the noncommutative (or free) entropy is defined as
\begin{align}\label{def:sigmu}
\Sigma(\mu)= \begin{cases}\iint \log |x-y|\mu( \dd x) \mu(\dd y), & \text { if } \int \log (|x|+1) \mu( \dd x)<\infty, \\ -\infty, & \text { otherwise. }
\end{cases}
\end{align}
The quantity $-\Sigma(\mu)$ is also known as the logarithmic energy of $\mu$. 
For $0<s\leq\infty$, we define the absolute moments of $\mu$ as
\begin{equation}\label{eq:mumom}
    m_s(\mu)=
\begin{cases}
\int |x|^s \mu(\dd x),&0<s<\infty,\\
\sup\{|x|:x\in\supp(\mu)\},&s=\infty.
\end{cases}
\end{equation}
For $0<q<p\leq\infty$ and $r_{p,q} < r<1$, let
\begin{align}\label{eq:phieq}
\Phi_{p,q}(r)
=\begin{cases}
\frac12\log\frac{b_p}{2}-\frac{1}{4p}
-\frac12 \sup_{\mu\in\mathcal M_1(\mathbb R)}\left\{
\Sigma(\mu):m_p(\mu)=1, m_q(\mu)=r^q
\right\}, & p<\infty,\\
-\frac12\log 2-\frac12 \sup_{\mu\in\mathcal M_1(\mathbb R)}\left\{
\Sigma(\mu):
 m_\infty(\mu)\leq1,
m_q(\mu)=r^q
\right\}, & p=\infty.
\end{cases}
\end{align}
The next result gives the exact exponential decay rate in the supercritical regime.

\begin{theorem}
\label{th:balldf}
Let $\beta\in\{1,2,4\}$, and let $t_{p,q}^\star$ and $c_{p,q}$ be defined by \eqref{eq:tcpq}. Assume either that $0<p<q\leq\infty$ and $t> t_{p,q}^*$, or that $0<q<p\leq\infty$ and $t_{p,q}^* < t<c_{p,q}^{-1}$.
Then the following asymptotics hold.
\begin{enumerate}
\item[(i)] Self-adjoint case.
$$
\lim_{n\to\infty}
\frac{1}{\beta n^{1+(p/q)\wedge1}}
\log
\vol\left(
\mathbb D_{p,\beta}^n
\setminus
t\mathbb D_{q,\beta}^n
\right)
=
-\Phi_{p,q}(c_{p,q}t).
$$
\item[(ii)] General case.
$$
\lim_{n\to\infty}
\frac{1}{\beta n^{1+(p/q)\wedge1}}
\log
\vol\left(
\mathbb{SD}_{p,\beta}^n
\setminus
t\mathbb{SD}_{q,\beta}^n
\right)
=
-2\Phi_{p,q}(c_{p,q}t).
$$
\end{enumerate}
Moreover, $\Phi_{p,q}$ is continuous and strictly increasing  with $\Phi_{p,q}(r_{p,q})=0$ and $0<\Phi_{p,q}(r)<\infty$ for $r>r_{p,q}$ in its respective domain.
\end{theorem}

For difference volumes with $q<p$, we only state a variational formula \eqref{eq:phieq} for the rate function.  One may invoke the logarithmic potential theory to find a closed expression for some special values of $p$ and $q$, which usually involves working with the multi-cut equilibrium measures. Since it is quite technical but lacks generality,  we do not include these special cases here. An explicit rate function for arbitrary $p$ and $q$ seems out of reach, as is already the case in the classical $\ell_p^n$ ball setting \cite{KPT19}. Indeed, for $1\leq q<p<\infty$, \cite{KPT19}*{Theorem~1.2} provided variational rates at speed $n$ for subcritical intersection volumes and supercritical difference volumes.
For $1\leq p<q<\infty$, \cite{KPT19}*{Theorem~1.3} established an explicit rate function at speed $n^{p/q}$ for supercritical difference volumes. While Theorem \ref{th:balldf} is the noncommutative analog of the classical setting for the difference volumes, there seems no counterpart of Theorem \ref{th:ballint} for the intersection volume rate function of classical $\ell_p^n$ balls.

It is worth noting that the speed and rate function of volume differences are more complicated than those of volume intersections. Let us comment briefly on the two different parameter regimes. Unlike the case $p<q$, where the large deviation of the set difference is caused by a single outlier, for $p>q$ it results from a collective shift of all eigenvalues and thus occurs at the empirical measure LDP speed $n^2$; see Remark \ref{re:speedd} for a discussion at the technical level. In short, the restriction on the $p$-moment automatically places a constraint on the $q$-moment for $p>q$, allowing only small changes for individual eigenvalues which in turn yields a global deviation of the spectrum, while for $p<q$ a big change in the extreme eigenvalue suffices to cause a large deviation in the $q$-moment.
Moreover, when $p>q$, the limiting rate function has a finite cutoff at $r=1$. Indeed, for every $A\in\mathbb B_{p,\beta}^n$,
$n^{\frac1p-\frac1q}\|A\|_q
\leq
\|A\|_p
\leq 1.$
Hence, for $r\geq1$,
$$
\mathbb B_{p,\beta}^n
\setminus
r n^{\frac1q-\frac1p}\mathbb B_{q,\beta}^n
=\varnothing.
$$
In contrast, when $p<q$, this set is empty only if $r\geq r_n:=n^{\frac1p-\frac1q}.$
Since $r_n\to\infty$, no finite cutoff appears in the limiting rate function. See Propositions \ref{thm:diffball} and \ref{thm:diffpgq} for this distinction.

Let us briefly describe the proof strategy. Our arguments build upon a series of inspiring works of Kabluchko, Prochno, and Th\"ale, including the volume asymptotics \cite{KPT20a}, the Sanov type theorem \cite{KPT20c} and the probabilistic representation \cite{KPT20b} for the Schatten balls. After we obtain the speed and variational formula of the rate function through large deviation techniques, we dig further to solve the optimization problem using the logarithmic potential theory \cite{ST97}. The main novelty of this work is to combine all these theory and techniques together to find the quantitative behavior. 

Among the technical contributions, we single out a large deviation result that is of independent interest. In our proof strategy, it serves as a step towards the difference volume asymptotics.
We first recall the definition of large deviation principle (LDP) (see, e.g., \cite{DZ98}). Let $(X_n)$ be random variables taking values in a Polish space $\mathbb X$, and let $a_n\to\infty$. 
We say that $(X_n)$ satisfies an LDP on $\mathbb X$ with speed $a_n$ and rate function $I:\mathbb X\to[0,\infty]$ if $I$ is lower semicontinuous and for every open set $O\subset\mathbb X$ and every closed set $F\subset\mathbb X$,
$$
\liminf_{n\to\infty}\frac{1}{a_n}\log\PP(X_n\in O)\ge -\inf_{x\in O}I(x),
\quad
\limsup_{n\to\infty}\frac{1}{a_n}\log\PP(X_n\in F)\le -\inf_{x\in F}I(x).
$$
A rate function $I$ is said to be good if the sublevel set $\{x\in X: I(x)\le M\}$ is compact for all $M\in \rz_+$.
The proof of Theorem \ref{th:balldf} for the case $p<q$ relies on the right tail large deviation probability (or LDP) of the absolute $q$-moment of the empirical measure of $\bt$-ensembles. The case $p\ge2$ and $q>p$ being an even integer has been obtained by Augeri \cite{Aug18}*{Theorem 1.1}, which requires the potential to be of the form $V(x)= b|x|^p+w(x)$ and $w(x)$ is a continuous convex function such that $w(x) = o(|x|^p)$ as $|x|\to\8$. It seems that the even integer $q$-moment assumption extends directly to any real absolute $q$-moment for $q>p\ge2$ using the truncation method in \cite{Aug18}, nevertheless the conditions on the convexity of $w$ and $p\ge2$ appear to be inevitable for the concentration of small terms there. 
The following result removes the convexity assumption on $w$ and allows arbitrary real $0<p<q<\infty$, requiring only that $V$ be continuous and
$V(x)/|x|^p\to b\in(0,\infty)$ as $|x|\to\infty$. Our improvement is achieved via estimating the small terms carefully while still using a truncation scheme.

\begin{theorem}\label{th:momldp}
Fix $\beta>0$. Let $(\lambda_1,\ldots,\lambda_n)$ have the joint density \begin{equation}\label{pdf}
\frac{1}{Z_n}e^{-\beta n\sum_{i=1}^n V(\lambda_i)}
\prod_{1\leq i<j\leq n}|\lambda_i-\lambda_j|^\beta
\qquad\text{on }\mathbb R^n,
\end{equation} where $Z_n$ is the normalizing constant and $V$ is continuous on $\rz$ and satisfies $
\lim_{|x|\to\infty}\frac{V(x)}{|x|^p}=b\in(0,\infty)
$ for some $p>0.$ Then for any $q\in(p, +\infty)$, the sequence $\bigl(n^{-1}\sum_{i=1}^n|\lambda_i|^q\bigr)_{n\geq1}$ satisfies an LDP on $\mathbb R$ with speed $\beta n^{1+p/q}$ and good rate function
\begin{equation}\label{eq:momrate}
J_q(t)=
\begin{cases}
b(t-m_q)^{p/q},&t\geq m_q,\\
+\infty,&t<m_q,
\end{cases}
\end{equation}
where $m_q=\int|x|^q\mu_V(\dd x)$ and
$\mu_V
=\argmin_{\mu\in\mathcal M_1(\mathbb R)}
\left\{\int V\dd\mu-\frac12\Sigma(\mu)\right\}$.
\end{theorem}

The rest of paper is organized as follows. After providing necessary preliminary results in logarithmic potential theory in Section \ref{se:lgt}, we investigate the intersection volumes and establish Theorem \ref{th:ballint} in Section \ref{sec:ball}. The difference of Schatten balls is studied and Theorem \ref{th:balldf} is proven in Section \ref{sec:diffball}. Theorem \ref{th:momldp} is derived in the same section near the end of the paper.

\noindent\textbf{Statement on AI use.} This paper was motivated by our investigation on large deviation problems of random matrices. We discovered the connection among Schatten balls, logarithmic potential theory and random matrices by ourselves, which is the main novelty of this paper. We had obtained the main results (Theorems \ref{th:ballint} and \ref{th:balldf}) by spring 2026 along with large deviation results in \cite{WZbeta}, and presented these results in several conferences and seminars. The help of large language models (ChatGPT and DeepSeek) was limited to figure generation, literature search, grammar and proof checks. Occasionally, we also used AI and Mathematica to evaluate some involved integrals, which were subsequently verified by us. After finishing the first version of this paper, we asked ChatGPT to give a review. It brought the paper \cite{Aug18} to our attention. We found that our method of proving Lemma \ref{le:snoutlier} (similar to that for Lemma \ref{le:sv-o}) can be used to improve \cite{Aug18}*{Theorem 1.1} and then we formulated the current Theorem \ref{th:momldp} out of Lemma \ref{le:snoutlier}. Our original proof of Theorem \ref{th:momldp} was directly generalized from that of Lemmas \ref{le:snoutlier} and \ref{le:sv-o}. ChatGPT suggested some simplifications in the proof. We rewrote and included it as the current proof of Theorem \ref{th:momldp} (and Lemma \ref{le:snoutlier} is now its consequence), but kept the original proof of Lemma \ref{le:sv-o} as a record of our efforts. In addition, we tried using ChatGPT to solve the variational problem \eqref{eq:phieq} for $p>q$. It only worked out for some specific values of $p$ and $q$. We did not include these special cases here.

\section{Preliminaries on logarithmic potential theory}\label{se:lgt}
Let us start with some notation. Throughout, for two sequences $(x_n)_{n\ge1}$ and $(a_n)_{n\ge1}$, we write $x_n=O(a_n)$ if $|x_n|\le C|a_n|$ for some $C>0$ and all sufficiently large $n$, and $x_n=o(a_n)$ if $x_n/a_n\to0$.
We write $x_n=\Theta(a_n)$ if there exist constants $c,C>0$ such that $c|a_n|\le |x_n|\le C|a_n|$ for all sufficiently large $n$. Unless specified otherwise, we denote by $C_{a}$ a constant depending only on (the list of parameters) $a$.
Let $D\subset \cz$ be a closed set. We denote by $\mathcal M_1(D)$ the set of Borel probability measures supported on $D$.  For $\mu\in\mathcal M_1(D)$ and $\varepsilon>0$, let
\begin{equation}\label{eq:BLball}
B(\mu,\varepsilon)=\{\nu\in\mathcal M_1(D): d(\mu,\nu)<\varepsilon\},
\end{equation}
where $d$ denotes the bounded Lipschitz metric. 
For a measurable map $T:\cz\to\cz$ and a probability measure $\mu$, we write $T_\#\mu =\mu\circ T^{-1}$ for the push-forward measure of $\mu$ under $T$.

For a lower semicontinuous function $V:D\to\mathbb R\cup\{+\infty\}$, define
\begin{align}\label{defiv}
I_V(\mu)= \begin{cases}\int V \dd \mu-\frac{1}{2} \Sigma(\mu)-C_V, & \text { if } \int V \dd \mu<\infty, \\ \infty, & \text { otherwise. }\end{cases}
\end{align}
where $C_V=\inf_{\mu \in \mathcal M_1(D)}\left\{\int V \dd \mu-\frac{1}{2} \Sigma(\mu)\right\} \in \rz$.

As mentioned above, a key step in our approach is to identify the minimizing measure. To this end, we employ the logarithmic potential theory \cite{ST97} and  record two  results, which are reproduced from \cite{WZbeta} for the reader's convenience. The first provides a characterization of the density and support of the one-dimensional constrained equilibrium measure. It is an immediate consequence of Theorems IV.1.10, IV.1.11, and IV.2.4, together with Corollary IV.1.12 in \cite{ST97}. Recall that $e^{-V}$ is an admissible weight on $D$ if $V:D\to(-\8,\8]$ is lower semicontinuous, $V(z)<\8$ on a set of positive logarithmic capacity and $\lim_{|z|\to\8, z\in D}V(z)-\log |z| = \8$ whenever $D$ is unbounded. In this case, by \cite{ST97}*{Theorem I.1.3} $I_V(\mu)$ admits a unique minimizer $\mu_V$ with compact support and $C_V\in \rz$, where $I_V$ and $C_V$ were defined in \eqref{defiv}.

\begin{lemma}\label{remr}
Let $D=[u_1,u_2]$, where $u_1=-\infty$ and $u_2=\infty$ are allowed, and let $\mu_D=\argmin_{\nu\in\mathcal M_1(D)} I_V(\nu).$ Assume that $w=e^{-V}$ is an admissible weight on $D$, and write $\mathcal S_D=\supp(\mu_D).$ 
Suppose that $\mathcal S_D=[v_1,v_2]\subset D$, that $w>0$ is absolutely continuous on $[v_1,v_2]$, and that $ w'\in L^p([v_1,v_2])$ for some $p>1.$ Then $\mu_D$ is locally absolutely continuous on $(v_1,v_2)$, and its density is given by
\begin{align}\label{eq:densd}
\sigma_D(x)
=
\frac{\pi-\ell(x)}{\pi^2\sqrt{(v_2-x)(x-v_1)}},
\quad
\ell(x):=\pv\int_{v_1}^{v_2}
\sqrt{(v_2-y)(y-v_1)} \frac{V'(y)}{x-y} \dd y.
\end{align}

Assume in addition that either $V$ is convex on $(u_1, u_2)$, or $x V^{\prime}(x)$ increases on $(u_1, u_2)$ when $\left[u_1, u_2\right) \subset[0, \infty)$. Then $\mathcal S_D$ is a compact interval, say $\mathcal S_D=[a,b]\subset [u_1,u_2].$
The support endpoints satisfy
\begin{align}
\frac1\pi\int_a^b V'(x)\sqrt{\frac{x-a}{b-x}}\dd x
&\begin{cases}=1,&b<u_2,\\ \leq1,&b=u_2,\end{cases}\label{eqb}\\
\frac1\pi\int_a^b V'(x)\sqrt{\frac{b-x}{x-a}}\dd x
&\begin{cases}=-1,&a>u_1,\\ \geq-1,&a=u_1.\end{cases}\label{eqa}
\end{align}
If $V$ is convex, these conditions characterize the support $[a,b]$. The candidate measure can also be verified using \eqref{eq:densd} and the equilibrium criterion in
Lemma \ref{lemr}.

Finally, if $D=\mathbb R$, $V$ is even, $V'$ exists on $(0,\infty)$, and $xV'(x)$ is positive and increasing on $(0,\infty)$, then $\mathcal S_D=[-a,a],$ where $a>0$ is characterized by 
\begin{align}\label{eqb2}
 \frac{2}{\pi}\int_0^1\frac{at V'(at)}{\sqrt{1-t^2}} \dd t =1. 
\end{align}
\end{lemma}

The second result gives a necessary and sufficient criterion for the extremal measure. In practice, once a candidate $\mu_\zeta$ is proposed, it suffices to check equilibrium condition \eqref{equcon}.
\begin{lemma}[Theorems I.1.3 and I.3.3 in \cite{ST97}] \label{lemr}
Let $e^{-V}$ be an admissible weight on the closed set $D$ and $\mu_V=\argmin_{\nu\in\mx_1(D)} I_V(\nu)$ the equilibrium measure associated with $V$ on $D$. Suppose that $\sigma\in \mathcal M_1(D)$ has compact support $S:=\supp(\sigma)$ and finite noncommutative entropy.
Then $\sigma=\mu_V$ if and only if there exists a constant $F_V\in\mathbb R$ such that
\begin{equation}\label{equcon}
V(z)- \int_D \log|z-t| \sigma(\dd t)
\begin{cases}
=F_V, & \text{q.e. on } S,\\
\ge F_V, & \text{q.e. on } D.
\end{cases}
\end{equation}
In particular, whenever \eqref{equcon} holds, one necessarily has $S=\supp(\mu_V)$.
\end{lemma}

The function $x\mapsto -\int \log |x-y|\mu(\dd y)$ is called the logarithmic potential of the compactly supported Borel probability measure $\mu$.
By \cite{ST97}*{Remark I.1.7}, measures of finite logarithmic energy do not charge sets of zero logarithmic capacity. Hence the equality in \eqref{equcon} holds $\mu_V$ almost everywhere and
\begin{equation}\label{eq:calsig}
\Sigma(\mu_V)
=
\int_D V(z) \mu_V(\dd z)-F_V.
\end{equation}
Moreover, if $\supp\mu_V=[\alpha,\beta]$ with $\alpha<\beta$,
averaging the equality in \eqref{equcon} against the arcsine law gives
\begin{equation}\label{eq:arcavg}
F_V
=\frac1\pi\int_\alpha^\beta
\frac{V(x)}{\sqrt{(\beta-x)(x-\alpha)}}\dd x
-\log\frac{\beta-\alpha}{4}.
\end{equation}
Here we used
$$
\int_\alpha^\beta
\frac{\log|x-y|}{\pi\sqrt{(\beta-x)(x-\alpha)}}\dd x
=\log\frac{\beta-\alpha}{4},
\qquad y\in[\alpha,\beta].
$$
\begin{remark}[\cite{WZbeta}*{Remark 2.1}]\label{rem:qepw}
In the setting of Lemma \ref{lemr}, assume that $D\subset\mathbb R$ is a closed interval with nonempty interior. 
Suppose that $V$ is finite and continuous in the interior and that its boundary values are its interior limits, possibly $+\infty$. Then the quasi-everywhere relations in \eqref{equcon} actually hold everywhere:
\begin{align}\label{eq:vlogmu}
 V(x)-\int\log|x-y| \mu_V(\dd y)
 \begin{cases}
  =F_V,&x\in \supp(\mu_V),\\
 \geq F_V,&x\in D.
 \end{cases}
\end{align}
\end{remark}

\section{Intersection of Schatten balls} \label{sec:ball}
In this section, we first state the properties of $\mathcal I_{p,q}$ and the four bounds needed for Theorem \ref{th:ballint}.
The auxiliary lemmas and proofs of the propositions are given in Section \ref{s:intaux}; the lemmas are then proved in Section \ref{s:intlem}.

Recall that, for $0< p<\infty$, $
\|A\|_{S_p}^p=\operatorname{tr}(A^*A)^{p/2},$ $
\|A\|_{S_\infty}=\|A\|_{\mathrm{op}}.$ 
For $A\in\mathcal H_n(\mathbb F_\beta)$ with eigenvalues
$\lambda_1(A),\dots,\lambda_n(A)$, we have
$
\|A\|_{S_p}^p=\sum_{j=1}^n |\lambda_j(A)|^p,
$
so
$$
\mathbb B_{p,\beta}^n
=
\left\{
A\in \mathcal H_n(\mathbb F_\beta):
\sum_{j=1}^n |\lambda_j(A)|^p\le 1
\right\}.
$$
Similarly, if $s_1(A),\dots,s_n(A)$ are the singular values of
$A\in\operatorname{Mat}_n(\mathbb F_\beta)$, then
$$
\mathbb{SB}_{p,\beta}^n
=
\left\{
A\in \operatorname{Mat}_n(\mathbb F_\beta):
\sum_{j=1}^n s_j(A)^p\le 1
\right\}.
$$
For $p=\infty$, replace the sums in these two ball definitions by
$\max_{1\leq j\leq n}|\lambda_j(A)|$ and
$\max_{1\leq j\leq n}s_j(A)$, respectively.

The limiting values of the constants in \eqref{eq:rpq} and
\eqref{eq:tcpq} are
\begin{gather}
b_\infty=\lim_{p\to\infty}b_p=1,\qquad
r_{p,\infty}=b_p,\qquad
s_{p,\infty}=2^{-1/p}b_p,\qquad
r_{\infty,q}=2^{1/q}b_q^{-1}, \notag \\
t_{\infty,q}^*=\left(\frac2{\sqrt e}\right)^{1/q},\qquad
t_{p,\infty}^*=e^{1/(2p)},\qquad
t_{\infty,\infty}^*=1, \label{eq:brslim} \\
c_{p,\infty}=e^{-1/(2p)}b_p,\qquad
c_{\infty,q}=e^{1/(2q)}b_q^{-1},\qquad
c_{\infty,\infty}=1.\notag
\end{gather}

\begin{lemma}\label{le:apq}
Let $0<p \neq q<\infty$. The following properties hold:
\begin{enumerate}[label=(\roman*)]
    \item For $s_{p,q}<r<r_{p,q}$, equation \eqref{eq:a-def} has exactly two positive solutions, and the
smaller one, denoted by $a_{p,q}(r)$, satisfies $a_{p,q}(r)\in (b_p\wedge b_qs_{p,q}, b_p\vee b_qs_{p,q}).$
\item $\lim_{r\downarrow s_{p,q}}a_{p,q}(r)=b_qs_{p,q}.$ 
\item At $r=r_{p,q}$, the smaller positive solution is $a_{p,q}(r_{p,q})=b_p$, and $\lim_{r\uparrow r_{p,q}} a_{p,q}(r)= b_p$.
\end{enumerate}
If $0<p<\infty$ and $0<q\leq\infty$, 
$\ix_{p,q}(r)\ge0$ is continuous and strictly decreasing for $0<r\le r_{p,q}$ with $\ix_{p,q}(r_{p,q})=0$.
\end{lemma}

\begin{remark}\label{rem:Ilim}
For every fixed $r\in(0,b_p)$, we have
\begin{equation}\label{eq:Iqlim}
\lim_{q\to\infty}\mathcal I_{p,q}(r)=\mathcal I_{p,\infty}(r).
\end{equation}
Indeed, if $0<r<2^{-1/p}b_p$, then $r<s_{p,q}$ for all sufficiently large $q$, and hence
$$
\mathcal I_{p,q}(r)
=
\frac12\log\frac{b_p}{r b_q}
+\frac14\left(\frac1q-\frac1p\right)
\to
\frac12\log\frac{b_p}{r}-\frac1{4p}
=
\mathcal I_{p,\infty}(r).
$$
If $2^{-1/p}b_p<r<b_p$, then $s_{p,q}<r<r_{p,q}$ for all sufficiently large $q$.
Let $a_{p,q}(r)$ be the smaller positive solution to \eqref{eq:a-def}. Since then $q>p$, one has
$b_q r<a_{p,q}(r)<b_p.$
Equation \eqref{eq:a-def} gives
$$
b_q^q\Bigl(\frac{r}{a_{p,q}(r)}\Bigr)^q
=
\frac{b_p^p}{a_{p,q}(r)^p}-\frac{q-p}{p+q}
>
\frac{2p}{p+q},
$$
which implies
$b_q<\frac{a_{p,q}(r)}{r}
<
b_q\Bigl(\frac{p+q}{2p}\Bigr)^{1/q}.$
Since $b_q\to1$, it follows that $a_{p,q}(r)\to r$, and substituting into the second branch yields
$\mathcal I_{p,q}(r)\to \mathcal I_{p,\infty}(r).$
At $r=2^{-1/p}b_p$, fix $0<u<r<v<b_p$. By Lemma \ref{le:apq},
$\mathcal I_{p,q}(v)\leq\mathcal I_{p,q}(r)\leq\mathcal I_{p,q}(u)$ for all sufficiently large $q$. Taking $q\to\infty$ and using \eqref{eq:Iqlim} already proved at $u$ and $v$, we deduce
$$
\mathcal I_{p,\infty}(v)
\leq\liminf_{q\to\infty}\mathcal I_{p,q}(r)
\leq\limsup_{q\to\infty}\mathcal I_{p,q}(r)
\leq\mathcal I_{p,\infty}(u).
$$
Letting $u\uparrow r$ and $v\downarrow r$ in these inequalities proves
\eqref{eq:Iqlim} at $r=2^{-1/p}b_p$, since $\mathcal I_{p,\infty}$
in \eqref{eq:ipinf} is continuous at this point.
\end{remark}

The next two results provide the upper and lower bounds for self-adjoint matrices.

\begin{proposition}\label{pr:saq-upper}
Let $0<p<\infty$ and $0<q\le \infty$. For every $0<r<r_{p,q}$,
$$
\limsup_{n\to\infty}
\frac{1}{\beta n^2}
\log
\frac{\vol\left(\mathbb B_{p,\beta}^n\cap r n^{\frac1q-\frac1p}\mathbb B_{q,\beta}^n\right)}
{\vol\left(\mathbb B_{p,\beta}^n\right)}
\le
-\mathcal I_{p,q}(r),
$$
where $\mathcal I_{p,q}$ is defined in \eqref{eq:ipqdef}. 
\end{proposition}
For the lower bound, the slightly stronger formulation with $r_n\to r$ will be useful when we apply the result to volume normalized balls.
\begin{proposition}\label{pr:saq-lower}
Let $0<p<\infty$ and $0<q\le \infty$. Let $r_n\to r \in(0,r_{p,q})$. Then
$$
\liminf_{n\to\infty}
\frac{1}{\beta n^2}
\log
\frac{\vol\left(\mathbb B_{p,\beta}^n\cap r_n n^{\frac1q-\frac1p}\mathbb B_{q,\beta}^n\right)}
{\vol\left(\mathbb B_{p,\beta}^n\right)}
\ge
-\mathcal I_{p,q}(r),
$$
where $\mathcal I_{p,q}$ is defined in \eqref{eq:ipqdef}. 
\end{proposition}

We now turn to the general case, which is developed in parallel with the self-adjoint case. 

\begin{proposition}\label{pr:nspqup}
Let $0<p<\infty$ and $0<q\le \infty$. For every $0<r<r_{p,q}$,
	\begin{align*}
\limsup_{n\to \infty}
\frac{1}{\beta n^2}
\log
\frac{\vol\left(\mathbb{SB}_{p,\beta}^n \cap r n^{\frac1q-\frac1p}\mathbb{SB}_{q,\beta}^n\right)}
{\vol\left(\mathbb{SB}_{p,\beta}^n\right)}
\leq
- 2\mathcal I_{p,q}(r).
\end{align*}
\end{proposition}

The matching lower bound is again stated for an arbitrary sequence $r_n\to r$.
\begin{proposition}\label{pr:nspqlw}
Let $0<p<\infty$ and $0<q\le \infty$. Let $r_n\to r$ with $r\in(0,r_{p,q})$. Then
$$
\liminf_{n\to \infty}
\frac{1}{\beta n^2}
\log
\frac{\vol\left(\mathbb{SB}_{p,\beta}^n \cap r_n n^{\frac1q-\frac1p}\mathbb{SB}_{q,\beta}^n\right)}
{\vol\left(\mathbb{SB}_{p,\beta}^n\right)}
\ge
-2\mathcal I_{p,q}(r).
$$
\end{proposition}

As mentioned in the introduction, the prefactor 2 in the general rate function is related to the Ullman distribution and its squared push-forward measure for optimizing the logarithmic energy.

\begin{proof}[Proof of Theorem \ref{th:ballint}]
We first prove (i). Let
\begin{align}\label{e:sarad}
\rho_{q,n}(t)
=
t\left(
\frac{\vol(\mathbb B_{p,\beta}^n)}
{\vol(\mathbb B_{q,\beta}^n)}
\right)^{1/d_{\beta,n}},
\qquad v_n=n^{1/p-1/q}\rho_{q,n}(t).
\end{align}
The volume normalization gives
$
\vol\left(
\mathbb D_{p,\beta}^n\cap t\mathbb D_{q,\beta}^n\right)
={\vol\left(
\mathbb B_{p,\beta}^n\cap v_n n^{1/q-1/p}\mathbb B_{q,\beta}^n
\right)}/{\vol(\mathbb B_{p,\beta}^n)}.
$
For $0<s\leq\infty$, by \cite{KPT20b}*{Theorem~3.1} we have
\begin{align}\label{eq:volkpt}
\lim_{n\to \infty} \frac{ \vol\left(\mathbb{B}_{s, \beta}^n\right)^{1/d_{\beta,n}}}{\left(\pi/\beta\right)^{1/2}e^{3/4} e^{-\frac1{2s}}b_s n^{-\frac1s-\frac12}}
=1.
\end{align}
It follows from \eqref{eq:tcpq} that $v_n\to r:=c_{p,q}t\in(0,r_{p,q}]$ .
If $r<r_{p,q}$, Proposition \ref{pr:saq-lower} applies with $r_n=v_n$. For $0<\varepsilon<r_{p,q}-r$, we have $v_n\leq r+\varepsilon$ for all sufficiently large $n$. Using Proposition \ref{pr:saq-upper} and the monotonicity of the intersection volumes, we deduce
\begin{align*}
-\mathcal I_{p,q}(r)
\leq\liminf_{n\to\infty}\frac1{\beta n^2}
\log\vol\left(
\mathbb D_{p,\beta}^n\cap t\mathbb D_{q,\beta}^n\right)
\leq\limsup_{n\to\infty}\frac1{\beta n^2}
\log\vol\left(
\mathbb D_{p,\beta}^n\cap t\mathbb D_{q,\beta}^n\right)
\leq-\mathcal I_{p,q}(r+\varepsilon).
\end{align*}
By Lemma \ref{le:apq}, $\mathcal I_{p,q}$ is continuous and $\mathcal I_{p,q}(r_{p,q})=0$. 
Letting $\varepsilon\downarrow0$, we have (i) for $r<r_{p,q}$.
At $r=r_{p,q}$, fix $0<u<r_{p,q}$. Since $v_n\geq u$ eventually,
Proposition \ref{pr:saq-lower} and
$\vol(\mathbb D_{p,\beta}^n)=1$ imply that
\begin{align}\label{e:intbd}
-\mathcal I_{p,q}(u)
\leq\liminf_{n\to\infty}\frac1{\beta n^2}
\log\vol\left(
\mathbb D_{p,\beta}^n\cap t_{p,q}^*\mathbb D_{q,\beta}^n\right)
\leq\limsup_{n\to\infty}\frac1{\beta n^2}
\log\vol\left(
\mathbb D_{p,\beta}^n\cap t_{p,q}^*\mathbb D_{q,\beta}^n\right)\leq0.
\end{align}
Letting $u\uparrow r_{p,q}$ in \eqref{e:intbd} gives the limit zero,
completing the proof of (i).

For (ii), let
\begin{equation}\label{e:marad}
\varrho_{q,n}(t)
=
t\left(
\frac{\vol(\mathbb{SB}_{p,\beta}^n)}
{\vol(\mathbb{SB}_{q,\beta}^n)}
\right)^{1/(\beta n^2)}.
\end{equation}
Then
$
\vol\left(
\mathbb{SD}_{p,\beta}^n\cap t\mathbb{SD}_{q,\beta}^n
\right)
={\vol\left(
\mathbb{SB}_{p,\beta}^n\cap\varrho_{q,n}(t)\mathbb{SB}_{q,\beta}^n
\right)}/{\vol(\mathbb{SB}_{p,\beta}^n)}.
$
The volume asymptotics in \cite{KPT20a}*{Theorem~1} and
\cite{DFGZ23}*{Remark~4.8} yield
\begin{align}\label{eq:vossba}
 \lim_{n\to\infty}
\frac{\vol\left(\mathbb{SB}_{s,\beta}^n\right)^{1/(\beta n^2)}}
{\left(\frac{\pi}{2\beta}\right)^{1/2}e^{3/4}e^{-1/(2s)} b_s n^{-1/s-1/2}}
=1.
\end{align}
Hence $n^{1/p-1/q}\varrho_{q,n}(t)\to c_{p,q}t=r$.
Using Propositions \ref{pr:nspqup}
and \ref{pr:nspqlw}, the same argument as for (i) proves (ii) with rate function $2\mathcal I_{p,q}$. 
\end{proof}

\begin{remark}\label{rem:pintq}
The $r\le s_{p,q}$ branch of the rate function $\ix_{p,q}$ is determined by the $q$-ball. Indeed, using the elementary inequality
\[
\Big(\sum_{i=1}^n |x_i|^q\Big)^{1/q}\le \Big(\sum_{i=1}^n |x_i|^p\Big)^{1/p}\le n^{\frac{1}{p}-\frac{1}{q}}\Big(\sum_{i=1}^n |x_i|^q\Big)^{1/q}, \quad 0<p<q\le\8,
\]
we know that $r n^{\frac1q-\frac1p}\mathbb B_{q,\beta}^n
\subset
\mathbb B_{p,\beta}^n$
whenever $r\le 1\wedge n^{\frac1p-\frac1q}.$
For $p\leq q$ and fixed $0<r\leq1$, this inclusion holds for every $n$. The volume asymptotic \eqref{eq:volkpt} yields
\begin{align*}
&\lim_{n\to\infty}
\frac{1}{\beta n^2}
\log
\frac{
\vol\left(
\mathbb B_{p,\beta}^n
\cap
r n^{\frac1q-\frac1p}\mathbb B_{q,\beta}^n
\right)
}{
\vol\left(\mathbb B_{p,\beta}^n\right)
}
=
\lim_{n\to\infty}
\frac{1}{\beta n^2}
\log
\frac{
\vol\left(
r n^{\frac1q-\frac1p}\mathbb B_{q,\beta}^n
\right)
}{
\vol\left(\mathbb B_{p,\beta}^n\right)
}
\\
 ={} &
-\frac12\log\frac{b_p}{r b_q}
-\frac14\left(\frac1q-\frac1p\right).
\end{align*}
For $p\ne q$, the proofs of Propositions \ref{pr:saq-upper} and \ref{pr:saq-lower} extend this limit to $0<r\leq s_{p,q}$.
But $s_{p,q}>1$ for $p<q$ (see Remark \ref{rem:ullman} below), the simple inclusion relation does not provide any information for $r>1$. Moreover, when $p>q$, the inclusion condition requires $r\leq n^{1/p-1/q}\to0$ and thus becomes trivial. 

On the other hand, the variational representation in the proofs of Propositions \ref{pr:saq-upper} and \ref{pr:saq-lower} tells us that
$$
\lim_{n\to\infty}
\frac{1}{\beta n^2}
\log
\frac{
\vol\left(
\mathbb B_{p,\beta}^n
\cap
r n^{\frac1q-\frac1p}\mathbb B_{q,\beta}^n
\right)
}{
\vol\left(\mathbb B_{p,\beta}^n\right)
}
=
\frac12
\sup_{\substack{
\mu\in\mathcal M_1(\mathbb R)\\
m_p(\mu)\le1, m_q(\mu)\le r^q
}}
\Sigma(\mu)-C_p.
$$
Here $C_p$ is defined in \eqref{eq:sarate}. 
For a measurable map $T:\cz\to\cz$ and a probability measure $\mu$, we write $T_\#\mu =\mu\circ T^{-1}$ for the push-forward measure of $\mu$ under $T$.  
For $q<\infty$, by Remark \ref{rem:ullman}, the measure $\mu_r=(x\mapsto rx)_\#\nu^{(q)}$ maximizes $\Sigma$ under $m_q(\mu)\leq r^q$. 
Since $m_q(\nu^{(p)})=r_{p,q}^q$, the relation \eqref{eq:rpq} yields that
$$
m_q(\mu_r)=r^q,\qquad
m_p(\mu_r)=r^pm_p(\nu^{(q)})
=\left(\frac{r}{s_{p,q}}\right)^p\leq1,
\qquad 0<r\leq s_{p,q}.
$$
As a result,
$$
\sup_{\substack{
\mu\in\mathcal M_1(\mathbb R)\\
m_p(\mu)\le1, m_q(\mu)\le r^q
}}
\Sigma(\mu)
=
\sup_{\substack{
\mu\in\mathcal M_1(\mathbb R)\\
m_q(\mu)\le r^q
}}
\Sigma(\mu),
\qquad 0<r\le s_{p,q}.
$$
Together with the volume asymptotics of the $q$-ball, we find
\begin{align*}
&\lim_{n\to\infty}
\frac{1}{\beta n^2}
\log
\frac{
\vol\left(
\mathbb B_{p,\beta}^n
\cap
r n^{\frac1q-\frac1p}\mathbb B_{q,\beta}^n
\right)
}{
\vol\left(\mathbb B_{p,\beta}^n\right)
}
=
\lim_{n\to\infty}
\frac{1}{\beta n^2}
\log
\frac{
\vol\left(
r n^{\frac1q-\frac1p}\mathbb B_{q,\beta}^n
\right)
}{
\vol\left(\mathbb B_{p,\beta}^n\right)
},
\qquad 0<r\le s_{p,q}.
\end{align*}
For $q=\infty$, the constraint becomes $m_\infty(\mu)\leq r$. Then the maximizer under this constraint is the arcsine law on $[-r,r]$ whose $p$-moment is $2r^p/b_p^p$; see \cite{ST97}*{Example I.3.5}. Hence, it satisfies $m_p\leq1$ if and only if $r\leq s_{p,\infty}=2^{-1/p}b_p$.
We conclude that the first branch is determined solely by the constraint of $q$-moment, and $s_{p,q}$ is the value of $r$ where the maximizing measure satisfies $m_p=1$.
\end{remark}

\subsection{Auxiliary lemmas and proofs of propositions}\label{s:intaux}

The proofs of Propositions \ref{pr:saq-upper} and \ref{pr:saq-lower} rely on explicit constrained equilibrium measures. We collect the required computations in the following two lemmas. Their proofs are given in Section \ref{s:intlem}. Recall the moments $m_s(\mu)$ from \eqref{eq:mumom}.

\begin{lemma}\label{lemipq}
Let $s>0$ and $p>0$, and consider the external potential $V(x)=s|x|^p$.
For $a>0,$ let  $\mu_a=\argmin_{\mu \in \mathcal{M}_1([-a, a])}\left\{-\frac{1}{2} \Sigma(\mu) +\int V(x) \mu(\dd x)\right\}.$

\begin{enumerate}
\item[(i)] If $0<a\le a_0:=\frac{b_p}{(2sp)^{1/p}}$, then $\supp(\mu_a)=[-a,a]$ and $\mu_a$ admits a density $\sigma_a$ given by
$$
\sigma_a(x)
=
\frac{1}{\pi\sqrt{a^2-x^2}}\Bigl(1-\frac{2spa^p}{b_p^p}\Bigr)
+\frac{2s pa^p }{ b_p^p} \cdot\frac{ p|x|^{p-1}}{\pi}\int_{|x|}^a \frac{\dd v}{v^{p} \sqrt{a^2-v^2}},
$$
which is a convex combination of the arcsine distribution and the Ullman distribution. Moreover, for every $q>0$,
\begin{align*}
m_q(\mu_a) =\frac{2a^q}{b_q^q}\left(1-\frac{2spqa^{p} }{(q+p)b_p^p}\right),
\quad \text{ and }\quad
\Sigma(\mu_a)=\log \frac{a}{2}-2 s^2 p \frac{a^{2 p}}{b_p^{2 p}}.
\end{align*}
\item[(ii)] If $a>a_0$, then $\mu_a=\mu_{a_0}$, the Ullman distribution with support $[-a_0,a_0]$, which is also the minimizer for the unconstrained problem with $a=\8$. In particular,
$$
m_q(\mu_{a_0})
=
\frac{2p}{q+p}\frac{b_p^q}{b_q^q}(2sp)^{-q/p},
\quad \text{ and }\quad
\Sigma(\mu_a)
=
\log\frac{a_0}{2}
-
\frac{1}{2p}.
$$
\end{enumerate}
\end{lemma}

For $p>0$, the Ullman distribution on $[-1,1]$ has density
\begin{equation}\label{eq:ullman}
h_p(x)=\frac{p}{\pi}
\left(\int_{|x|}^1\frac{t^{p-1}}{\sqrt{t^2-x^2}}\dd t\right)
\mathbf1_{\{|x|\le1\}}(x).
\end{equation}
For $p=2$, $h_p$ is the semicircle density.

\begin{remark}\label{rem:ullman}
Let $0<p<\infty$. Taking $s=1/(2p)$ in Lemma \ref{lemipq} shows that the measure
$$
\nu^{(p)}
:=\argmin_{\mu\in\mathcal M_1(\mathbb R)}
\left\{
\frac1{2p}\int_{\mathbb R}|x|^p \mu(\dd x)
-\frac12\Sigma(\mu)
\right\}
$$
has density
\begin{equation}\label{eq:nueta}
\nu^{(p)}(\dd x)
=\frac1{b_p}h_p\left(\frac{x}{b_p}\right)
\mathbf1_{[-b_p,b_p]}(x)\dd x.
\end{equation}
In particular, for $0<q<\infty$,
$$
m_p(\nu^{(p)})=1,
\qquad
m_q(\nu^{(p)})=r_{p,q}^q,
\qquad
m_p(\nu^{(q)})^{-1}=s_{p,q}^p,
$$
where $r_{p,q}$ and $s_{p,q}$ are defined in \eqref{eq:rpq}. In particular, $s_{p,q}>1$ for $p<q$ and $s_{p,q}<1$ for $p>q$.
This density agrees with \cite{HP00}*{Proposition 5.3.4}.
Moreover, $\nu^{(p)}$ maximizes $\Sigma(\mu)$ over $\mu\in\mathcal M_1(\mathbb R)$ with $m_p(\mu)\le1$.
\end{remark}

\begin{remark}\label{re:her2sd} 

We can draw connection to random matrix theory from Lemma \ref{lemipq}. In fact, thanks to Weyl's integral formula, the volume of $\mathbb{B}_{p,\beta}^n \cap r n^{\frac1q-\frac1p}\mathbb{B}_{q,\beta}^n$ for $p=2,q=\8$ and $\bt=1$ is proportional to 
$$
\int_{ \|\boldsymbol x\|_\8\le rn^{-1/2}, \|\boldsymbol  x\|_2\le 1 } \prod_{i<j}|x_i-x_j|^\bt \dd \boldsymbol  x.
$$
Note that $\|\boldsymbol x\|_2\le1$ if and only if $\frac1n\sum_{i=1}^n (\sqrt{n}x_i)^2\le 1$.
By Laplace's method, finding the volume asymptotics amounts to solving the optimization problem
$$
\sup_{\mu\in \mx_1([-r,r])} \frac12\Sigma(\mu)-s\int x^2\mu(\dd x)
$$
after introducing the Lagrange multiplier $s$. This is exactly the goal of Lemma \ref{lemipq}. It is also the large deviation rate function of the confined probability of Gaussian $\bt$-ensemble considered in \cite{WZbeta}*{Example 5}. 
\end{remark}

\begin{lemma}\label{lemspq}
Let $s,t,p,q>0$ and consider the external potential $V(x)=s|x|^p+t|x|^q.$
Define
$$
\mu_{s,t}=
\argmin_{\mu\in\mathcal M_1(\mathbb R)}
\left\{-\frac12\Sigma(\mu)+\int V(x)\mu(\dd x)\right\}.
$$
Then $\supp(\mu_{s,t})=[-a,a]$, where $a>0$ is uniquely determined by
$$\frac{2sp}{b_p^p}a^{p}+\frac{2tq}{b_q^q}a^{q}=1.$$
Moreover, $\mu_{s,t}$ admits a density $\sigma_{s,t}$ given by
\begin{align*}
\sigma_{s,t}(x)
= \frac{2s p^2|x|^{p-1}a^{p}}{\pi b_p^p}
\int_{|x|}^{a} \frac{\dd w}{w^{p}\sqrt{a^2-w^2}}
+\frac{2t q^2|x|^{q-1}a^{q}}{\pi b_q^q}
\int_{|x|}^{a} \frac{\dd w}{w^{q}\sqrt{a^2-w^2}}.
\end{align*}
For $r>0$, the absolute moments satisfy
$$
m_r(\mu_{s,t})
=
\frac{4sp^2}{r+p}
\frac{a^{r+p}}{b_p^p b_r^r}
+
\frac{4tq^2}{r+q}
\frac{a^{r+q}}{b_q^q b_r^r}
$$
and the noncommutative entropy is
$$
\Sigma(\mu_{s,t})
=
\log\frac{a}{2}
-\frac{2s^2pa^{2p}}{b_p^{2p}}
-\frac{2t^2qa^{2q}}{b_q^{2q}}
-\frac{8stpq}{p+q}
\frac{a^{p+q}}{b_p^p b_q^q}.
$$
\end{lemma}

\begin{lemma}\label{le:logsiri}
Let $\nu$ be a probability measure supported on compact interval $[a,b]$.
Assume that its distribution function $F(x)=\nu((-\infty,x])$ is continuous and strictly increasing on $[a,b]$. 
Set
$\ell_i=F^{-1}\left(\frac{i-\frac12}{n}\right)$,
$
s_i=F^{-1}\left(\frac{i}{n}\right)$ for  $1\le i\le n.
$
Then
$$
\lim_{n\to\8}\frac{1}{n^2}
\sum_{1\le i<j\le n}
\log(\ell_j-s_i)
=
\frac12
\iint \log|x-y| \nu(\dd x)\nu(\dd y).
$$
\end{lemma}

\begin{proof}[Proof of Proposition \ref{pr:saq-upper}]
Let $Z_n$ be a random matrix uniformly distributed on the Schatten unit ball
$\mathbb{B}_{p,\beta}^n$. Then
$$
\mathbb{P}\left(Z_n \in rn^{\frac{1}{q}-\frac{1}{p}}\mathbb{B}_{q,\beta}^n\right)
=
\frac{\vol\bigl(\mathbb{B}_{p,\beta}^n \cap rn^{\frac{1}{q}-\frac{1}{p}}\mathbb{B}_{q,\beta}^n\bigr)}
{\vol\bigl(\mathbb{B}_{p,\beta}^n\bigr)}, \quad 0<q \leq \infty.
$$
Thus it suffices to estimate this probability. 

 By \cite{KPT20c}*{Theorem 1.1}, the empirical measure $\mu_n=\frac{1}{n} \sum_{i=1}^n \delta_{n^{1 / p} \lambda_i\left(Z_n\right)} $ satisfies an LDP on $\mathcal{M}_1(\mathbb{R})$ with speed $\beta n^2$ and good rate function
\begin{equation}\label{eq:sarate}
\mathcal{K}(\mu)=
-\frac{1}{2} \Sigma(\mu)+C_p,
\qquad
C_p=\frac{1}{2}\log\frac{b_p}{2}-\frac{1}{4p},
\end{equation}
for probability measures $\mu$ satisfying $\int |x|^p\mu(\dd x)\le 1$, and $\mathcal{K}(\mu)=+\infty$ otherwise. 
For $0<q<\infty$,
$$
n^{\frac1p-\frac1q}\|\lambda(Z_n)\|_q
=
\Bigl(\frac1n\sum_{i=1}^n |n^{1/p}\lambda_i(Z_n)|^q\Bigr)^{1/q}
=
m_q(\mu_n)^{1/q},
$$
while for $q=\infty$, $
n^{1/p}\|\lambda(Z_n)\|_\infty=\sup\{|x|:x\in \supp(\mu_n)\}.
$
Hence
$
\mathbb P\left(\|\lambda(Z_n)\|_q\le rn^{\frac1q-\frac1p}\right)=\mathbb P(\mu_n \in F_q(r)),
$ 
where
$$F_q(r)=\left\{\mu \in \mathcal{M}_1(\mathbb{R}): m_p(\mu)\le 1, m_q(\mu)\le r^q\right\}, \quad 0<q<\infty,$$
and
$
F_{\infty}(r)=\left\{\mu \in \mathcal{M}_1(\mathbb{R}): \supp(\mu) \subset[-r, r], m_p(\mu)\le 1\right\}.
$
By Fatou's lemma and Skorokhod's representation theorem, $F_q(r)$ is closed for every $0<q \leq \infty$. The LDP for $\mu_n$ yields
\begin{equation}\label{saub}
\limsup_{n \rightarrow \infty} \frac{1}{\beta n^2}
\log \mathbb{P}\left(\mu_n \in F_q(r)\right)
\leq-\inf_{\mu \in F_q(r)} \mathcal{K}(\mu).
\end{equation}
It remains to show that $\inf_{\mu\in F_q(r)}\mathcal K(\mu)=\mathcal I_{p,q}(r)$, with $\mathcal I_{p,q}$ defined in \eqref{eq:ipqdef}.
We treat the cases $q=\infty$ and $0<q<\infty$ separately, and in each case characterize the corresponding minimizing measure, thereby obtaining the value of the infimum.

\noindent{\bf (i) The case $q=\infty$}.
In this case,
\begin{align}\label{defjr}
\inf_{\mu \in F_\infty(r)} \mathcal{K}(\mu)=\inf_{\mu \in \mathcal M_1(\mathbb R)} \left\{ \mathcal{K}(\mu): \supp(\mu) \subset[-r, r], m_p(\mu)\leq1\right\}.
\end{align}

\noindent\emph{(a) The regime $0<r\le \frac{b_p}{2^{1/p}}.$} 
Let $G_\infty(r)=
\{\mu\in\mathcal M_1(\mathbb R): \supp(\mu)\subset[-r,r]\}.$
By \cite{ST97}*{Example I.3.5}, the arcsine distribution
\begin{align*}
\mu_r^{\mathrm{as}}(\dd x)=\frac{1}{\pi \sqrt{r^2-x^2}} \mathbf{1}_{\{|x|<r\}} \dd x
\end{align*}
 uniquely minimizes $\mathcal K$ over $G_\infty(r)$. Since
$
m_p(\mu_r^{\mathrm{as}})=\frac{2r^p}{b_p^p}\le1,
$ 
we have $\mu_r^{\mathrm{as}}\in F_\infty(r)\subset G_\infty(r)$. Therefore,
$\inf_{\mu\in F_\infty(r)}\mathcal K(\mu)
=
\mathcal K(\mu_r^{\mathrm{as}}).$
Trigonometric substitution yields $\Sigma(\mu_r^{\mathrm{as}})=\log r-\log2.$
Hence, 
\begin{align} \label{eq:calk1}
\inf_{\mu \in F_\infty(r)} \mathcal{K}(\mu)= \mathcal K(\mu_r^{\mathrm{as}})=-\frac12 \Sigma(\mu_r^{\mathrm{as}})+C_p=\frac{1}{2} \log \frac{b_p}{r}-\frac{1}{4 p}.	
\end{align}

\noindent\emph{(b) The regime $2^{-1 / p} b_p<r<b_p.$} 
 In this case, the variational problem \eqref{defjr} reduces to
$$\inf_{\mu \in F_\infty(r)} \mathcal{K}(\mu)=\inf_{\mu \in \mathcal M_1(\mathbb R)} \left\{-\frac{1}{2} \Sigma(\mu)+C_p: \supp(\mu) \subset[-r, r], m_p(\mu)=1\right\}.$$
Indeed, suppose that a minimizer $\mu^*$ of \eqref{defjr} satisfies $m_p(\mu^*)<1.$ Then there exists $\eta>0$ such that the convex combination $\mu_\eta:=(1-\eta) \mu^*+\eta \mu_r^{\mathrm{as}}$ still satisfies $m_p(\mu_\eta) \leq 1$. 
Since $\mu_r^{\mathrm{as}}$ has a logarithmic energy strictly smaller than any other probability measure supported on $[-r,r]$, we have $ \Sigma(\mu_r^{\mathrm{as}})>\Sigma(\mu^*).$ By the strict concavity of $\Sigma$ (see, e.g., \cite{AGZ10}), 
\begin{align*}
- \Sigma(\mu_\eta)<- (1-\eta)\Sigma(\mu^*)-\eta\Sigma(\mu_r^{as})<-\Sigma(\mu^*),
\end{align*}
which contradicts the minimality of $\mu^*$.

In the following, we employ the method of Lagrange multipliers. For $s\ge0$, define
$$
\mathcal L_s^p(\mu)
=
-\frac12\Sigma(\mu)+s\bigl(m_p(\mu)-1\bigr)+C_p,
$$
and denote the unique minimizer by
$\mu_s^p=\argmin\Bigl\{
\mathcal L_s^p(\mu): \supp(\mu)\subset[-r,r]
\Bigr\}.$ 
Applying Lemma \ref{lemipq}(i) with external field $V(x)=s|x|^p$, we find that if $r^p\le \frac{b_p^p}{2s p}$,
then $\supp(\mu_{s}^p)=[-r,r]$, and 
$$m_p(\mu_s^p)
=
\frac{2r^p}{b_p^p}
-
\frac{2s p r^{2p}}{b_p^{2p}}.$$
Imposing the constraint $m_p(\mu_s^p)=1$ gives the unique multiplier
$
s^*
=
\frac{b_p^p}{pr^p}
-
\frac{b_p^{2p}}{2pr^{2p}},
$ 
and one checks that $s^*\ge0$ and
$
r^p\le \frac{b_p^p}{2s^*p}
$
for $b_p/2^{1/p}<r<b_p$. Hence the above formula applies to $\mu_{s^*}^p$, and in particular $m_p(\mu_{s^*}^p)=1.$

Since $m_p(\mu_{s^*}^p)=1$ and $\mu_{s^*}^p$ minimizes $\mathcal L_{s^*}^p$ over measures supported on $[-r,r]$, every $\mu\in F_\infty(r)$ satisfies
$
\mathcal K(\mu_{s^*}^p)
=\mathcal L_{s^*}^p(\mu_{s^*}^p)
\le\mathcal L_{s^*}^p(\mu)
=\mathcal K(\mu)+s^*(m_p(\mu)-1)
\le\mathcal K(\mu).
$
Thus $\mu_{s^*}^p$ minimizes $\mathcal K$ over $F_\infty(r)$. 
Applying Lemma \ref{lemipq}(i) again, we obtain
$$\Sigma\left(\mu_{ s^*}^p\right)=\log \frac{r}{2}-2 s^{* 2} p \frac{r^{2 p}}{b_p^{2 p}}=\log \frac{r}{2}-\frac{2}{p}+\frac{2 b_p^p}{p r^p}-\frac{b_p^{2 p}}{2 p r^{2 p}}.$$
By the above computation,
\begin{align}\label{eq:calk2}
\inf_{\mu \in F_{\infty}(r)} \mathcal{K}(\mu) =\mathcal{K}(\mu_{ s^*}^p) =-\frac{1}{2} \Sigma(\mu_{s^*})+C_p = \frac{1}{2} \log \frac{b_p}{r} +\frac{3}{4 p}-\frac{b_p^p}{p r^p}+\frac{b_p^{2 p}}{4 p r^{2 p}}.	
\end{align}

\noindent{\bf (ii) The case $0<q<\infty$}. In this case
\begin{align}\label{eq:inpq}
\inf_{\mu \in F_q(r)} \mathcal{K}(\mu)
= 
\inf_{\mu \in \mathcal M_1(\R)}\Bigl\{
-\tfrac12\Sigma(\mu)+C_p:
m_p(\mu)\le1, m_q(\mu)\le r^q
\Bigr\}.	
\end{align}

\noindent\emph{(a) The regime $0<r\le s_{p,q}$.} 
Let $\nu^{(p)}$ be the probability measure defined by \eqref{eq:nueta}. By Remark \ref{rem:ullman}, $\nu^{(p)}$ maximizes $\Sigma$ under $m_p\le1$ and
$m_p(\nu^{(p)})=1$, 
$ m_q(\nu^{(p)})=r_{p,q}^q>r^q.$
Arguing as in part (i)(b), the variational problem \eqref{eq:inpq} reduces to 
$$
\inf_{\mu \in F_q(r)} \mathcal{K}(\mu)=
\inf_{\mu \in \mathcal{M}_1(\mathbb{R})}\Bigl\{
-\tfrac12\Sigma(\mu)+C_p:
m_p(\mu)\le1,
m_q(\mu)=r^q
\Bigr\}.
$$
Indeed, otherwise one could construct a convex combination
$\mu_\eta=(1-\eta)\mu^*+\eta\nu^{(p)}$ that yields a strictly smaller value,
contradicting optimality.

For fixed $t\ge0$, define
$
\mathcal L_t^q(\mu)
=
-\frac12\Sigma(\mu)+t\bigl(m_q(\mu)-r^q\bigr)+C_p,
$
and let
$\mu_t^q
=
\argmin_{\mu\in\mathcal M_1(\mathbb R)}
\mathcal L_t^q(\mu).$
By Lemma \ref{lemipq}(ii), the unique minimizer $\mu_t^q$ is the Ullman distribution with support
$$
\supp(\mu_t^q)=\Bigl[-\frac{b_q}{(2tq)^{1/q}},\frac{b_q}{(2tq)^{1/q}}\Bigr],
$$
and satisfies
$$
m_q(\mu_t^q)=\frac{1}{2tq},
\quad
\Sigma(\mu_t^q)=
\log\frac{b_q}{2(2tq)^{1/q}}-\frac{1}{2q}.
$$
Imposing $m_q( \mu_t^q)=r^q$ yields $t^*=\frac{1}{2r^q q}$.
A direct computation shows that if $r\le s_{p,q},$
$$
m_p(\mu_{t^*}^q)
=
\frac{2q}{p+q}\frac{b_q^p}{b_p^p}r^p
\le1.
$$
For $r\le s_{p,q},$ the same argument as in part (i)(b) yield that
$
\mathcal{L}_{t^*}^q\left(\mu_{t^*}^q\right)=\mathcal{K}(\mu_{t^*}^q)= \inf_{\mu \in F_q(r)} \mathcal{K}(\mu).
$
Using $\Sigma(\mu_{t^*}^q)=\log \frac{r b_q}{2} -\frac{1}{2q},$
we find
\begin{align}\label{eq:calk3}
	\inf_{\mu \in F_q(r)} \mathcal{K}(\mu)=
 \mathcal{K}(\mu_{t^*}^q)
=
-\tfrac12\Sigma(\mu_{t^*}^q)+C_p
=
\frac12\log \frac{b_p}{b_q}
+\frac14\left(\frac1q-\frac1p\right)
-\frac12\log r.
\end{align}

\noindent\emph{(b) The regime $s_{p,q}<r<r_{p,q}$.} In this regime, $m_p(\mu_{t^*}^q)>1$ and $\mu_{t^*}^q$ is a minimizer of $\kx(\mu)$ with $m_q(\mu)\le r^q$.
Arguing as in part (i)(b), 
$$
\inf_{\mu \in F_q(r)} \mathcal{K}(\mu)=\inf_{\mu \in \mathcal{M}_1(\mathbb{R})}\Bigl\{
-\tfrac12\Sigma(\mu)+C_p:
m_p(\mu)=1,
m_q(\mu)=r^q
\Bigr\}.
$$
Introducing Lagrange multipliers, we consider for $s,t>0,$
$$\mu_{s,t}:= \argmin_{\mu \in \mathcal M_1(\mathbb R)} \left\{ {-\frac{1}{2} \Sigma(\mu)+s\int|x|^p\dd \mu +t\int|x|^q\dd \mu +C_p} \right\}.$$
By Lemma \ref{lemspq}, the measure $\mu_{s,t}$ admits a density $\sigma_{s,t}$
supported on a compact interval $[-a,a]$, where $a>0$ satisfies the boundary
condition
$$
2sp\frac{a^{p}}{b_p^p}
+
2tq\frac{a^{q}}{b_q^q}
=1.
$$
Imposing the moment constraints
$m_p( \mu_{s,t})=1$
and
$m_q( \mu_{s,t})=r^q,$
the parameters $s$ and $t$ are determined by the system of equations
\begin{align*}
2 sp \frac{a^{2p}}{b_p^{2p}}
+
\frac{4tq^2}{p+q}
\frac{a^{p+q}}{b_p^p b_q^q}
=1, 
\quad
2 tq \frac{a^{2q}}{b_q^{2q}}
+
\frac{4sp^2}{p+q}
\frac{a^{p+q}}{b_p^p b_q^q}
=r^q.
\end{align*}
Solving this system together with the boundary condition gives
$$
s^*=\frac{b_p^{p}(p+q)}{2p(p-q)^2a^{2p+q}}
\Bigl(
2pa^p b_q^q r^q
-2(p-q)a^{p+q}
-(p+q)b_p^p a^q
\Bigr)=
\frac{b_p^{p}(p+q)}{2p(p-q)a^{2p}}
\Bigl(
b_p^p-\frac{2qa^p}{p+q}
\Bigr),
$$
$$
t^*=\frac{b_q^{q}(p+q)}{2q(p-q)^2a^{p+2q}}
\Bigl(
2qa^q b_p^p
+2(p-q)a^{p+q}
-(p+q)a^p b_q^q r^q
\Bigr)=
\frac{b_q^{q}(p+q)}{2q(p-q)a^{p+q}}
\Bigl(
a^p-b_p^p
\Bigr),
$$
where $a>0$ satisfies \eqref{eq:a-def}.
Hence the nonnegativity condition $s^*,t^*\geq0$ is equivalent to $a\in
\bigl[b_p\wedge b_qs_{p,q}, b_p\vee b_qs_{p,q}\bigr].$
By Lemma \ref{le:apq}, the smaller positive solution $a_{p,q}(r)$ lies in the corresponding open interval. 
Therefore the Lagrange multiplier problem selects $a=a_{p,q}(r).$

The same argument as in part (i)(b) yields $\mathcal K(\mu_{s^*,t^*})=\mathcal L_{s^*,t^*}(\mu_{s^*,t^*})
\le\mathcal L_{s^*,t^*}(\mu)=\mathcal K(\mu).$ Consequently,
\begin{align}\label{eq:calk4}
\inf_{\mu \in F_q(r)} \mathcal{K}(\mu)& = \mathcal{K}(\mu_{s^*,t^*}) 
=
-\tfrac12\Sigma(\mu_{s^*,t^*})+C_p\notag\\
&=
\frac12\log\frac{b_p}{a}
+\frac{p+q}{4pq}\Bigl(\frac{b_p^p}{a^p}-\frac{p+2q}{p+q}\Bigr)^2
-\frac{1}{4p}
+\frac{1}{4(p+q)},
\end{align}
where we used
$$
\Sigma(\mu_{s^*,t^*})
=
\log\frac a2
-\frac{p+q}{2pq}\Bigl(\frac{b_p^p}{a^p}-\frac{p+2q}{p+q}\Bigr)^2
-\frac{1}{2(p+q)}.
$$
Altogether, \eqref{eq:calk1}, \eqref{eq:calk2}, \eqref{eq:calk3}, and \eqref{eq:calk4} show that
$\inf_{\mu\in F_q(r)}\mathcal K(\mu)=\mathcal I_{p,q}(r).$
Combining this identity with \eqref{saub}, we obtain the desired upper bound and complete the proof.
\end{proof}

\begin{proof}[Proof of Proposition \ref{pr:saq-lower}]
Let $Z_n$ be uniformly distributed on $\mathbb{B}_{p,\beta}^n$, with ordered eigenvalues $\la_1(Z_n)\le\dots\le \la_n(Z_n)$. Let $\pi$ be an independent uniform random permutation of $\{1,\dots,n\}$. By \cite{KPT20b}*{Corollary 4.3},
\begin{align}\label{eq:lamzn}
 \bigl(\lambda_{\pi(1)}(Z_n),\ldots,\lambda_{\pi(n)}(Z_n)\bigr)
\stackrel d= U^{1/d_{\beta,n}}\frac{\bX}{\|\bX\|_p},
\end{align}
where $U$ is uniformly distributed on $[0,1]$ and, independent of $U$, the vector $\bX=(X_1,\dots,X_n)$ has joint density
$$ Z_{n,p,\beta}^{-1}
\exp\left(-\sum_{i=1}^n|x_i|^p\right)
\prod_{1\le i<j\le n}|x_i-x_j|^\beta
$$
with normalizing constant $Z_{n,p,\bt}$. Change of variables $x_i=(\beta n)^{1/p}y_i$ gives for $0<q<\infty$,
\begin{align*}
&\mathbb P\left(Z_n\in r_n n^{\frac1q-\frac1p}\mathbb B_{q,\beta}^n\right)=\ez\left[\pz\left(U\le \Big(\frac{r_n n^{-1/p}\|\bX\|_p}{n^{-1/q}\|\bX\|_q}\Big)^{d_{\bt,n}}\Big|\bX\right)\right]\\
={} &\frac{1}{\widetilde{Z}_{n, p, \beta}} \int_{\mathbb{R}^n} \min \Big\{\Big(\frac{r_n n^{-1 / p}\|\by\|_p}{n^{-1 / q}\|\by\|_{q}}\Big)^{d_{\beta, n}}, 1\Big\} 
e^{-\beta n \sum_{i=1}^n\left|y_i\right|^p} \prod_{1 \leq i<j \leq n}\left|y_i-y_j\right|^\beta \dd \by,
\end{align*}
where $\by=(y_1,\dots,y_n)$, $\dd \by=\dd y_1\dots\dd y_n$ and
\begin{align*}
\widetilde{Z}_{n, p, \beta}=\frac{Z_{n, p, \beta}}{(\beta n)^{\frac{n}{p} +\frac{\beta n(n-1)}{2 p}}}=\int_{\mathbb{R}^n} e^{-\beta n \sum_{i=1}^n\left|y_i\right|^p} \prod_{1 \leq i<j \leq n}\left|y_i-y_j\right|^\beta \dd \by.
\end{align*}
For $q=\infty$, the same formula holds with $n^{-1/q}\|\by\|_q$ replaced by $\|\by\|_\infty$.
Using the asymptotic formula for $\log Z_{n,p,\beta}$ from \cite{KPT20c}*{Eq.~(7)}, we obtain
\begin{align}\label{logtz}
\lim_{n \rightarrow \infty}\left(-\frac{1}{\beta n^2} \log \widetilde{Z}_{n, p, \beta}\right)=\frac{1}{2}\left(1+\frac{1}{p}\right) \log 2+\frac{3}{4 p}-\frac{1}{2} \log b_p+\frac{1}{2 p} \log p.
\end{align}

We now take $\mu^*$ to be one of the minimizing measures arising in
\eqref{eq:calk1}, \eqref{eq:calk2}, \eqref{eq:calk3}, and \eqref{eq:calk4},
namely,
$\mu_r^{\mathrm{as}}, 
\mu_{s^*}^p, 
\mu_{t^*}^q, 
\mu_{s^*,t^*},$
respectively. 
Let $c=\bigl(2pm_p(\mu^*)\bigr)^{1/p},$ and $\nu^*(A)=\mu^*(cA).$
Then
$$
m_p(\nu^*)=\frac{1}{2p},
\qquad
m_q(\nu^*)
=
\frac{m_q(\mu^*)}
{\bigl(2pm_p(\mu^*)\bigr)^{q/p}},
\qquad
\Sigma(\nu^*)
=
\Sigma(\mu^*)-\frac{1}{p}\log\bigl(2pm_p(\mu^*)\bigr).
$$
As we shall see, $\nu^*$ leads to the correct lower bound of $\mathbb P(Z_n\in r_n n^{\frac1q-\frac1p}\mathbb B_{q,\beta}^n)$ as the empirical measure of $\boldsymbol{y}$. 
Motivated by \eqref{eq:lamzn}, we may guess that $\nu^*$ differs from $\mu^*$ by a rescaling factor and then the value of $c$ can be found by maximizing $\frac12\Sigma(\mu^*(c\cdot))-m_p(\mu^*(c\cdot))$.  
After this rescaling, the measures $\nu^*$ under consideration have densities of the following explicit forms: 
$$
\frac{1}{\pi\sqrt{a^2-x^2}},
$$
or
$$
\sigma_a(x)
=
\frac{1}{\pi \sqrt{a^2-x^2}}
\left(1-\frac{2spa^p}{b_p^p}\right)
+
\frac{2sp^2|x|^{p-1}a^p}{\pi b_p^p}
\int_{|x|}^a
\frac{1}{v^{p}\sqrt{a^2-v^2}} \dd v,
$$
or
$$
\sigma_{s,t}(x)
=
\frac{2sp^2|x|^{p-1}a^p}{\pi b_p^p}
\int_{|x|}^a
\frac{\dd w}{w^p\sqrt{a^2-w^2}}
+
\frac{2tq^2|x|^{q-1}a^q}{\pi b_q^q}
\int_{|x|}^a
\frac{\dd w}{w^q\sqrt{a^2-w^2}}.
$$
The possible singularities of these densities occur only at the endpoints $\pm a$ and possibly at the origin. 
At the endpoints, the worst singularity is of order $(a^2-x^2)^{-1/2}$; at the origin, the term with exponent $r-1$ for $0<r<1$ has at worst the behavior $|x|^{r-1},$ or a logarithmic singularity when $r=1$. These singularities are integrable.
Consequently, there exist constants $C>0$ and $\alpha>0$, depending only on the density, such that for every interval $J\subset[-a,a]$ with Lebesgue measure $|J|<\frac12$,
\begin{align}\label{eq:cjb}
\nu^*(J)\le C |J|^\alpha.	
\end{align}
Indeed, the endpoint singularity contributes at most $C|J|^{1/2}$, while the
singularity at the origin contributes at most $C|J|^r$ if $0<r<1$, and at
most $C|J||\log |J||$ if $r=1$, which is bounded by $C_\alpha |J|^\alpha$ for any $\alpha<1$.

Let
$F(x)=\nu^*((-\infty,x])$
be the distribution function of $\nu^*$, and let
$$
F^{-1}(u)=\inf\{x\in \rz:F(x)\ge u\},
\qquad u\in[0,1],
$$
which coincides with the usual inverse function of $F$ as it is strictly increasing on the support of $\nu^*$. Let $T=\sup\{|x|: x\in \supp(\mu^*)\}.$ The support of $\nu^*$ is exactly the interval $[-T/c,T/c]$. For $1\le i\le n$, define $
\ell_i=F^{-1}\left(\frac{i-\frac12}{n}\right)$,
$s_i=F^{-1}\left(\frac{i}{n}\right).
$
Then $\ell_1<s_1<\ell_2<s_2<\cdots<\ell_n<s_n.$ 
Let $\Delta_n=\{y\in\mathbb R^n: \ell_i\le y_i\le s_i, 1\le i\le n\}.$ 
Then for every $y\in\Delta_n$ and $s>0$,
$$
L_{s,n}:=\frac1n\sum_{i=1}^n \inf_{x\in[\ell_i,s_i]}|x|^s\le \frac1n\sum_{i=1}^n |y_i|^s\le \frac1n\sum_{i=1}^n \sup_{x\in[\ell_i,s_i]}|x|^s=:U_{s,n}.
$$
In particular, if $0<q<\infty$, then
$$
\frac{r_n n^{-1/p}\|\by\|_p}{n^{-1/q}\|\by\|_q}
\ge
\frac{r_n L_{p,n}^{1/p}}{U_{q,n}^{1/q}}
=:W_n^{(q)}.
$$
For $q=\infty$, since $\Delta_n\subset[-T/c,T/c]^n$, we have
$$
\|\by\|_\infty\le T/c,
\quad
\frac{r_n n^{-1/p}\|y\|_p}{\|\by\|_\infty}
\ge
\frac{r_n L_{p,n}^{1/p}}{T/c}
=:W_n^{(\infty)}.
$$
Moreover,
$$
\prod_{1\le i<j\le n}|y_i-y_j|^\beta
\ge
\prod_{1\le i<j\le n}(\ell_j-s_i)^\beta,
\quad y\in\Delta_n.
$$
Restricting the integral to $\Delta_n$, we obtain for $0<q<\infty$,
$$
\mathbb P\left(Z_n\in r_n n^{\frac1q-\frac1p}\mathbb B_{q,\beta}^n\right)
\ge
\frac{1}{\widetilde Z_{n,p,\beta}}
\min\{(W_n^{(q)})^{d_{\beta,n}},1\}
\prod_{1\le i<j\le n}(\ell_j-s_i)^\beta
e^{-\beta n^2 U_{p,n}}
\prod_{i=1}^n (s_i-\ell_i),
$$
and the analogous bound with $W_n^{(\infty)}$ in the case $q=\infty$.

For each $i$, $\nu^*([\ell_i,s_i])=\frac{1}{2n}.$ Using \eqref{eq:cjb}, we get
$(2n)^{-1}
=
\nu^*([\ell_i,s_i])
\le
C(s_i-\ell_i)^\alpha,$
and therefore
$s_i-\ell_i
\ge
C^{-1/\alpha}(2n)^{-1/\alpha}.$
Since $s_i-\ell_i\le 2a$,
$$
\log(2a)n\geq \sum_{i=1}^n\log(s_i-\ell_i)
\ge
-C_\al n\log n.
$$
Consequently,
$
\frac{1}{\beta n^2}\sum_{i=1}^n \log(s_i-\ell_i)\to0.
$
Next, the empirical measures $ \frac{1}{n} \sum_{i=1}^n \delta_{\ell_i}$ and $\frac{1}{n} \sum_{i=1}^n \delta_{s_i}$ both converge weakly to $\nu^*$. Since the support of $\nu^*$ is compact, for every $s>0$ the function $x\mapsto |x|^s$ is continuous and bounded on $[-a,a]$. Hence
$$
\frac1n\sum_{i=1}^n |\ell_i|^s\to m_s(\nu^*),
\qquad
\frac1n\sum_{i=1}^n |s_i|^s\to m_s(\nu^*).
$$
Noting that $F^{-1}$ is uniformly continuous, we have
$\lim_{n\to\8}\max_{1\le i\le n}(s_i-\ell_i)=0.$
Using the continuity of $x\mapsto |x|^s$ on $[-a,a]$, we derive $U_{s,n}-L_{s,n}\to0$ and thus
$$
L_{p,n}\to m_p(\nu^*),
\qquad
U_{p,n}\to m_p(\nu^*),
\qquad
U_{q,n}\to m_q(\nu^*).
$$
By Lemma \ref{le:logsiri}, we obtain
$$
\frac{1}{\beta n^2}\sum_{1\le i<j\le n}\log(\ell_j-s_i)^\beta
\to
\frac12
\iint \log|x-y| \nu^*(\dd x)\nu^*(\dd y)
=
\frac12\Sigma(\nu^*).
$$
Therefore,
$$
\lim_{n\to\infty}\frac1{\beta n^2}\log\min\{(W_n^{(q)})^{d_{\beta,n}},1\}
=
\frac12\log\min\left\{
r\frac{m_p(\nu^*)^{1/p}}{m_q(\nu^*)^{1/q}},
1\right\},
\qquad 0<q<\infty,
$$
while for $q=\infty$,
$$
\lim_{n\to\infty}\frac1{\beta n^2}\log\min\{(W_n^{(\infty)})^{d_{\beta,n}},1\}
=
\frac12\log\min\left\{
r\frac{m_p(\nu^*)^{1/p}}{T/c},
1\right\}.
$$
These limits and \eqref{logtz} imply
\begin{align*}
  & \liminf_{n\to\infty}\frac{1}{\beta n^2}
\log
\mathbb{P}\left(
Z_n\in r_n n^{1/q-1/p}\mathbb{B}_{q,\beta}^n
\right) \\
\ge{} &
- \lim_{n \rightarrow \infty}\frac{ \log \widetilde{Z}_{n, p, \beta}}{\beta n^2}
+\frac{1}{2} \log \min \left\{r \frac{m_p\left(\nu^*\right)^{1 / p}}{m_q\left(\nu^*\right)^{1 / q}}, 1\right\}
+\frac{1}{2}\Sigma\left(\nu^*\right) - m_p\left(\nu^*\right) \\
\ge{} &
\frac12\Sigma(\mu^*)
-\frac{1}{2p}\log m_p(\mu^*)
+\frac12
\log
\min\left\{
r\frac{m_p(\mu^*)^{1/p}}{m_q(\mu^*)^{1/q}},
1
\right\}
+\frac12\log \frac{2}{b_p}
+\frac{1}{4p}. 
\end{align*}
In the case $q=\infty$, the term $m_q(\mu^*)^{1/q}$ is replaced by $\sup\{|x|:x\in \supp(\mu^*)\}.$

Finally, the minimizing measures arising in
\eqref{eq:calk1}, \eqref{eq:calk2}, \eqref{eq:calk3}, and \eqref{eq:calk4}
are precisely
$\mu_r^{\mathrm{as}}, 
\mu_{s^*}^p, 
\mu_{t^*}^q, 
\mu_{s^*,t^*},$
respectively. For each of these measures, the proof of Proposition \ref{pr:saq-upper} provides the explicit values of
$m_p(\mu^*)$ and $m_q(\mu^*)$, leading to
$$
-\frac{1}{2p}\log m_p(\mu^*)
+\frac12
\log
\min\left\{
r\frac{m_p(\mu^*)^{1/p}}{m_q(\mu^*)^{1/q}},
1
\right\}=0
$$
in each case. Therefore,
\[
\liminf_{n\to\infty}
\frac{1}{\beta n^2}
\log
\frac{\vol\left(\mathbb B_{p,\beta}^n\cap r n^{\frac1q-\frac1p}\mathbb B_{q,\beta}^n\right)}
{\vol\left(\mathbb B_{p,\beta}^n\right)}
\ge
-\mathcal I_{p,q}(r).\qedhere
\]
\end{proof}

Let us turn to the general case.  It requires constrained equilibrium measures on $\rz_+$. 

\begin{lemma}\label{lemhalfp}
Let $s, p>0$ and consider the external potential
$V(x)=sx^{p/2}$ for $x\ge 0.$
For $a>0$, let
$$
\mu_a=
\argmin_{\mu\in\mx_1([0,a])}
\left\{-\frac12\Sigma(\mu)+\int V(x)\mu(\dd x)\right\}.
$$

\begin{enumerate}
 \item[(i)] If $0<a\le a_0:=(\frac{b_p^p}{sp})^{2/p},$ then $\supp(\mu_a)=[0,a]$ and $\mu_a$ has density
$$
\si_a(x)
=
\frac{1}{\pi\sqrt{x(a-x)}}
\left(1-\frac{sp}{b_p^p}a^{p/2}\right)
+\frac{sp^2}{2\pi b_p^p}a^{p/2}x^{\frac{p-2}{2}}
\int_x^a \frac{y^{-\frac{p+1}{2}}}{\sqrt{a-y}}\dd y,
\quad x\in(0,a).
$$
Moreover, for every $q>0$,
$$
m_{q/2}(\mu_a)
=
\frac{2a^{q/2}}{b_q^q}
\left(1-\frac{spq}{q+p}\frac{a^{p/2}}{b_p^p}\right),
\qquad \text{and} \qquad
\Sigma(\mu_a)
=
\log\frac{a}{4}
-s^2p \frac{a^{p}}{b_p^{2p}}.
$$
\item[(ii)] If $a>a_0$, then $\mu_a=\mu_{a_0}$. In particular,
$$
m_{q/2}(\mu_{a_0})
=
\frac{2p}{q+p}
\frac{1}{b_q^q}
\left(\frac{b_p^p}{sp}\right)^{q/p},
\quad \text{ and }\quad
\Sigma(\mu_{a_0})
=
2\log\frac{b_p}{2}-\frac{2}{p}\log(sp)
-\frac{1}{p}.
$$
\end{enumerate}
\end{lemma}

\begin{remark}\label{rem:squlm}
Let $0<p<\infty$. Taking $s=1/p$ in Lemma \ref{lemhalfp} shows that the measure
$$
\eta^{(p)}
:=\argmin_{\mu\in\mathcal M_1([0,\infty))}
\left\{
\frac1p\int_0^\infty x^{p/2} \mu(\dd x)
-\frac12\Sigma(\mu)
\right\}
$$
has density
\begin{equation}\label{eq:peta}
\eta^{(p)}(\dd x)
=\frac1{b_p\sqrt x}h_p\left(\frac{\sqrt x}{b_p}\right)
\mathbf1_{(0,b_p^2)}(x)\dd x.
\end{equation}
Here $h_p$ is the Ullman density in \eqref{eq:ullman}.
Equivalently, \eqref{eq:peta} gives $\eta^{(p)}=(x\mapsto x^2)_\#\nu^{(p)}$. In particular, for $0<q<\infty$,
$$
m_{p/2}(\eta^{(p)})=1,
\qquad
m_{q/2}(\eta^{(p)})=r_{p,q}^q,
\qquad
m_{p/2}(\eta^{(q)})^{-1}=s_{p,q}^p,
$$
where $r_{p,q}$ and $s_{p,q}$ are defined in \eqref{eq:rpq}.
This density agrees with \cite{HP00}*{Proposition 5.3.5}.
Moreover, $\eta^{(p)}$ maximizes $\Sigma(\mu)$ over
$\mu\in\mathcal M_1([0,\infty))$ with $m_{p/2}(\mu)\le1$.
\end{remark}

\begin{lemma}\label{le:halfst}
Let $s,t,p,q>0$ and $V(x)=sx^{p/2}+tx^{q/2}$ for $x\ge0$. Define
$$
\mu_{s,t}=
\argmin_{\mu\in\mathcal M_1(\mathbb R_+)}
\left\{-\frac12\Sigma(\mu)+\int V\dd \mu \right\}.
$$
Then $\mu_{s,t}$ is supported on $[0,a^2]$, where $a>0$ is uniquely determined by $$\frac{ps}{b_p^p}a^{p}+\frac{qt}{b_q^q}a^{q}=1.$$ Moreover, $\mu_{s,t}$ admits a density given by
\begin{align*}
\sigma_{s,t}(x)
= \frac{sp^2}{2\pi b_p^p}a^{p}
x^{\frac{p-2}{2}}
\int_x^{a^2} \frac{y^{-\frac{p+1}{2}}}{\sqrt{a^2-y}}\dd y
+\frac{tq^2}{2\pi b_q^q}a^{q}
x^{\frac{q-2}{2}}
\int_x^{a^2} \frac{y^{-\frac{q+1}{2}}}{\sqrt{a^2-y}}\dd y.
\end{align*}
For every $r>0$, the absolute moments satisfy
$$
m_{r/2}(\mu_{s,t})
=
\frac{2sp^2}{p+r}
\frac{a^{p+r}}{b_p^p b_r^r}
+
\frac{2tq^2}{q+r}
\frac{a^{q+r}}{b_q^q b_r^r},
$$
and the noncommutative entropy is
$$
\Sigma(\mu_{s,t})
=
2\log\frac{a}{2}
-\frac{ps^2}{b_p^{2p}}a^{2p}
-\frac{qt^2}{b_q^{2q}}a^{2q}
-\frac{4pq}{p+q}
\frac{st}{b_p^p b_q^q}
a^{p+q}.
$$
\end{lemma}

The following external potential arises in the probabilistic representation of singular values of random matrices. 
We consider
$$
V_n(x)=\frac{x^{p/2}}p
-\left(\frac1{2n}-\frac1{\beta n}\right)\log x,
\qquad x>0.
$$
These are variants of Freud weights.  Note that  $V_n$ is not bounded below for $\bt<2$ and thus the classic theory (see, e.g., \cite{HP00}*{Theorem~5.5.1}, \cite{AGZ10}*{Theorem 2.6.1}) does not apply.

\begin{lemma} \label{x:powemp}
Fix $\beta>0$ and $0<p<\infty$.  Let $(\lambda_1,\ldots,\lambda_n)$ have joint density
\begin{equation}\label{e:powpot}
\frac{1}{\bar Z_{n,p,\beta}^{(n)}}e^{-\beta n\sum_{i=1}^n V_n(\lambda_i)}
\prod_{1\leq i<j\leq n}|\lambda_i-\lambda_j|^\beta
\qquad\text{on }(0,\infty)^n,
\end{equation}
where $\bar Z_{n,p,\beta}^{(n)}$ is the normalizing constant. Then the empirical measure $L_n=n^{-1}\sum_{i=1}^n\delta_{\lambda_i}$ satisfies an LDP on $\mathcal M_1(\mathbb R_+)$ with speed $\beta n^2$ and good rate function $I_{W_p}$ given by \eqref{defiv}, with $W_p(x)=x^{p/2}/p$ for $x\geq0$. In particular,
\begin{equation}\label{eq:powZ}
\lim_{n\to\infty}\frac{1}{\beta n^2}\log\bar Z_{n,p,\beta}^{(n)}=C_{W_p}
=\log\frac{b_p}{2}-\frac{3}{2p}.
\end{equation}
\end{lemma}

\begin{proof}[Proof of Proposition \ref{pr:nspqup}]
Let $Z_n$ be uniformly distributed on $\mathbb{SB}_{p,\beta}^n$. Then
$$
\mathbb P\left(Z_n\in r n^{\frac1q-\frac1p} \mathbb{SB}_{q,\beta}^n\right)
=
\frac{\vol\left(\mathbb{SB}_{p,\beta}^n\cap r n^{\frac1q-\frac1p} \mathbb{SB}_{q,\beta}^n\right)}
{\vol\left(\mathbb{SB}_{p,\beta}^n\right)}.
$$
The proof of \cite{KPT20c}*{Theorem~1.5} uses the LDP of empirical measure in \cite{HP00}*{Theorem~5.5.1}.
The density \eqref{e:powpot} contains the factor $\prod_{i=1}^n\lambda_i^{\beta/2-1}$.
When $\beta=1$, the exponent $\beta/2-1$ violates the nonnegativity assumption in that LDP. Instead, using Lemma \ref{x:powemp} as the input for \cite{KPT20c}*{Theorem~1.5}, we still have an LDP for $\mu_n:=n^{-1}\sum_{j=1}^n\delta_{n^{2/p}s_j^2(Z_n)}$ on $\mathcal M_1(\mathbb R_+)$ with speed $\beta n^2$ and good rate function 
\begin{align}\label{eq:defmathy}
\mathcal{Y}(\mu)=-\frac{1}{2} \iint_{\mathbb{R}_{+}^2} \log |x-y| \mu(\dd x) \mu(\dd y)+2C_p,
\end{align}
for probability measures $\mu$ satisfying $m_{p/2}(\mu)\leq 1$, and $\mathcal{Y}(\mu)=+\infty$ otherwise. 
Moreover,
$$
\sum_{j=1}^n s_j(Z_n)^q
=
n^{1-q/p}m_{q/2}(\mu_n),
\quad 0<q<\infty.
$$
Hence
$$\limsup_{n \rightarrow \infty} \frac{1}{\beta n^2}
\log \mathbb P\left(
Z_n \in r n^{\frac1q-\frac1p}\mathbb{SB}_{q,\beta}^n
\right)
\leq - \inf_{\mu \in T_q(r)}\mathcal{Y}(\mu),$$
where
$$
T_q(r)=\left\{\mu \in \mathcal{M}_1(\mathbb{R}_+): m_{p/2}(\mu)\le1, 
m_{q/2}(\mu)\le r^q\right\},
\quad 0<q<\infty,
$$
and
$$
T_{\infty}(r)=\left\{\mu \in \mathcal{M}_1(\mathbb{R_+}): \supp(\mu)\subset[0,r^2],\
m_{p/2}(\mu)\le1\right\}.
$$
It remains to prove that
$
\inf_{\mu\in T_q(r)} \mathcal Y(\mu)=2\mathcal I_{p,q}(r)
$
as stated in Theorem \ref{th:ballint}. 

\noindent{\bf (i) The case $q=\infty$}. In this case, 
$$
\inf_{\mu\in T_\infty(r)} \mathcal Y(\mu) =
\inf_{\mu \in \mathcal M_1(\mathbb R_+) }\Big\{
-\tfrac12 \Sigma(\mu)+2C_p:
 \mathrm{supp}(\mu)\subset[0,r^2],\
m_{p/2}(\mu)\le1
\Big\}.
$$

\noindent\emph{(a) The case $0<r\le \frac{b_p}{2^{1/p}}.$} By \cite{ST97}*{Example I.3.5}, the minimizer of 
 $$\inf_{\mu \in \mathcal{M}_1(\mathbb{R_+})}\left\{-\frac{1}{2} \Sigma(\mu)+C_p: \supp(\mu) \subset[0, r^2]\right\}$$
 is the arcsine distribution
$
\nu_r^{\mathrm{as}}(\dd x)
=\frac{1}{\pi\sqrt{x(r^2-x)}}\mathbf1_{(0,r^2)}(x)\dd x.
$
A direct computation yields
$$
m_{p/2}(\nu_r^{\mathrm{as}})=\int x^{p/2}\dd \nu_r^{\mathrm{as}}
=
\frac{2r^p}{b_p^p} \leq 1,
\quad
\Sigma(\nu_r^{\mathrm{as}})
=\log\frac{r^2}{4}.
$$
Therefore, for $r\le b_p/2^{1/p}$,
\begin{align}\label{eq:caly1}
	\inf_{\mu\in T_\infty(r)} \mathcal Y(\mu) 
=
\mathcal Y(\nu_r^{\mathrm{as}}) 
=
-\tfrac12 \Sigma\left(\nu_r^{\text {as }}\right)+2C_p
=
\log\frac{b_p}{r} -\frac1{2p}.
\end{align}

\noindent\emph{(b) The case $b_p/2^{1/p}<r<b_p$.}
 In this case, $m_{p/2}(\nu_r^{\mathrm{as}}) >1.$ Arguing as before,
 $$
\inf_{\mu\in T_\infty(r)} \mathcal Y(\mu) =
\inf_{\mu \in \mathcal M_1(\mathbb R_+) }\Big\{
-\tfrac12 \Sigma(\mu)+2C_p:
 \mathrm{supp}(\mu)\subset[0,r^2],\
\int x^{p/2}\mu(\dd x) = 1
\Big\}.
$$
 Introducing a Lagrange multiplier $s\ge0$, consider
$$
\nu_s^p
=
\argmin_{\mathrm{supp}(\mu)\subset[0,r^2]}
\Big\{
-\tfrac12\Sigma(\mu)
+
s\Bigl(\int x^{p/2}\mu(\dd x)-1\Bigr)
\Big\}.
$$
By Lemma \ref{lemhalfp}(i), for $r^2<\bigl(\frac{b_p^p}{sp}\bigr)^{2/p}$, we have $
m_{p/2}(\nu_s^p) =\frac{2r^p}{b_p^p}
\left(1-\frac{sp}{2}\frac{r^p}{b_p^p}\right).
$ 
Imposing $m_{p/2}(\nu_s^p)=1$ gives $
s^*
=
\frac{2b_p^p}{pr^p}
\left(1-\frac{b_p^p}{2r^p}\right).
$
For $b_p/2^{1/p}<r<b_p$, we have $r^2<\left(\frac{b_p^p}{s^* p}\right)^{2 / p}$. Furthermore, Lemma \ref{lemhalfp}(i) gives
$\Sigma(\nu_{s^*}^p)
= \log \frac{r^2}{4}
-\frac{4}{p}\left(1-\frac{b_p^p}{2 r^p}\right)^2.
$
Hence,
\begin{align}\label{eq:caly2}
\inf_{\mu\in T_\infty(r)} \mathcal Y(\mu) 
=
\mathcal Y(\nu_{s^*}^p) 
=
-\tfrac12\Sigma(\nu_{s^*}^p)+2C_p
=
\log\frac{b_p}{r}
+
\frac{2}{p}
\left(
\frac{b_p^p}{2r^p}-1
\right)^2
-
\frac1{2p}.
\end{align}

\noindent{\bf (ii) The case $q<\infty$.} 
Let $\eta^{(p)}$ be the probability measure defined in \eqref{eq:peta}. By Remark \ref{rem:squlm}, $\eta^{(p)}$ maximizes $\Sigma$ under $m_{p/2}\le1$ and 
$m_{p/2}(\eta^{(p)}) =1, m_{q/2}(\eta^{(p)}) =r_{p,q}^q>r^q.$
By the same argument as in the proof of Proposition \ref{pr:saq-upper} part (ii),
$$\inf_{\mu\in T_q(r)} \mathcal Y(\mu)=
\inf_{\mu \in \mathcal M_1(\mathbb R_+) }
\Bigl\{
-\tfrac12\Sigma(\mu)+2C_p:
\int x^{p/2}\dd \mu\le1,
\int x^{q/2}\dd \mu = r^q
\Bigr\}. $$

\noindent\emph{(a) The case $0<r\le s_{p,q}.$} Consider
$$
\nu_t^q= \argmin_{ \mu\in \mathcal M_1(\mathbb R_+) } \Bigl\{
-\tfrac12\Sigma(\mu) +t \Big(\int x^{q/2}\dd \mu-r^q\Big)+2C_p
\Bigr\}.
$$
By Lemma \ref{lemhalfp}(ii), the minimizer $\nu_t^q$ is supported on
$\left[0,\left(\frac{b_q^q}{tq}\right)^{2/q}\right],$
and satisfies
$$m_{p / 2}(\nu_{t}^q)=\frac{2 q}{(p+q) b_p^p}\left(\frac{b_q^q}{t q}\right)^{p / q}, \quad
\Sigma(\nu_t^q)
=
2\log\frac{b_q}{2}-\frac{2}{q}\log(tq)
-\frac{1}{q}.$$
Imposing $m_{q/2}(\nu_t^q)=r^q$ yields $t^*=1/(q r^q)$. A direct computation gives
$m_{p/2}(\nu_{t^*}^q) =\frac{2qb_q^p}{(p+q)b_p^p}r^p.$
Thus, when $r^p\le \frac{p+q}{2q}\frac{b_p^p}{b_q^p}
=s_{p,q}^p,$ we have $m_{p/2}(\nu_{t^*}^q)\leq 1$ and $m_{q/2}(\nu_{t^*}^q) =r^q.$ 
Therefore,
\begin{align}\label{eq:caly3}
\inf_{\mu\in T_q(r)} \mathcal Y(\mu)=\mathcal Y(\nu_{t^*}^q)
=-\frac{1}{2} \Sigma\left(\nu_{t^*}^q\right)+2 C_p
=\log\frac{b_p}{b_qr}
+\frac{1}{2q}
-\frac{1}{2p}.
\end{align}

\noindent\emph{(b) The case $s_{p,q}<r<r_{p,q}$.} In this case, $ m_{p / 2}\left(\nu_{t^*}^q\right) > 1$, so the $p$-moment constraint becomes active and
$$\inf_{\mu\in T_q(r)} \mathcal Y(\mu)
=
\inf_{\mu \in \mathcal M_1(\mathbb R_+) }\Big\{
-\tfrac12\Sigma(\mu)+2C_p:
\int x^{p/2}\mu(\dd x)=1,
\int x^{q/2}\mu(\dd x)=r^q
\Big\}.$$
Introducing Lagrange multipliers $s,t\ge0$, let
$$\nu_{s,t}= \argmin_{\mu \in \mathcal M_1(\mathbb R_+) }\Big\{ -\tfrac12\Sigma(\mu)
+s\int x^{p/2}\mu(\dd x)+t\int x^{q/2}\mu(\dd x)\Big\}.$$
By Lemma \ref{le:halfst}, the minimizer $\nu_{s,t}$ is supported on $[0,a^2]$, where $a>0$ satisfies 
\begin{align}\label{conda}
	\frac{p s}{b_p^p}a^{p}
+\frac{q t}{b_q^q}a^{q}
=1.
\end{align}
Moreover,
$$m_{p/2}(\nu_{s,t})
=sp\frac{a^{2p}}{b_p^{2p}}
+
\frac{2tq^2}{p+q}
\frac{a^{p+q}}{b_p^p b_q^q}, \quad
m_{q/2}(\nu_{s,t})
=
\frac{2sp^2}{p+q}
\frac{a^{p+q}}{b_p^p b_q^q}
+
tq\frac{a^{2q}}{b_q^{2q}}.$$
Solving the system $m_{p/2}(\nu_{s,t})=1,
m_{q/2}(\nu_{s,t})=r^q$ 
yields 
$$s^*
=
\frac{b_p^p(p+q)}{p a^p(p-q)^2}
\left(
(p+q)\frac{b_p^p}{a^p}
-
2q\frac{b_q^q}{a^q}r^q
\right), \quad
t^*
=
\frac{b_q^q(p+q)}{q a^q(p-q)^2}
\left(
(p+q)\frac{b_q^q}{a^q}r^q
-
2p\frac{b_p^p}{a^p}
\right).$$
As in the proof of Proposition \ref{pr:saq-upper}(ii)(b), the nonnegativity constraints $s^*,t^*\ge0$, together with \eqref{conda} and Lemma \ref{le:apq}, imply that $a$ must satisfy \eqref{eq:a-def}. We denote by $a_{p,q}(r)$ the smaller positive
solution to this equation. As in the proof of Proposition \ref{pr:saq-upper},
the constraints $s^*,t^*\ge0$ force this choice of $a$.
Plugging $(s^*,t^*)$ into
$\Sigma(\nu_{s,t})$ as given in Lemma \ref{le:halfst}, we obtain
$$
\Sigma(\nu_{s^*,t^*})
=
2\log\frac{a_{p,q}(r)}{2}
-
\frac{1}{p+q}
-
\left(\frac1p+\frac1q\right)
\left(
\frac{b_p^p}{a_{p,q}(r)^p}
-
\frac{p+2q}{p+q}
\right)^2.
$$
Finally,
\begin{align}\label{eq:caly4}
\inf_{\mu \in T_q(r)} \mathcal{Y}(\mu)= \mathcal{Y}(\nu_{s^*,t^*})
= 
-\frac12\Sigma(\nu_{s^*,t^*})
+2C_p.
\end{align}
After simplification, this gives precisely the rate function $2 \mathcal I_{p,q}(r)$ stated in Theorem \ref{th:ballint}.
\end{proof}

\begin{proof}[Proof of Proposition \ref{pr:nspqlw}] 
The argument for the lower bound is similar to that of Proposition \ref{pr:saq-lower}. We start with the stochastic representation from \cite{KPT20c}*{Proposition 8.2} for the unordered singular values of $Z_n$:
\begin{align}\label{eq:storex}
 \left(s_1^2\left(Z_n\right), \ldots, s_n^2\left(Z_n\right)\right) \stackrel{d}{=} U^{\frac{1}{m_{\beta, n }}} \frac{\bX}{\left\|\bX\right\|_{p / 2}}
\end{align}
where $m_{\beta, n} = \frac{\beta n^2}{2}$ and the random vector $\bX=(X_1,...,X_n) \in \mathbb{R}_+^n$ has joint density 
\begin{align*}
\frac{1}{Z_{n,p,\beta}} \prod_{1 \leq i<j \leq n}\left|x_i-x_j\right|^\beta \prod_{i=1}^n x_i^{\frac{\beta}{2}-1} e^{-n \sum_{i=1}^n\left|x_i\right|^{p / 2}} 
\end{align*}
and, independently of $\bX$, $U$ is uniformly distributed on $(0,1)$. 
Applying the change of variables $x_i=(\beta/p)^{2/p} y_i, i=1,\dots,n,$ we obtain
\begin{align*}
 &\mathbb P\left(Z_n\in r_n n^{1/q-1/p}\mathbb{SB}_{q,\beta}^n \right) \\
&= 
\frac{1}{\bar Z_{n,p,\beta}^{(n)}}
\int_{\mathbb R_+^n}
\min\Big\{
\Big(
 \frac{r_n^2 n^{-2/p} \|\by\|_{p/2}}{n^{-2/q}\|\by\|_{q/2}}
\Big)^{m_{\beta,n}},
1
\Big\} \prod_{i<j}|y_i-y_j|^\beta
\prod_{i=1}^n y_i^{\frac{\beta}{2}-1}
e^{-\beta n\sum_{i=1}^n y_i^{p/2}/p}
\dd \by,
\end{align*}
where
$$\bar Z_{n,p,\beta}^{(n)}=(p/\beta)^{\beta n^2/p} Z_{n,p,\beta}
=
\int_{\mathbb R_+^n}
\prod_{i<j}|y_i-y_j|^\beta
\prod_{i=1}^n y_i^{\frac{\beta}{2}-1}
e^{-\beta n\sum_{i=1}^n y_i^{p/2}/p}
\dd\by,$$
and we interpret the formula as before when $q=\8$.

Assume that $\mu^*$ is one of the minimizing measures appearing in \eqref{eq:caly1}, \eqref{eq:caly2}, \eqref{eq:caly3}, and \eqref{eq:caly4}, namely,
$
\nu_r^{\mathrm{as}},
\nu_{s^*}^p,
\nu_{t^*}^q,
\nu_{s^*,t^*}.
$
Define
$c=m_{p/2}(\mu^*)^{2/p},
\nu^*(A)=\mu^*(cA).$
Then
$$
m_{q/2}(\nu^*)
=
\frac{m_{q/2}(\mu^*)}{m_{p/2}(\mu^*)^{q/p}},
\quad
\Sigma(\nu^*)
=
\Sigma(\mu^*)
-
\frac{2}{p}\log p
-
\frac{2}{p}\log m_{p/2}(\mu^*).$$
Arguing as in the proof of the lower bound in Proposition \ref{pr:saq-lower}, we construct empirical configurations associated with $\nu^*$ and obtain
\begin{align*}
 & \liminf_{n\to\infty}\frac{1}{\beta n^2}
\log \mathbb{P}\left(
Z_n\in r_n n^{1/q-1/p}\mathbb{SB}_{q,\beta}^n
\right) \\
\ge{}&
-
\lim_{n\to\infty}
\frac{1}{\beta n^2}
\log \bar Z_{n,p,\beta}^{(n)}
+
\frac{1}{2}\Sigma(\nu^*)
+
\frac{1}{2}
\log
\min\left\{
\frac{
\left(m_{p/2}(\nu^*)\right)^{2/p} r^2
}{
\left(m_{q/2}(\nu^*)\right)^{2/q}
},
1
\right\}
-
\frac{1}{p}m_{p/2}(\nu^*).
\end{align*}
Plugging in \eqref{eq:powZ} and using the scaling relations between $\nu^*$ and $\mu^*$, we obtain
\begin{align}
 & \liminf_{n\to\infty}\frac{1}{\beta n^2}
\log \mathbb{P}\left(
Z_n\in r_n n^{1/q-1/p}\mathbb{SB}_{q,\beta}^n
\right) \notag \\
\ge{} &
\frac{1}{2}\Sigma(\mu^*)
+
\log\frac{2}{b_p}+
\frac{1}{2p}
-
\frac{1}{p}\log m_{p/2}(\mu^*)
+
\frac{1}{2}
\log
\min\left\{
r^2
\frac{
\left(m_{p/2}(\mu^*)\right)^{2/p}
}{
\left(m_{q/2}(\mu^*)\right)^{2/q}
},
1
\right\}. \label{lbsbq}
\end{align}
In the case $q=\infty$, the term $\left(m_{q / 2}\left(\mu^*\right)\right)^{2 / q}$ is replaced by $\sup \{x: x \in \left.\supp\left(\mu^*\right)\right\}$.
For each of these measures, the proof of Proposition \ref{pr:nspqup} provides the explicit values of
$m_{p/2}(\mu^*)$, $m_{q/2}(\mu^*)$, and $\Sigma(\mu^*)$. Substituting these expressions into the right hand side of \eqref{lbsbq} reveals that the last two summands cancel out, yielding that
$$
\liminf_{n\to \infty}
\frac{1}{\beta n^2}
\log
\frac{\vol\left(\mathbb{SB}_{p,\beta}^n \cap r n^{\frac1q-\frac1p}\mathbb{SB}_{q,\beta}^n\right)}
{\vol\left(\mathbb{SB}_{p,\beta}^n\right)}
\ge
-2\mathcal I_{p,q}(r).
$$
This completes the proof.
\end{proof}

\subsection{Proofs of auxiliary lemmas}\label{s:intlem}
We state the following elementary inequalities without proof.
\begin{lemma}
For $z,w\in\mathbb C$,
\begin{equation}\label{e:pairbd}
 |z-w|\leq(1+|z|)(1+|w|), \qquad
 \log_+|z-w|\leq\log(1+|z|)+\log(1+|w|),
\end{equation}
where $\log_+ x =\log x \mathbf1_{\{x\ge1\} }$. For every $\alpha,\varepsilon>0$, there is $C>0$ such that
\begin{align}
 \log(1+t)&\leq\varepsilon t^\alpha+C,
 &t\geq0,\label{eq:Frelb}\\
 \log|z-w|&\leq\varepsilon(|z|^\alpha+|w|^\alpha)+C,
 &z,w\in\mathbb C.\label{e:logpair}
\end{align}
\end{lemma}

\begin{proof}[Proof of Lemma \ref{le:apq}]
(i) For fixed $r>0$, let
$$
F_r(a)
=
\frac{b_q^q}{a^q}r^q-\frac{b_p^p}{a^p}-\frac{p-q}{p+q},
\qquad a>0.
$$
Then \eqref{eq:a-def} is equivalent to $F_r(a)=0$. 
The derivative $F_r'(a)=\frac{\partial F_r(a)}{\partial a}$ has exactly one positive zero. Hence $F_r$ has at
most two positive zeros. For $s_{p,q}<r<r_{p,q}$, we have
$$
F_r(b_p)
=
\left(\frac{b_qr}{b_p}\right)^q-\frac{2p}{p+q}<0, \quad 
F_r(b_qs_{p,q})
=
\left(\frac{r}{s_{p,q}}\right)^q-1>0.
$$
If $q<p$, then $b_p<b_qs_{p,q}$. Since
$$
F_r(b_p)<0,\qquad
F_r(b_qs_{p,q})>0,
\qquad
\lim_{a\to\infty}F_r(a)
=
-\frac{p-q}{p+q}<0,
$$
there is one positive zero in $(b_p,b_qs_{p,q})$ and another one in
$(b_qs_{p,q},\infty)$. Since $F_r$ has at most two positive zeros, these
are the only two positive zeros. In particular, the smaller positive zero lies
in $(b_p,b_qs_{p,q})$.
If $p<q$, then $b_qs_{p,q}<b_p$. Since
$$
F_r(b_qs_{p,q})>0,\qquad
F_r(b_p)<0,
\qquad
\lim_{a\to\infty}F_r(a)
=
\frac{q-p}{p+q}>0,
$$
there is one positive zero in $(b_qs_{p,q},b_p)$ and another one in $(b_p,\infty)$. Again, these
are the only two positive zeros. Hence the smaller positive zero lies in
$(b_qs_{p,q},b_p)$.

(ii) Write $a_0=b_qs_{p,q}$. By definition, $F_r(a_0)>0, r>s_{p,q}$ and $F_{s_{p,q}}(a_0)=0.$ 
Assume first that $q<p$. Then $F_{s_{p,q}}'(a_0)>0$. Hence, for every sufficiently small $\varepsilon>0$, $F_{s_{p,q}}(a_0-\varepsilon)<0.$
Since $r\mapsto F_r(a_0-\varepsilon)$ is continuous and $F_{s_{p,q}}(a_0-\varepsilon)<0$, there exists $\delta>0$ such that
$$
F_r(a_0-\varepsilon)<0<F_r(a_0)
=\left(\frac{r}{s_{p,q}}\right)^q-1,
\qquad s_{p,q}<r<s_{p,q}+\delta.
$$
Thus $F_r$ has a zero in $(a_0-\varepsilon,a_0)$. Since the smaller
positive zero has already been shown to lie in $(b_p,b_qs_{p,q})=(b_p,a_0)$,
we get $a_0-\varepsilon<a_{p,q}(r)<a_0$ for all $r>s_{p,q}$ sufficiently close to $s_{p,q}$. Therefore $\lim_{r\downarrow s_{p,q}}a_{p,q}(r)=a_0.$

Assume next that $p<q$. Then $F_{s_{p,q}}'(a_0)<0$ and $F_{s_{p,q}}(a_0+\varepsilon)<0$ for every sufficiently small $\varepsilon>0$.  Again by continuity, for all $r>s_{p,q}$ close enough to $s_{p,q}$, $F_r(a_0)>0>F_r(a_0+\varepsilon).$
Thus $F_r$ has a zero in $(a_0,a_0+\varepsilon)$. Since the smaller
positive zero lies in
$(b_qs_{p,q},b_p)=(a_0,b_p)$, we get
$
a_0<a_{p,q}(r)<a_0+\varepsilon
$
for all $r>s_{p,q}$ sufficiently close to $s_{p,q}$. Hence $\lim_{r\downarrow s_{p,q}}a_{p,q}(r)=a_0.$


(iii) At $r=r_{p,q}$,
$$
F_{r_{p,q}}(b_p)=0,\qquad
F'_{r_{p,q}}(b_p)
=\left.\frac{\partial F_r(a)}{\partial a}\right|_{r=r_{p,q},\,a=b_p}
=\frac{p(p-q)}{(p+q)b_p}\neq0.
$$
Thus $F_{r_{p,q}}(a)$ has the sign of $p-q$ for $a>b_p$ sufficiently close to $b_p$, whereas $$\lim_{a\to+\infty}F_{r_{p,q}}(a)=-(p-q)/(p+q)$$ has the opposite sign.
The intermediate value theorem gives another root larger than $b_p$.
Since there are at most two positive roots, $b_p$ is the smaller root, and the implicit function theorem shows that $\lim_{r\uparrow r_{p,q}} a_{p,q}(r) = b_p$.

We now study the properties of $\mathcal I_{p,q}$ for $0<p<\infty$ and $0<q\leq\infty$. If $q<\infty$ and $q\ne p$, from part (iii) and \eqref{eq:ipqdef} we deduce
$$
\lim_{r\uparrow r_{p,q}}\mathcal I_{p,q}(r)
=\mathcal I_{p,q}(r_{p,q})
=\frac{q}{4p(p+q)}-\frac1{4p}+\frac1{4(p+q)}=0.
$$
Since $\lim_{r\to s_{p,q}^+}a_{p,q}(r) = a_0$, sending $r\downarrow s_{p,q}$ in the second branch of
\eqref{eq:ipqdef} yields
\begin{align*}
\lim_{r\downarrow s_{p,q}}\mathcal I_{p,q}(r)
=
\frac12\log\frac{b_p}{b_q s_{p,q}}
+
\frac14\left(\frac1q-\frac1p\right).
\end{align*}
This is exactly the value of the first branch at $r=s_{p,q}$. Hence,
$\mathcal I_{p,q}$ is continuous at $s_{p,q}$.

For $s_{p,q}<r<r_{p,q}$, by part (i) the two positive roots of $F_r(a)$ are simple. Using the implicit function theorem, $a_{p,q}$ is continuously differentiable on $(s_{p,q},r_{p,q})$.
Differentiating \eqref{eq:ipqdef}, since $a_{p,q}(r)>b_p$ for $q<p$ and $a_{p,q}(r)<b_p$ for $p<q$, we find
$$
\mathcal I'_{p,q}(r)
=\frac{(p+q)b_q^q r^{q-1}}{2(p-q)a_{p,q}(r)^q}
\left(\frac{b_p^p}{a_{p,q}(r)^p}-1\right)<0.
$$
On the first branch $r<s_{p,q}$, $\mathcal I'_{p,q}(r)=-1/(2r)<0$.
Together with the continuity at $s_{p,q}$ and $r_{p,q}$,
we conclude that $\ix_{p,q}(r)$ is strictly decreasing on $(0,r_{p,q}]$.

For $p=q$, the stated properties of $\mathcal I_{p,p}$ follow from the explicit formula $
\mathcal I_{p,p}(r)=-\frac12\log r,$  $0<r\leq r_{p,p} = 1.$ When $q=\infty$, $\ix_{p,\8}$ has a closed form \eqref{eq:ipinf}, which is continuous on $(0,b_p]$, $\mathcal I_{p,\infty}(b_p)=0$, and the derivative
$$
\mathcal I'_{p,\infty}(r)
=\begin{cases}
-\frac1{2r},&0<r<2^{-1/p}b_p,\\
-\frac1{2r}\left(\frac{b_p^p}{r^p}-1\right)^2,
&2^{-1/p}b_p<r<b_p.
\end{cases}
$$
Finally, $\ix_{p,q}\ge0$ follows from strict monotonicity and
$\mathcal I_{p,q}(r_{p,q})=0$ in all cases. 
\end{proof}

\begin{proof}[Proof of Lemma \ref{lemipq}]
By Lemma \ref{remr}, the minimizer $\mu_a$ is one-cut and
the natural endpoint of the unconstrained equilibrium problem is determined by the standard endpoint condition \eqref{eqb2}, which reduces to 
$\frac{2spa^p}{b_p^p}=1$ here.
Thus, $a_0=\frac{b_p}{(2sp)^{1/p}}$ and $\supp(\mu_a)=[-a,a]$ for $a\le a_0$. In this case, Lemma \ref{remr} suggests that we should consider, after rescaling $x=au$, a density of the form
$$
\sigma_a(a u)=\frac{1}{\pi a \sqrt{1-u^2}}-\frac{s p a^{p-1}}{\pi^2} L_p(u), \quad u \in[-1,1],
$$
where $$L_p(u)=\pv \int_{-1}^1 \frac{\sqrt{1-v^2}}{\sqrt{1-u^2}} \frac{v|v|^{p-2}}{u-v} \dd v.$$
To evaluate this singular integral, we invoke a result of \cite{Rak84}, which gives for $u \in[0,1],$
$$K_p(u):=\pv \int_{-1}^1 \frac{\sqrt{1-u^2}}{\sqrt{1-v^2}} \frac{v|v|^{p-2}}{u-v} \dd v=-\frac{2 \pi p u^{p-1}}{b_p^p} \int_u^1 \frac{\dd v}{v^{p} \sqrt{1-v^2}}.$$
Moreover,
\begin{align*}
& L_p(u)-K_p(u) = \pv \int_{-1}^1\left(\frac{\sqrt{1-v^2}}{\sqrt{1-u^2}}-\frac{\sqrt{1-u^2}}{\sqrt{1-v^2}}\right) \frac{v|v|^{p-2}}{u-v} \dd v \\
 ={}& \frac{1}{\sqrt{1-u^2}} \int_{-1}^1 \frac{(u+v) v|v|^{p-2}}{\sqrt{1-v^2}} \dd v
=\frac{1}{\sqrt{1-u^2}} \int_{-1}^1 \frac{|v|^p}{\sqrt{1-v^2}} \dd v =\frac{2\pi}{b_p^p}\frac{1}{\sqrt{1-u^2}}.
\end{align*}
Combining the above identities yields for $u\in[0,1]$,
\begin{align*}
\sigma_a(au)
&=
\frac{1}{\pi a\sqrt{1-u^2}}
\left(1-\frac{2spa^p}{b_p^p}\right)
+
\frac{2sp^2a^{p-1}}{\pi b_p^p}
u^{p-1}
\int_u^1 v^{-p}\frac{\dd v}{\sqrt{1-v^2}}.
\end{align*}
Since $V$ is even, $\sigma_a$ is symmetric; returning to $x=au$ yields the stated density on $[-a,a]$. Note that the second summand is proportional to the Ullman density after the change of variable $v=au/s$:
$$
\frac{2sp^2a^{p-1}}{\pi b_p^p}
u^{p-1}
\int_u^1 v^{-p}\frac{\dd v}{\sqrt{1-v^2}}= \frac{2sp a^p}{b_p^p} \frac{p}{\pi a^p}\int_{au}^a \frac{s^{p-1} \dd s}{\sqrt{s^2-a^2u^2}}.
$$

For $a>a_0$ the constraint is inactive and $\mu_a=\mu_{a_0}$, which is exactly the case when the first summand of $\si_a$ vanishes. Note that $|x|^p$ is not differentiable at $x=0$ for $0<p<1$ so that \eqref{eq:densd} may not be directly applicable. However, the density $\si_a$ yields
\begin{align*}
 -\int_{-a}^a \log|x-y|\si_a(\dd y )+ V(x)\ge - \log\frac{a}{2}+\frac{2sa^p}{b_p^p}=:F_a
\end{align*}
with equality when $|x|\le a\le a_0$, which verifies the condition in Lemma \ref{lemr}; see \cite{ST97}*{Theorem IV.5.1} for the computations and for the case $a>a_0$.

For the $q$-th absolute moment, we write
$$m_q(\mu_a)=\int_{-a}^a |x|^q\mu_a(\dd x)=2\int_0^a x^q\sigma_a(x)\dd x =:M_1+M_2,$$
where the first term is given by the $q$-th moment of the arcsine distribution
\begin{align*}
M_1=\frac{2}{\pi}\Bigl(1-\frac{2spa^p}{b_p^p}\Bigr)\int_0^a\frac{x^q}{\sqrt{a^2-x^2}}\dd x 
=\frac{2a^q}{b_q^q}\Bigl(1-\frac{2spa^p}{b_p^p}\Bigr),
\end{align*}
and by Fubini's theorem,
\begin{align*}
M_2
&=\frac{4sp^2a^p}{\pi b_p^p}\int_0^a x^{q+p-1}\Bigl(\int_x^a v^{-p}\frac{\dd v}{\sqrt{a^2-v^2}}\Bigr)\dd x
=\frac{4sp^2a^p}{\pi b_p^p}\int_0^a v^{-p}\frac{\dd v}{\sqrt{a^2-v^2}}\int_0^v x^{q+p-1}\dd x \\
&=\frac{4sp^2a^p}{\pi b_p^p(q+p)}\int_0^a\frac{v^q}{\sqrt{a^2-v^2}}\dd v
=\frac{4sp^2}{(q+p)b_p^p b_q^q}a^{p+q}.
\end{align*}
Adding $M_1$ and $M_2$ gives the stated formula for the $q$-th moment.
We finally conclude that
\[
\Sigma\left(\mu_a\right)
=-\Big[F_a-\int s|x|^p\mu_a(\dd x)\Big]=
\log \frac{a}{2}
-\frac{2s^2 p a^{2p}}{b_p^{2p}}. \qedhere
\]
\end{proof}

\begin{proof}[Proof of Lemma \ref{lemspq}]
By Lemma \ref{remr}, for $s,t>0$ and $p,q>0$ the minimizer $\mu_{s,t}$ is a one-cut measure supported on $[-a,a]$ 
which is uniquely determined by \eqref{eqb2}. A direct computation gives
$$
\frac{2sp}{b_p^p}a^p+\frac{2tq}{b_q^q}a^q=1.
$$
The left hand side is continuous and strictly increasing in $a>0$, and hence this equation has a unique solution.

For $x\in[-a,a]$, define
\begin{align*}
\sigma_{s,t}(x)
=
\frac{2s p^2|x|^{p-1}a^{p}}{\pi b_p^p}
\int_{|x|}^{a} \frac{\dd w}{w^{p}\sqrt{a^2-w^2}}
+
\frac{2t q^2|x|^{q-1}a^{q}}{\pi b_q^q}
\int_{|x|}^{a} \frac{\dd w}{w^{q}\sqrt{a^2-w^2}}.
\end{align*}
The same computation using Fubini's theorem as in Lemma \ref{lemipq} yields 
$$
\int_{-a}^a
\frac{2r^2 |x|^{r-1}a^r}{\pi b_r^r}
\left(
\int_{|x|}^{a}\frac{\dd w}{w^r\sqrt{a^2-w^2}}
\right)\dd x
=
\frac{2ra^r}{b_r^r},
\qquad r>0.
$$
It follows that
$$
\int_{-a}^a \sigma_{s,t}(x)\dd x
=
\frac{2spa^p}{b_p^p}
+
\frac{2tqa^q}{b_q^q}
=1.
$$
Thus $\sigma_{s,t}(x)\dd x$ is a probability measure supported on $[-a,a]$. Substituting the above density into the
logarithmic potential gives
$$
-\int_{-a}^a \log|x-y| \sigma_{s,t}(y)\dd y
+
s|x|^p+t|x|^q
\ge
-\log\frac a2
+
\frac{2sa^p}{b_p^p}
+
\frac{2ta^q}{b_q^q}=:F,
\qquad x\in\mathbb R,
$$
with equality for $|x|\le a$. This verifies the condition in Lemma
\ref{lemr}; see \cite{ST97}*{Theorem IV.5.1} for the corresponding
computations. Consequently,
$\sigma_{s,t}(x)\dd x=\mu_{s,t}(\dd x),
$
where
$$
\mu_{s,t}=
\argmin_{\mu\in\mathcal M_1(\mathbb R)}
\left\{-\frac12\Sigma(\mu)+\int V(x)\mu(\dd x)\right\}.
$$

The formulas for the $p$-th and $q$-th moments follow from the Fubini theorem as used in Lemma \ref{lemipq}. It remains to compute $\Sigma(\mu_{s,t})$. By \eqref{eq:calsig},
\begin{align*}
    \Sigma(\mu_{s,t})
={}&
\int V \dd\mu_{s,t}-F
=
s m_p(\mu_{s,t})
+
t m_q(\mu_{s,t})
+
\log\frac a2
-
\frac{2s a^p}{b_p^p}
-
\frac{2t a^q}{b_q^q}\\
={}& \log\frac{a}{2}
-\frac{2s^2pa^{2p}}{b_p^{2p}}
-\frac{2t^2qa^{2q}}{b_q^{2q}}
-\frac{8stpq}{p+q}
\frac{a^{p+q}}{b_p^p b_q^q},
\end{align*}
where, in the last step, we substituted the moment formula with $r=p$ and $r=q$ and used the endpoint equation. This completes the proof.
\end{proof}

\begin{proof}[Proof of Lemma \ref{lemhalfp}]
For $V(x)=sx^{p/2}$, we have
$
V'(x)=\frac{sp}{2}x^{p/2-1},x>0.
$
Although $V'$ is singular at $0$ when $0<p<2$, the weight
$w=e^{-V}$ is absolutely continuous on $[0,a]$ and $w'\in L^r([0,a])$
for some $r>1$. Hence we can still apply Lemma \ref{remr}. Since
$xV'(x)=\frac{sp}{2}x^{p/2}$ is increasing on $(0,\infty)$, the minimizer is supported in a bounded interval. In the  case $a\le a_0$, $\supp(\mu_a)=[0,a]$, and its density is
$$
\sigma_a(x)
=
\frac{\pi-\ell(x)}{\pi^2\sqrt{x(a-x)}}, \qquad
\ell(x)
=
\frac{sp}{2} \pv\int_0^a
\sqrt{(a-y)y} 
\frac{y^{p/2-1}}{x-y} \dd y,
\qquad x\in(0,a).
$$

The unconstrained right endpoint is determined by \eqref{eqb2} that simplifies to
$
\frac{sp}{b_p^p}a^{p/2}=1.
$
It follows that
$
a_0=\left(\frac{b_p^p}{sp}\right)^{2/p}.
$
Thus the constraint is active (or in effect) for $0<a\le a_0$, while for $a>a_0$ it is inactive and $\mu_a=\mu_{a_0}$.

For the explicit density, we rewrite $\ell(x)$ as
\begin{align*}
    \ell(x)= 
\frac{sp}{2}
\pv \int_{0}^{a}
\frac{(a-y)y^{\frac{p-1}{2}}}{(x-y)\sqrt{a-y}}\dd y
= \frac{sp}{2}
\bigg[
(a-x) \pv \int_0^a
\frac{y^{\frac{p-1}{2}}}{(x-y)\sqrt{a-y}}\dd y
+
\int_0^a
\frac{y^{\frac{p-1}{2}}}{\sqrt{a-y}}\dd y
\bigg].
\end{align*}
Using the principal value identity of \cite{Rak84}
$$\pv \int_0^1
\frac{t^{\frac{p-1}{2}}}{\sqrt{1-t}(t-z)}\dd t
=
\frac{\sqrt{\pi} p \Gamma\left(\frac{p+1}{2}\right)}{2\Gamma\left(\frac{p}{2}+1\right)}
\frac{z^{\frac{p-1}{2}}}{\sqrt{1-z}}
\int_z^1 t^{-\frac{p+1}{2}}\frac{\dd t}{\sqrt{1-t}},$$
we deduce
$$\ell(x)
=
-\frac{s p^{2}}{2}
\frac{\pi a^{p/2}}{b_p^{p}}
x^{\frac{p-1}{2}}\sqrt{a-x}
\int_x^a
\frac{y^{-\frac{p+1}{2}}}{\sqrt{a-y}}\dd y + \frac{s p\pi}{b_p^{p}}a^{p/2}.$$
The stated formula of $\si_a$ follows from here.

For the moment calculation, using the Beta integral
$$\int_0^a \frac{x^{q/2}}{\pi\sqrt{x(a-x)}}\dd x
=
\frac{1}{\pi}\int_0^a x^{\frac q2-\frac12}(a-x)^{-\frac12}\dd x
=
\frac{2a^{q/2}}{b_q^q},$$
together with Fubini's theorem,
\begin{align*}
 \int_0^a x^{\frac q2}x^{\frac p2-1}
\left(\int_x^a \frac{y^{-\frac{p+1}{2}}}{\sqrt{a-y}} \dd y\right) \dd x
=
\int_0^a
\left(\int_0^y x^{\frac{q+p-2}{2}}\dd x\right)
\frac{y^{-\frac{p+1}{2}}}{\sqrt{a-y}}\dd y =
\frac{4 \pi a^{q/2}}{(q+p)b_q^q},
\end{align*}
we obtain
$$m_{q/2}(\mu_a)
=
\frac{2a^{q/2}}{b_q^q}
\left(1-\frac{spq}{q+p}\frac{a^{p/2}}{b_p^p}\right).$$

For the logarithmic term, compute
$$
F_a
:=
\int_0^a \log x\sigma_a(x)\dd x =
\left(1-\frac{sp}{b_p^p}a^{p/2}\right)\log\frac{a}{4}
+
\frac{sp^2}{2b_p^p}a^{p/2}\cdot\frac{1}{\pi}J =
\log\frac{a}{4}-\frac{2s}{b_p^p}a^{p/2},
$$
where we used the identity
$$
J=\int_0^a x^{\frac{p-2}{2}} \log x
\left(\int_x^a \frac{y^{-\frac{p+1}{2}}}{\sqrt{a-y}}\dd y\right)\dd x
=
\pi\left(\frac{2}{p}\log\frac{a}{4}-\frac{4}{p^2}\right)$$
which in turn follows from the Fubini theorem. 
Applying Remark \ref{rem:qepw} to $V(x)=sx^{p/2}$ on
$[0,a]$ gives the equilibrium identity at $0\in\supp\mu_a$.
Since $V(0)=0$, the constant in the equilibrium condition \eqref{equcon} is $-F_a$.
By \eqref{eq:calsig}, we find
$$\Sigma(\mu_a)=sm_{p/2}(\mu_a)+F_a
=
\log\frac{a}{4}-\frac{s^2p}{b_p^{2p}}a^{p}.$$

If $a>a_0$, the constraint becomes inactive and thus $\mu_a=\mu_{a_0}$. Replacing $a$ by $a_0$ in the above expressions, we obtain the desired results.
\end{proof}

\begin{proof}[Proof of Lemma \ref{le:halfst}]
Let $A=a^2$. Since $xV'(x)$
is increasing on $(0,\infty)$, by Lemma \ref{remr} the equilibrium measure is 
supported on $[0,A]$. The endpoint $A=a^2$ is determined by \eqref{eqb}, i.e.,
$$
\frac{sp}{b_p^p}a^p+\frac{tq}{b_q^q}a^q=1,
$$
which has a unique positive solution as the left hand side is strictly increasing in $a$.
Note that the density formula \eqref{eq:densd} is linear in $V'$. We compute the principal values for the terms $sx^{p/2}$ and $tx^{q/2}$ as in the proof of Lemma \ref{lemhalfp} and obtain
$x\in(0,a^2)$,
$$
\sigma_{s,t}(x)
=
\frac{sp^2}{2\pi b_p^p}a^p
x^{\frac{p-2}{2}}
\int_x^{a^2}
\frac{y^{-\frac{p+1}{2}}}{\sqrt{a^2-y}} \dd y 
+
\frac{tq^2}{2\pi b_q^q}a^q
x^{\frac{q-2}{2}}
\int_x^{a^2}
\frac{y^{-\frac{q+1}{2}}}{\sqrt{a^2-y}} \dd y.
$$
Let us turn to the moments. For $r,u>0$, using the Fubini theorem as before, we find
$$
\int_0^{a^2}
x^{r/2}x^{\frac{u-2}{2}}
\left(
\int_x^{a^2}
\frac{y^{-\frac{u+1}{2}}}{\sqrt{a^2-y}} \dd y
\right)
\dd x
=
\frac{4\pi a^r}{(r+u)b_r^r}.
$$
It follows that for every $r>0$,
$$
m_{r/2}(\mu_{s,t})
=
\frac{2sp^2}{p+r}
\frac{a^{p+r}}{b_p^p b_r^r}
+
\frac{2tq^2}{q+r}
\frac{a^{q+r}}{b_q^q b_r^r}.
$$
Finally, we compute the noncommutative entropy.
Applying Remark \ref{rem:qepw} to $V(x)=sx^{p/2}+tx^{q/2}$ on
$[0,\infty)$ and evaluating \eqref{eq:vlogmu} at $x=0$ yield
\begin{align*}
F=V(0)-\int_0^{a^2}\log x\,\mu_{s,t}(\dd x)
=-2\log\frac{a}{2}
+\frac{2s}{b_p^p}a^p+\frac{2t}{b_q^q}a^q,
\end{align*}
where we used  the endpoint equation and the Fubini theorem again.
By \eqref{eq:calsig}, we have
\begin{align*}
	\Sigma(\mu_{s,t})
&=
\int V \dd\mu_{s,t}-F
=
s m_{p/2}(\mu_{s,t})
+
t m_{q/2}(\mu_{s,t})
-
F, \\
&=	2\log\frac{a}{2}
-
\frac{ps^2}{b_p^{2p}}a^{2p}
-
\frac{qt^2}{b_q^{2q}}a^{2q}
-
\frac{4pq}{p+q}
\frac{st}{b_p^p b_q^q}
a^{p+q}.
\end{align*}
This completes the proof.
\end{proof}

\begin{proof}[Proof of Lemma \ref{le:logsiri}]
Let
$$
Q(u)=F^{-1}(u),\quad u\in[0,1], \quad
S_n
=
\frac1{n^2}
\sum_{1\le i<j\le n}
\log(\ell_j-s_i).
$$
For $M>0$, define $\phi_M(x,y)=\log|x-y|\vee(-M).$
Let
$
S_{n,M}
=
n^{-2}
\sum_{1\le i<j\le n}
\phi_M(\ell_j,s_i).
$
Since $Q$ is continuous, the function $
g_M(u,v):=\phi_M(Q(u),Q(v)), (u,v)\in[0,1]^2,
$
is bounded and continuous. Therefore
\begin{equation}\label{eq:trriem}
S_{n,M}
=
\frac1{n^2}
\sum_{1\le i<j\le n}
g_M\left(\frac{j-\frac12}{n},\frac{i}{n}\right)
\to
\iint_D g_M(u,v) \dd u \dd v,
\end{equation}
where $D=\{(u,v)\in[0,1]^2:v<u\}.$ Since $g_M$ is symmetric and the set $\{(u,v):u=v\}$ has zero
two-dimensional Lebesgue measure, the integral over the triangle
$D=\{(u,v):v<u\}$ is one half of the integral over $[0,1]^2$. Thus, 
$$
\iint_D g_M(u,v) \dd u \dd v = \frac12 \iint_{[0,1]^2} g_M(u,v) \dd u \dd v
=
 \frac12 \iint \phi_M(x,y) \nu(\dd x)\nu(\dd y).
$$
Moreover, from \eqref{eq:trriem},
\begin{equation}\label{eq:SM-limit}
\lim_{n\to \infty}S_{n,M} = 
\frac12
\iint \phi_M(x,y) \nu(\dd x)\nu(\dd y)=:I_M.
\end{equation}

It remains to remove the truncation. Define
$h_M(r)=(-M-\log r)_+, r>0.$
Then
$$
0\le S_{n,M}-S_n
=
\frac1{n^2}
\sum_{1\le i<j\le n}
h_M(\ell_j-s_i).
$$
Let
$
d_{ij}=\ell_j-s_i
=
Q\left(\frac{j-\frac12}{n}\right)
-
Q\left(\frac{i}{n}\right)
$
and
$$
U_i=\left(\frac{i}{n},\frac{i+\frac12}{n}\right], 1\le i\le n-1, \quad V_j=\left(\frac{j-1}{n},\frac{j-\frac12}{n}\right], 2\le j\le n.
$$
If $j\ge i+1$, then for every $(u,v)\in U_i\times V_j$,
$$
|Q(v)-Q(u)|
\le
Q\left(\frac{j-\frac12}{n}\right)
-
Q\left(\frac{i}{n}\right)
=d_{ij}.
$$
Since $h_M$ is decreasing, this implies
$$
h_M(d_{ij})
\le
h_M(|Q(v)-Q(u)|),
\qquad (u,v)\in U_i\times V_j.
$$
As $|U_i\times V_j|=1/(4n^2)$, we get
$$
\frac1{n^2}h_M(d_{ij})
\le
4\int_{U_i\times V_j}
h_M(|Q(v)-Q(u)|) \dd u \dd v.
$$
Therefore,
$$
0\le S_{n,M}-S_n
\le
4\int_0^1\int_0^1
h_M(|Q(v)-Q(u)|) \dd u \dd v=
4\iint
h_M(|x-y|) \nu(\dd x)\nu(\dd y).
$$
It follows from \eqref{eq:SM-limit} that
$
\limsup_{n\to\infty}S_n
\le
I_M.$
Since $\phi_M$ is bounded above, sending $M\to\infty$ gives
$$
\limsup_{n\to\infty}S_n
\le
\frac12
\iint \log|x-y| \nu(\dd x)\nu(\dd y).
$$
If the right-hand side is $-\8$, we are done. Now assume
$\iint \log|x-y| \nu(\dd x)\nu(\dd y)>-\infty.$
 Since
$h_M(|x-y|)\downarrow0$ as $M\to\infty$, dominated convergence gives
\begin{align}\label{eq:discrete-tail}
 \lim_{M\to\infty}
\limsup_{n\to\infty}
(S_{n,M}-S_n)=0.
\end{align}
Using $
S_n= S_{n,M}-(S_{n,M}-S_n),
$ \eqref{eq:SM-limit}, \eqref{eq:discrete-tail} and letting $M\to\infty$, we obtain
$$
\liminf_{n\to\infty}S_n
\ge
\frac12
\iint \log|x-y| \nu(\dd x)\nu(\dd y).
$$
The proof is complete.
\end{proof}

\begin{proof}[Proof of Lemma \ref{x:powemp}]
We apply \cite{WZbeta}*{Theorem 1.1}
by taking $D=[0,\infty)$ and $\mathcal O=(0,\infty)$. 
Alternatively, one can also use \cite{GZ19} and verify the conditions there.
Write $V_n=W_n-(\beta n)^{-1}\log\omega_n$, where
\begin{equation}\label{e:powwgt}
W_n(x)=W_p(x)=\frac{x^{p/2}}p,\qquad
\omega_n(x)=x^{\beta/2-1},\qquad x>0,
\end{equation}
and $W_p(0)=0$ and $\omega_n(0)=1$. 
Since $\beta n\to\8$ and $W_p$ is continuous with polynomial growth, Assumptions (A0)--(A2) in \cite{WZbeta} hold.
Moreover, $\log\omega_n= (\frac\bt2-1)\log x=o(\beta n)$ locally uniformly
on $\mathcal O=(0,\infty)$, we have (A3)(i).
For (A3)(ii), let $\theta_n=(\beta n)^{-1/2}$.
Using \eqref{eq:Frelb} with $\alpha=p/2$ and $\varepsilon=1/(2p)$, we obtain
$C>0$ such that 
\begin{equation}\label{e:powgro}
W_p(x)-\log(1+x)\geq\frac{x^{p/2}}{2p}-C,\qquad x\geq0.
\end{equation}
For $x_i\geq0$, \eqref{e:pairbd} gives
$\prod_{i<j}|x_i-x_j|^{\theta_n\beta}
\leq\prod_{i=1}^n(1+x_i)^{\theta_n\beta(n-1)}.$
By \eqref{e:powgro}, we deduce
\begin{align*}
\mathcal R_n(\theta_n)
&:=\int_{\mathbb R_+^n}
e^{-\theta_n\beta n\sum_{i=1}^nW_p(x_i)}
\prod_{1\leq i<j\leq n}|x_i-x_j|^{\theta_n\beta}
\prod_{i=1}^n\omega_n(x_i)\dd\bx\\
&\leq\left[\int_0^\infty
e^{-\theta_n\beta n[W_p(x)-\frac{n-1}{n}\log(1+x)]}
x^{\beta/2-1}\dd x\right]^n\\
&\leq e^{C\theta_n\beta n^2}
\left[\int_0^\infty
e^{-\frac{\theta_n\beta n}{2p}x^{p/2}}x^{\beta/2-1}\dd x\right]^n
=e^{C\theta_n\beta n^2}
\left[\frac2p\Gamma\left(\frac\beta p\right)
\left(\frac{\theta_n\beta n}{2p}\right)^{-\beta/p}\right]^n.
\end{align*}
Since $p$ and $\beta$ are fixed,
$\log\mathcal R_n(\theta_n)\leq C\theta_n\beta n^2+O(n)
=o(\beta n^2)$, proving (A3)(ii).
Since $D$ is a convex set and $W_p$ is continuous on $D$, \cite{WZbeta}*{Lemma 2.3} verifies Assumption B of that paper for every $\mu\in\mathcal M_1(\mathbb R_+)$ with $I_{W_p}(\mu)<\infty$.
Therefore, \cite{WZbeta}*{Theorem 1.1}  gives the asserted LDP.

By Remark \ref{rem:squlm}, $\eta^{(p)}$ is the equilibrium measure for $W_p$ and $m_{p/2}(\eta^{(p)})=1$. Taking $s=1/p$ in Lemma \ref{lemhalfp}(ii), we find
$\Sigma(\eta^{(p)})=2\log(b_p/2)-1/p$. Thus
\[
\lim_{n\to\infty}\frac{1}{\beta n^2}\log\bar Z_{n,p,\beta}^{(n)}
=\frac12\Sigma(\eta^{(p)})-\frac1p m_{p/2}(\eta^{(p)})
=\log\frac{b_p}{2}-\frac{3}{2p},
\]
which proves \eqref{eq:powZ}.
\end{proof}

\section{Difference of Schatten balls}\label{sec:diffball}
We first state the asymptotics for the difference volumes of unnormalized balls in the two parameter regimes. Recall $\Phi_{p,q}$ as defined in \eqref{eq:phipq}.

\begin{proposition}\label{thm:diffball}
Let  $\beta \in \{1,2,4\}$ and $0<p<q\leq\infty$. For $r>r_{p,q}$, we have the following:
\begin{enumerate}
 \item[(i)] Self-adjoint case.
$$ \lim_{n\to\infty}
\frac{1}{\beta n^{1+p/q}}
\log \frac{
\vol\left(
\mathbb B_{p,\beta}^n
\setminus
r n^{\frac1q-\frac1p}\mathbb B_{q,\beta}^n
\right)
}{
\vol\left(\mathbb B_{p,\beta}^n\right)
}
=
-\Phi_{p,q}(r),$$

\item[(ii)] General case.
$$\lim_{n\to\infty}
\frac{1}{\beta n^{1+p/q}}
\log
\frac{
\vol\left(
\mathbb{SB}_{p,\beta}^n
\setminus
r n^{\frac1q-\frac1p}\mathbb{SB}_{q,\beta}^n
\right)
}{
\vol\left(\mathbb{SB}_{p,\beta}^n\right)
}
=
-2\Phi_{p,q}(r).$$
\end{enumerate}
Moreover, $\Phi_{p,q}(r)$ is continuous and strictly increasing with values in $(0,\8)$ whenever $r>r_{p,q}$.
\end{proposition}

\begin{proposition}
\label{thm:diffpgq}
Let $\beta\in\{1,2,4\}$ and $0<q<p\leq\infty$.  Then for $r_{p,q}< r<1$ we have
$$
 \lim_{n\to\infty}\frac1{\beta n^2}
 \log
 \frac{
 \vol\left(
 \mathbb B_{p,\beta}^n
 \setminus
 r n^{\frac1q-\frac1p}\mathbb B_{q,\beta}^n
 \right)
 }{
 \vol(\mathbb B_{p,\beta}^n)
 }
 =-\Phi_{p,q}(r).
$$
The same conclusion holds for the general case with the rate function $2\Phi_{p,q}(r)$.
Moreover, $\Phi_{p,q}(r)$ is continuous and strictly increasing with values in $(0,\8)$ whenever $r_{p,q}<r<1$.
\end{proposition}

\begin{proof}[Proof of Theorem \ref{th:balldf}]
We first prove the self-adjoint case. Let $Z_n$ be uniformly distributed on $\mathbb B_{p,\beta}^n$. 
Let $\rho_{q,n}(t)$ be defined in \eqref{e:sarad}, $r_n=n^{1/p-1/q}\rho_{q,n}(t)$ and $T_n=n^{1/p-1/q}\|\lambda(Z_n)\|_q$.
The volume normalization gives
$$
\vol\left(
\mathbb D_{p,\beta}^n
\setminus
t\mathbb D_{q,\beta}^n
\right)
=
\frac{
\vol\left(
\mathbb B_{p,\beta}^n
\setminus
r_n n^{\frac1q-\frac1p}\mathbb B_{q,\beta}^n
\right)
}{
\vol\left(\mathbb B_{p,\beta}^n\right)
} =\mathbb P(T_n>r_n).
$$
By the volume asymptotic formula \eqref{eq:volkpt},
$
r_n\to r:=c_{p,q}t.
$
Since $t_{p,q}^*=r_{p,q}/c_{p,q}$, the assumptions imply $r>r_{p,q}$ when $p<q$, and $r_{p,q}<r<1$ when $p>q$.
Choose $\varepsilon>0$ so that for all sufficiently large $n$,
$
\{T_n>r+\varepsilon\}
\subseteq\{T_n>r_n\}
\subseteq\{T_n>r-\varepsilon\}.
$
Proposition \ref{thm:diffball} or \ref{thm:diffpgq} therefore gives
\begin{align*}
-\Phi_{p,q}(r+\varepsilon)
&\leq\liminf_{n\to\infty}
\frac1{\beta n^{1+(p/q)\wedge1}}
\log\mathbb P(T_n>r_n)\\
&\leq\limsup_{n\to\infty}
\frac1{\beta n^{1+(p/q)\wedge1}}
\log\mathbb P(T_n>r_n)
\leq-\Phi_{p,q}(r-\varepsilon).
\end{align*}
By continuity of $\Phi_{p,q}$ at $r$, letting $\varepsilon\downarrow0$ yields 
$$
\lim_{n\to\infty}
\frac{1}{\beta n^{1+(p/q)\wedge1}}
\log
\vol\left(
\mathbb D_{p,\beta}^n
\setminus
t\mathbb D_{q,\beta}^n
\right)
=
-\Phi_{p,q}(c_{p,q}t).
$$

For the general case, we use $\varrho_{q,n}(t)$ from \eqref{e:marad}. The volume asymptotic \eqref{eq:vossba} implies that
$$
r_n:=n^{\frac1p-\frac1q}\varrho_{q,n}(t) \to c_{p,q}t>r_{p,q}.
$$
Moreover,
$$
\vol\left(
\mathbb{SD}_{p,\beta}^n
\setminus
t\mathbb{SD}_{q,\beta}^n
\right)
=
\frac{
\vol\left(
\mathbb{SB}_{p,\beta}^n
\setminus
r_n n^{\frac1q-\frac1p}\mathbb{SB}_{q,\beta}^n
\right)
}{
\vol\left(\mathbb{SB}_{p,\beta}^n\right)
}.
$$
The rest of proof follows verbatim that of self-adjoint case.
\end{proof}

\begin{remark}\label{rem:Plim}
For fixed $0<p<\infty$, the rate function $\Phi_{p,q}(r)$ is discontinuous at
$q=\infty$. Indeed, for every fixed $r>b_p$,
\begin{align*}
\lim_{q\to\infty}\Phi_{p,q}(r)
&=
\lim_{q\to\infty}
\frac{1}{2p}
\left(r^q-r_{p,q}^q\right)^{p/q}=
\frac{r^p}{2p} \neq  \Phi_{p,\infty}(r).
\end{align*}
Let us explain heuristically how the phase transition arises in the self-adjoint case. Using the probabilistic representation \eqref{eq:lamzn}, the volume ratio $\vol\left(
\mathbb B_{p,\beta}^n
\setminus
r n^{\frac1q-\frac1p}\mathbb B_{q,\beta}^n
\right)/
\vol\left(\mathbb B_{p,\beta}^n\right)
$ is ultimately determined by the event $(\frac1n\sum_{i=1}^n |Y_i|^q)^{1/q}>r$ for  $q\in(p,\8)$ where $(Y_1,\dots,Y_n)$ is rescaled from $(X_1,\dots,X_n)$; see the proof of Proposition \ref{thm:diffball}. This event 
is realized by an outlier of order $n^{1/q}r$. It causes a `cost' of order $n^{1+ p/q}r^p$ from the potential $n |x|^p$, whereas its logarithmic interaction with the remaining particles is only of order $n \log n$. Hence, the interaction term is negligible at the speed $n^{1+p/q}$ and we only get a term in $r^p$ for the rate function. In contrast, when $q=\infty$, the relevant event is $\max_{1\leq i\leq n}|Y_i|>r,$ and the outlier is of constant order so that the costs from the potential and interaction are both of order $n$.  The large deviation speed is then $n$. These interaction terms give rise to the hypergeometric and logarithmic terms appearing in $\Phi_{p,\infty}$. Consequently, continuity of $\Phi_{p,q}$ at $q=\infty$ should not be expected. 
\end{remark}

The proofs of Propositions \ref{thm:diffball} and \ref{thm:diffpgq} rely on four technical lemmas. Lemmas  \ref{le:pc-Freud} and \ref{le:snoutlier} will be used for the self-adjoint case, while
Lemmas \ref{le:sv-pc} and \ref{le:sv-o} provide the corresponding estimates for the general case. The proofs of these lemmas are postponed to the end of
this section. 

\begin{lemma}\label{le:pc-Freud}
Let $V$ and $\beta$ be as in Theorem \ref{th:momldp}.
For $s>0$, write
\begin{equation}\label{eq:defzkp}
Z_k^{(s)}
=\int_{\mathbb R^k}e^{-\beta s\sum_{i=1}^kV(x_i)}
\prod_{1\leq i<j\leq k}|x_i-x_j|^\beta\dd\bx,
\qquad Z_0^{(s)}=1.
\end{equation}
There exists $C>0$ such that uniformly for $0\leq m\leq n$,
\begin{equation}\label{eq:Zratio}
\log Z_{n-m}^{(n)}-\log Z_n^{(n)}\leq C\beta nm.
\end{equation}
Moreover,
\begin{equation}\label{eq:widezz}
\log Z_{n-1}^{(n)}-\log Z_n^{(n)}\geq-C\beta n.
\end{equation}
\end{lemma}

\begin{lemma}\label{le:snoutlier}
Let $\beta\in\{1,2,4\}$ and $0<p<q\leq\infty$, and let $\bY=(Y_1,\ldots,Y_n)$ have joint density 
\begin{align}\label{eq:den-sa}
\widehat Z_{n,p,\beta}^{-1} \exp\left\{-\frac{\beta n}{2p}\sum_{i=1}^n |y_i|^p\right\} \prod_{1\le i<j\le n}|y_i-y_j|^\beta.
\end{align}
Let $B_{\infty,n}= \max_{1\le i\le n}|Y_i|$ and $B_{q,n}=
\left(\frac1n\sum_{i=1}^n |Y_i|^q\right)^{1/q}
$ for $0< q <\infty .$ Then, for every $s>r_{p,q}$,
\begin{equation}\label{eq:oneoutq}
\lim_{n\to\infty}
\frac{1}{\beta n^{1+p/q}}
\log
\mathbb P(B_{q,n}>s)
=
-\Phi_{p,q}(s),
\end{equation}
where $\Phi_{p,q}(s)$ is given in \eqref{eq:phipq}.
\end{lemma}

Let
\begin{align}\label{eq:barzk}
\bar Z_{k,p,\beta}^{(n)}
=
\int_{[0,\infty)^k}
\prod_{1\le i<j\le k}|y_i-y_j|^\beta
\prod_{i=1}^k y_i^{\frac{\beta}{2}-1}
\exp\left\{
-\frac{\beta n}{p}\sum_{i=1}^k y_i^{p/2}
\right\}
 \dd\by,
\end{align}
with the convention $\bar Z_{0,p,\beta}^{(n)}=1$. 

\begin{lemma}\label{le:sv-pc}
Let $\beta\in\{1,2,4\}$ and $0<p<\infty$. There exists a constant $C_0>0$ such that, uniformly for
$0\le m\le n$,
$$
\log \bar Z_{n-m,p,\beta}^{(n)}- \log \bar Z_{n,p,\beta}^{(n)}
\le C_0\beta nm.
$$
Moreover, there exists a constant $C_{p,\beta}>0$ such that
\begin{align}\label{eq:barzz}
\log \bar Z_{n-1,p,\beta}^{(n)}
- \log \bar Z_{n,p,\beta}^{(n)}
\ge
-C_{p,\beta}n.
\end{align}
\end{lemma}

\begin{lemma}\label{le:sv-o}
Let $\beta\in\{1,2,4\}$ and $0<p<q\le \infty$ and let $\bY=(Y_1,\ldots,Y_n)$ have the distribution in Lemma \ref{x:powemp}.
Let $\tilde B_{\infty,n}=\max_{1\le i\le n}Y_i$ and
$
\tilde B_{q,n}=
\left(\frac1n\sum_{i=1}^n Y_i^{q/2}\right)^{2/q},
 0<q<\infty.
$
Then, for every $s>r_{p,q}$,
$$
\lim_{n\to\infty}
\frac{1}{\beta n^{1+p/q}}
\log
\mathbb P(\tilde B_{q,n}>s^2)
=
-2\Phi_{p,q}(s),
$$
where $\Phi_{p,q}(s)$ is given in \eqref{eq:phipq}.
\end{lemma}

The following LDP is needed for the case $q=\8$.

\begin{lemma}\label{x:powmax}
For the ensemble in Lemma \ref{x:powemp},  $\lambda_{\max}=\max_{1\leq i\leq n}\lambda_i$ satisfies an LDP on $\mathbb R$ with speed $\beta n$ and good rate function
\begin{equation}\label{e:powmax}
J_{W_p}(x)=
\begin{cases}
2\Phi_{p,\infty}(\sqrt x),&x\geq b_p^2,\\
+\infty,&x<b_p^2,
\end{cases}
\end{equation}
where $\Phi_{p,\infty}$ is defined in \eqref{eq:phipq}.
Moreover, the rate function $J_{W_p}$ is continuous and strictly increasing on $[b_p^2,\infty)$.
\end{lemma}

\begin{proof}[Proof of Proposition \ref{thm:diffball}]
The asserted properties of $\Phi_{p,q}$ follow from \eqref{eq:phipq} for $p<q<\infty$ and from Lemma  \ref{x:powmax} for $q=\infty$. 

\noindent\emph{(i) Self-adjoint case}. 
Let $Z_n$ be uniformly distributed on $\mathbb B_{p,\beta}^n$, and set $T_n=n^{\frac1p-\frac1q}\|\lambda(Z_n)\|_q,$ where $T_n=n^{1/p}\|\lambda(Z_n)\|_\infty$ for $q=\infty$.
Then
$$
\frac{
\vol\left(
\mathbb B_{p,\beta}^n
\setminus
r n^{\frac1q-\frac1p}\mathbb B_{q,\beta}^n
\right)
}{
\vol\left(\mathbb B_{p,\beta}^n\right)
}
=
\mathbb P(T_n>r).
$$
By the representation of $Z_n$ in \eqref{eq:lamzn}, after the change of variables
$
Y_i:=\left(\frac{2p}{\beta n}\right)^{1/p}X_i, 1\le i\le n,
$
the vector $\bY=(Y_1,\ldots,Y_n)$ has joint density \eqref{pdf} with $D=\mathbb R$ and $V_n=V_p$.
Define
$$
A_n=\left(\frac1n\sum_{i=1}^n |Y_i|^p\right)^{1/p},
\quad B_{\infty,n}= \max_{1\leq i\leq n} |Y_i|, \quad
B_{q,n}=\left(\frac1n\sum_{i=1}^n |Y_i|^q\right)^{1/q}, \quad q<\infty.
$$
Then $T_n
\stackrel d=
U^{1/d_{\beta,n}}\frac{B_{q,n}}{A_n}.
$

By \cite{KPT20c}*{Theorem 4.1}, the random elements
$
\left(L_n^Y,A_n^p\right)
=
\left(
\frac1n\sum_{i=1}^n\delta_{Y_i},
\frac1n\sum_{i=1}^n |Y_i|^p
\right)
$
satisfy an LDP on $\mathcal M_1(\mathbb R)\times[0,\infty)$ with speed $\beta n^2$ and a good rate function $\mathcal I_1$, 
whose unique minimizer is $(\nu^{(p)},1)$ satisfying $m_p(\nu^{(p)})=1$. Here $\nu^{(p)}$ is given in \eqref{eq:nueta}.
Fix $\delta>0$, and let
$$
D_\delta
=
\left\{
m\in[0,\infty): |m^{1/p}-1|\ge \delta
\right\},
\qquad
F_\delta
=
\mathcal M_1(\mathbb R)\times D_\delta.
$$
Since $\mathcal I_1$ is a good rate function with unique minimizer $(\nu^{(p)},1)$, 
$
c_\delta
:=
\inf_{(\mu,m)\in F_\delta}\mathcal I_1(\mu,m)
>0.
$
Noting that $A_n^p=\frac1n\sum_{i=1}^n |Y_i|^p$, we have
$
\{|A_n-1|>\delta\}
\subseteq
\left\{
(L_n^Y,A_n^p)\in F_\delta
\right\}.
$
Applying the LDP upper bound to the closed set $F_\delta$ yields $
\limsup_{n\to\infty}
\frac1{n^2}
\log
\mathbb P(|A_n-1|>\delta)
\le
-c_\delta.
$
Consequently, for every $\delta>0,$
\begin{align}\label{eq:supexpan}
\limsup_{n\to\infty}
\frac{1}{\beta n^{1+p/q}}
\log
\mathbb P(|A_n-1|>\delta)
=
-\infty.	
\end{align}
Moreover, $
\mathbb P\bigl(U^{1/d_{\beta,n}}<1-\delta\bigr)
=(1-\delta)^{d_{\beta,n}}.$
Since $d_{\beta,n}=\Theta(n^2)$ and $1+p/q<2$, we find
\begin{align}\label{eq:supexpun}
\limsup_{n\to\infty}
\frac{1}{\beta n^{1+p/q}}
\log
\mathbb P\bigl(U^{1/d_{\beta,n}}<1-\delta\bigr)
=
-\infty.	
\end{align}

Fix $\delta\in(0,1)$ such that $r(1-\delta)>r_{p,q}.$ 
Since $U^{1/d_{\beta,n}}\le1$, on the event $\{A_n\ge1-\delta\}$ we have
$$
T_n
\stackrel{d}=
U^{1/d_{\beta,n}}\frac{B_{q,n}}{A_n}
\le
\frac{B_{q,n}}{1-\delta}.
$$
Therefore, $\pz(T_n>r)\le  \pz(B_{q,n}>r(1-\delta))+\pz(
A_n<1-\delta).$ 
Using \eqref{eq:supexpan} and \eqref{eq:oneoutq}, we get
$$
\limsup_{n\to\infty}
\frac{1}{\beta n^{1+p/q}}
\log \mathbb P(T_n>r)
\le
-\Phi_{p, q}(r(1-\delta)).
$$
Letting $\delta\downarrow0$, the continuity of $\Phi_{p,q}$ on
$(r_{p,q},\infty)$ gives
$$
\limsup_{n\to\infty}
\frac{1}{\beta n^{1+p/q}}
\log \mathbb P(T_n>r)
\le
-\Phi_{p, q}(r).
$$

For the lower bound, observing that
$$
\pz\left(\left\{
B_{q,n}>r\frac{1+\delta}{1-\delta}
\right\}
\cap
\{A_n\le1+\delta\}
\cap
\{U^{1/d_{\beta,n}}\ge1-\delta\}\right)
\le
\pz(T_n>r),
$$
we obtain
\begin{align}\label{eq:psnr}
\mathbb P(T_n>r)
\ge
\mathbb P\left(
B_{q,n}>r\frac{1+\delta}{1-\delta}
\right)
-
\mathbb P(A_n>1+\delta)
-
\mathbb P(U^{1/d_{\beta,n}}<1-\delta).	
\end{align}
Applying \eqref{eq:supexpan}, \eqref{eq:supexpun}, and \eqref{eq:oneoutq} to \eqref{eq:psnr}, we arrive at 
$$
\liminf_{n\to\infty}
\frac{1}{\beta n^{1+p/q}}
\log \mathbb P(T_n>r)
\ge
-\Phi_{p, q}\left(r\frac{1+\delta}{1-\delta}\right).
$$
As $\delta\downarrow0$, continuity yields the desired lower bound and the proof is complete.

\noindent\emph{(ii) General case.} 
Let $Z_n$ be uniformly distributed on $\mathbb{SB}_{p,\beta}^n$, and let
$s(Z_n)=(s_1(Z_n),\ldots,s_n(Z_n))$ denote its singular value vector. 
Let $ T_n=n^{\frac1p-\frac1q}\|s(Z_n)\|_q
$
for $0<q\le\infty$. Then
$$
\frac{
\vol\left(
\mathbb{SB}_{p,\beta}^n
\setminus
r n^{\frac1q-\frac1p}\mathbb{SB}_{q,\beta}^n
\right)
}{
\vol\left(\mathbb{SB}_{p,\beta}^n\right)
}
=
\mathbb P(T_n>r).
$$
Under the representation \eqref{eq:storex}, the rescaled vector $\bY=(p/\beta)^{2/p}\bX$ has the distribution in Lemma \ref{x:powemp}.
Define
$$
A_n
=
\left(\frac1n\sum_{i=1}^nY_i^{p/2}\right)^{2/p}, \quad \tilde B_{q,n}
=
\left(\frac1n\sum_{i=1}^n Y_i^{q/2}\right)^{2/q},\quad
 0<q<\infty, \quad \tilde B_{\infty,n}=\max_{1\le i\le n}Y_i.
$$
Thus, $T_n^2
\stackrel d=
U^{2/(\beta n^2)}
\frac{\tilde B_{q,n}}{A_n}.
$
By abuse of notation, we simply write $T_n^2$ for the right-hand side representation in the following proof.

Let $L_n^Y=n^{-1}\sum_{i=1}^n\delta_{Y_i}$.
We follow the joint LDP proof in \cite{KPT20c}*{Section~8}. The technical input \cite{HP00}*{Theorem~5.5.1} is replaced by Lemma \ref{x:powemp} so that the case $\bt=1$ is now valid.
Thus $(L_n^Y,A_n^{p/2})$ satisfies an LDP with speed $\beta n^2$ and a good rate function $\mathcal J$ whose unique zero is $(\eta^{(p)},1)$. As in the proof of \eqref{eq:supexpan}, we obtain for every $\delta>0$ 
\begin{equation}\label{eq:svAexp}
\limsup_{n\to\infty}
\frac1{\beta n^{1+p/q}}
\log
\mathbb P(|A_n-1|>\delta)
=
-\infty.
\end{equation}
Alternatively, one can directly use the fact that $A_n^{p/2}$ follows a Gamma distribution to derive the above concentration result.
Furthermore, since $U$ is uniformly distributed on $(0,1)$ and $1+p/q<2$, we also have
\begin{equation}\label{eq:svUexp}
\limsup_{n\to\infty}
\frac1{\beta n^{1+p/q}}
\log
\mathbb P\left(U^{2/(\beta n^2)}<1-\delta\right)
=
-\infty.
\end{equation}

We first prove the upper bound. Fix $\delta\in(0,1)$ such that $r\sqrt{1-\delta}>r_{p,q}.$
Since $U^{2/(\beta n^2)}\le1$, on the event $\{A_n\ge1-\delta\}$ we have
$T_n^2
\le
\frac{\tilde B_{q,n}}{1-\delta}.$
Hence
$$
\{T_n>r\}
=
\{T_n^2>r^2\}
\subset
\left\{
\tilde B_{q,n}>r^2(1-\delta)
\right\}
\cup
\{A_n<1-\delta\}.
$$
By \eqref{eq:svAexp} and Lemma \ref{le:sv-o} applied with $s=r\sqrt{1-\delta},$  we get
$$
\limsup_{n\to\infty}
\frac1{\beta n^{1+p/q}}
\log\mathbb P(T_n>r)
\le
-2\Phi_{p,q}\bigl(r\sqrt{1-\delta}\bigr).
$$
As $\delta\downarrow0$, the continuity of $2\Phi_{p,q}$ on $(r_{p,q},\infty)$ yields the desired upper bound.

We next prove the lower bound. For $\delta\in(0,1)$, observe that
$$
\left\{
\tilde B_{q,n}>r^2\frac{1+\delta}{1-\delta}
\right\}
\cap
\{A_n\le1+\delta\}
\cap
\left\{
U^{2/(\beta n^2)}\ge1-\delta
\right\}
\subset
\{T_n>r\}.
$$
It follows that
$$
 \mathbb P(T_n>r)
\ge
\mathbb P\left(
\tilde B_{q,n}>r^2\frac{1+\delta}{1-\delta}
\right)
-
\mathbb P(A_n>1+\delta) -
\mathbb P\left(
U^{2/(\beta n^2)}<1-\delta
\right).
$$
Applying Lemma \ref{le:sv-o} with
$s=r\sqrt{\frac{1+\delta}{1-\delta}},$
and using the estimates \eqref{eq:svAexp} and
\eqref{eq:svUexp}, we obtain
$$
\liminf_{n\to\infty}
\frac1{\beta n^{1+p/q}}
\log\mathbb P(T_n>r)
\ge
-2\Phi_{p,q}\left(
r\sqrt{\frac{1+\delta}{1-\delta}}
\right).
$$
Letting $\delta\downarrow0$, continuity yields the lower bound and the proof is complete.
\end{proof}

\begin{proof}[Proof of Proposition \ref{thm:diffpgq}]
\noindent\emph{(i) Self-adjoint case.}
Let $Z_n$ be uniformly distributed on $\mathbb B_{p,\beta}^n$.
By \cite{KPT20c}*{Theorems 1.1 and 1.3}, the empirical measure
$
\mu_n=\frac1n\sum_{i=1}^n\delta_{n^{1/p}\lambda_i(Z_n)}
$
satisfies an LDP with speed $\beta n^2$ and good rate function
$\mathcal K$ given by \eqref{eq:sarate}. For
$p<\infty$, write 
$$\mathcal D_p=\{\mu\in\mathcal M_1(\mathbb R):
m_p(\mu)\leq1\},
$$ 
where $m_p(\mu)$ was defined in \eqref{eq:mumom}. For $p=\infty$, let
$\mathcal D_\infty=\{\mu\in\mathcal M_1(\mathbb R):
\operatorname{supp}(\mu)\subseteq[-1,1]\}$ and
$C_\infty=-\frac12\log2$. Then
$\mathcal K(\mu)=-\frac12\Sigma(\mu)+C_\infty$ on
$\mathcal D_\infty$ and $\mathcal K(\mu)=+\infty$ otherwise. 
 
Let $T_n=n^{\frac1p-\frac1q}\|\lambda(Z_n)\|_q$. Then
\begin{equation}\label{eq:savolt}
 \frac{
 \vol\left(
 \mathbb B_{p,\beta}^n
 \setminus
 r n^{\frac1q-\frac1p}\mathbb B_{q,\beta}^n
 \right)
 }{
 \vol(\mathbb B_{p,\beta}^n)
 }
 =\mathbb P\bigl(\|\lambda(Z_n)\|_q
 >rn^{\frac1q-\frac1p}\bigr)
 =\mathbb P(T_n>r)
 =\mathbb P\bigl(\mu_n\in\widetilde F_q(r)\bigr),
\end{equation}
where
$
\widetilde F_q(r)=\{\mu\in\mathcal D_p:m_q(\mu)>r^q\}.
$

Since $\mu_n\in\mathcal D_p$ almost surely and
$\mathcal K=+\infty$ on $\mathcal D_p^c$, the above LDP can be restricted to $\mathcal D_p$ endowed with the topology induced by the weak topology on $\mathcal M_1(\mathbb R)$. Note that the map
$
m_q:\mathcal D_p \to [0,1]
$
is continuous. Indeed, let $\mu_k,\mu\in\mathcal D_p$ with $\mu_k\Rightarrow\mu$. If $p<\infty$, since $p>q$,  we have for every $R>0$,
\begin{equation*}
\sup_k
\int_{\{|x|>R\}}|x|^q \mu_k(\dd x)
\leq
R^{q-p}\sup_km_p(\mu_k)
\leq
R^{q-p} \to 0, \quad \text{as } R\to\infty.
\end{equation*}
Thus the $q$-moments are uniformly integrable, and hence
$m_q(\mu_k)\to m_q(\mu)$. For $p=\infty$, the same conclusion follows
from the common compact support $[-1,1]$.
Therefore, for $r_{p,q}<r<1$, the set $\widetilde F_q(r)$ is open
in $\mathcal D_p$, whereas
$\widehat F_q(r)=\{\mu\in\mathcal D_p:m_q(\mu)\geq r^q\}$ is closed.
The LDP lower and upper bounds yield
\begin{align}
 -\inf_{\mu\in\widetilde F_q(r)}
 \mathcal K(\mu)
 \leq
 \liminf_{n\to\infty}\frac1{\beta n^2}
 \log\mathbb P(T_n>r)
 \leq
 \limsup_{n\to\infty}\frac1{\beta n^2}
 \log\mathbb P(T_n>r)
\leq
 -\inf_{\mu\in\widehat F_q(r)}
 \mathcal K(\mu).
 \label{eq:satldp}
\end{align}

We next show that
\begin{equation}\label{eq:sainfs}
\inf_{\mu\in\widetilde F_q(r)}\mathcal K(\mu)
=
\inf_{\mu\in\widehat F_q(r)}\mathcal K(\mu).
\end{equation}
Here the inequality ``$\geq$'' is trivial. For the reverse inequality, it suffices to approximate a measure $\mu\in\mathcal D_p$ satisfying $\mathcal K(\mu)<\infty$ and $m_q(\mu)=r^q$.
Let $U_r$ be the uniform probability measure
\begin{equation}\label{eq:Ur}
U_r(\dd x)=\frac4{1-r}
\mathbf1_{[(1+r)/2,(3+r)/4]}(x)\dd x.
\end{equation}
Then $m_q(U_r)>r^q$ and $m_p(U_r)<1$ if $p<\infty$, while $\operatorname{supp}(U_r)\subseteq[-1,1]$ if $p=\infty$. Hence, for $0<\varepsilon<1$, the measure $\mu_\varepsilon=(1-\varepsilon)\mu+\varepsilon U_r$ belongs to $\widetilde F_q(r)$. Since $U_r$ has bounded density and compact support, its logarithmic potential is bounded below and grows at most logarithmically.  
By concavity of $\Sigma$,  $\Sigma(\mu_\eps)$ is finite, and bilinearity gives $\mathcal K(\mu_\varepsilon)\to\mathcal K(\mu)$ as $\varepsilon\downarrow0$. This proves \eqref{eq:sainfs}.

It remains to identify the variational formula. First suppose that
$0<q<p<\infty$. By
\eqref{eq:sarate},
$$
\inf_{\mu\in\widehat F_q(r)}\mathcal K(\mu)
=
C_p-\frac12
\sup\left\{
\Sigma(\mu):
\mu\in\mathcal M_1(\mathbb R),\
m_p(\mu)\leq1, m_q(\mu)\geq r^q
\right\}.
$$
We have to show
$$
\sup_{\substack{\mu\in\mathcal M_1(\mathbb R)\\
m_p(\mu)=1, m_q(\mu)=r^q}}\Sigma(\mu)
\ge
\sup_{\substack{\mu\in\mathcal M_1(\mathbb R)\\
m_p(\mu)\leq1, m_q(\mu)\geq r^q}}\Sigma(\mu).
$$
To this end, fix
$\mu\in\mathcal M_1(\mathbb R)$ such that
$m_p(\mu)\leq1, m_q(\mu)\geq r^q$ and $\Sigma(\mu)>-\infty.$
Let
$a=m_p(\mu)^{-1/p}\geq1$ and
$\widetilde\mu=(x\mapsto ax)_\#\mu$. The linear scaling gives
$$m_p(\widetilde\mu)=1, \qquad
m_q(\widetilde\mu)=a^qm_q(\mu)\geq r^q, \qquad 
\Sigma(\widetilde\mu)=\Sigma(\mu)+\log a\geq\Sigma(\mu).$$
Thus the $p$-moment constraint can be fulfilled without decreasing $\Sigma.$
If $m_q(\widetilde\mu)=r^q$, there is nothing further to prove.
Otherwise, recalling the measure
$\nu^{(p)}$ defined in \eqref{eq:nueta}, 
let
\[
\mu_t=(1-t)\widetilde\mu+t\nu^{(p)}, \quad 
t=\frac{m_q(\widetilde\mu)-r^q}{
m_q(\widetilde\mu)-r_{p,q}^q}\in(0,1).
\]
Then
$m_p(\mu_t)=1$ and $m_q(\mu_t)=r^q$. Since $\nu^{(p)}$ maximizes
$\Sigma$ over $\mathcal D_p$, it follows from concavity that
$$\Sigma(\mu_t)\geq(1-t)\Sigma(\widetilde\mu)
+t\Sigma(\nu^{(p)})\geq\Sigma(\widetilde\mu)
\geq\Sigma(\mu).$$ 
Hence, $
\inf_{\mu\in\widehat F_q(r)}\mathcal K(\mu)
=\Phi_{p,q}(r).
$

Finally, suppose that $p=\infty$. Let
\begin{equation}\label{eq:arclaw}
\nu^{(\infty)}(\dd x)
=
\frac{\mathbf 1_{(-1,1)}(x)}
{\pi\sqrt{1-x^2}} \dd x
\end{equation}
be the arcsine law on $[-1,1]$. It maximizes $\Sigma$ over
$\mathcal D_\infty$ and satisfies
$m_q(\nu^{(\infty)})=r_{\infty,q}^q<r^q$. Thus, if
$\mu\in\mathcal D_\infty$ and $m_q(\mu)>r^q$, a suitable
convex combination of $\mu$ and $\nu^{(\infty)}$ has $q$-moment
$r^q$, remains in $\mathcal D_\infty$, and has free entropy at least $\Sigma(\mu)$. Consequently,
$\inf_{\mu\in\widehat F_q(r)}\mathcal K(\mu)
=\Phi_{\infty,q}(r)$.

Combining
$\inf_{\mu\in\widehat F_q(r)}\mathcal K(\mu)=\Phi_{p,q}(r)$ with
\eqref{eq:sainfs},
\eqref{eq:satldp}, and
\eqref{eq:savolt} proves part~(i).

\noindent\emph{(ii) General case.}
Let $Z_n$ be uniformly distributed on $\mathbb{SB}_{p,\beta}^n$.
As explained in the proof of Proposition \ref{pr:nspqup}, by \cite{KPT20c} the empirical measure
$\eta_n=\frac{1}{n}\sum_{i=1}^n
\delta_{n^{2/p}s_i^2(Z_n)}$ satisfies an LDP on
$\mathcal M_1(\mathbb R_+)$ with speed $\beta n^2$ and good rate function $\mathcal Y$ given by \eqref{eq:defmathy}. Let
\[
\widehat{\mathcal D}_p=
\begin{cases}
    \{\eta\in\mathcal M_1([0,\infty)):
m_{p/2}(\eta)\leq1\},& p<\8,\\
\{\eta\in\mathcal M_1([0,\infty)):
\operatorname{supp}(\eta)\subseteq[0,1]\},& p=\8.
\end{cases}
\]
Then, $\mathcal Y(\eta)=-\frac12\Sigma(\eta)+2C_p$ for $\eta\in \widehat {\mathcal D}_p$ and $\mathcal Y(\eta)=\8$ otherwise. 
Since $\eta_n\in\widehat{\mathcal D}_p$ almost surely, we regard this as an LDP on $\widehat{\mathcal D}_p$ with the relative weak topology.

Let $S_n=n^{\frac1p-\frac1q}\|s(Z_n)\|_q$. Then
$S_n^q=m_{q/2}(\eta_n)$, and
\begin{equation}\label{eq:svvolt}
 \frac{
 \vol\left(
 \mathbb{SB}_{p,\beta}^n
 \setminus
 r n^{\frac1q-\frac1p}\mathbb{SB}_{q,\beta}^n
 \right)
 }{
 \vol(\mathbb{SB}_{p,\beta}^n)
 }
 =\mathbb P(S_n>r)
 =\mathbb P\bigl(\eta_n\in\widetilde G_q(r)\bigr),
\end{equation}
where $\widetilde G_q(r)=\{\eta\in\widehat{\mathcal D}_p:
m_{q/2}(\eta)>r^q\}$.

As in part~(i), the map
$m_{q/2}:\widehat{\mathcal D}_p\to[0,1]$ is continuous in the
relative weak topology. Hence $\widetilde G_q(r)$ is open, whereas
$\widehat G_q(r)=\{\eta\in\widehat{\mathcal D}_p:
m_{q/2}(\eta)\geq r^q\}$ is closed. The LDP bounds yield
\begin{align}
 -\inf_{\eta\in\widetilde G_q(r)}\mathcal Y(\eta)
 &\leq
 \liminf_{n\to\infty}\frac1{\beta n^2}\log\mathbb P(S_n>r)
 \leq
 \limsup_{n\to\infty}\frac1{\beta n^2}\log\mathbb P(S_n>r)
 \leq
 -\inf_{\eta\in\widehat G_q(r)}\mathcal Y(\eta).
 \label{eq:svtldp}
\end{align}
Moreover,
\begin{equation}\label{eq:svinfs}
 \inf_{\eta\in\widetilde G_q(r)}\mathcal Y(\eta)
 =
 \inf_{\eta\in\widehat G_q(r)}\mathcal Y(\eta).
\end{equation}
Indeed, this follows from the argument used for
\eqref{eq:sainfs}, with $U_r$ replaced by
$\rho_r$, the uniform probability measure on
$[((1+r)/2)^2,((3+r)/4)^2]$. Its $(q/2)$-moment is greater than
$r^q$, and it belongs to $\widehat{\mathcal D}_p$ in both cases. If
$\eta$ is a boundary point with $\mathcal Y(\eta)<\infty$, then
$\eta_\varepsilon=(1-\varepsilon)\eta+\varepsilon\rho_r$ belongs to
$\widetilde G_q(r)$ and
$\mathcal Y(\eta_\varepsilon)\to\mathcal Y(\eta)$ as
$\varepsilon\downarrow0$.

To simplify the variational problem \eqref{eq:svinfs}, suppose first that $p<\infty$. Fix
$\eta\in\mathcal M_1([0,\infty))$ with $m_{p/2}(\eta)\leq1$,
$m_{q/2}(\eta)\geq r^q$, and $\Sigma(\eta)>-\infty$. Then
$m_{p/2}(\eta)>0$. Let $a=m_{p/2}(\eta)^{-2/p}\geq1$ and
$\widetilde\eta=(x\mapsto ax)_\#\eta$. Scaling yields
$m_{p/2}(\widetilde\eta)=1$,
$m_{q/2}(\widetilde\eta)=a^{q/2}m_{q/2}(\eta)\geq r^q$, and
$\Sigma(\widetilde\eta)=\Sigma(\eta)+\log a\geq\Sigma(\eta)$.
If $m_{q/2}(\widetilde\eta)>r^q$, let
$$
\eta_t=(1-t)\widetilde\eta+t\eta^{(p)}, \quad
t=\frac{m_{q/2}(\widetilde\eta)-r^q}{m_{q/2}(\widetilde\eta)-r_{p,q}^q)}\in(0,1).
$$
Remark \ref{rem:squlm} gives
$m_{p/2}(\eta_t)=1$ and $m_{q/2}(\eta_t)=r^q$.
The maximality stated there and concavity imply that
$\Sigma(\eta_t)\geq\Sigma(\widetilde\eta)\geq\Sigma(\eta)$.
If $m_{q/2}(\widetilde\eta)=r^q$, simply take
$\eta_t=\widetilde\eta$. Consequently, for $p<\infty$,
\begin{equation}\label{eq:gcsb}
\inf_{\eta\in\widehat G_q(r)}\mathcal Y(\eta)
=2C_p-\frac12\sup\left\{
\Sigma(\eta):
\eta\in\mathcal M_1([0,\infty)),\
m_{p/2}(\eta)=1, m_{q/2}(\eta)=r^q
\right\}.
\end{equation}

For $p=\infty$, let
$$
\eta^{(\infty)}(\dd y)
=\frac{\mathbf 1_{(0,1)}(y)}{\pi\sqrt{y(1-y)}} \dd y.
$$
Equivalently,
$\eta^{(\infty)}=(x\mapsto x^2)_\#\nu^{(\infty)}$, where
$\nu^{(\infty)}$ is defined in \eqref{eq:arclaw}. This
measure maximizes $\Sigma$ over probability measures supported on
$[0,1]$ and satisfies
$m_{q/2}(\eta^{(\infty)})=r_{\infty,q}^q<r^q$. Then the same argument as above yields \eqref{eq:gcsb} for $p=\8$ with $m_{p/2}(\eta)=1$ replaced by $\eta\in \mx_1([0,1])$.

Let us determine the relationship between the rate functions in the self-adjoint and general cases. For $p<\infty$, let
\begin{align*}
A&=
\sup\left\{
\Sigma(\eta):
\eta\in\mathcal M_1([0,\infty)),\
m_{p/2}(\eta)=1, m_{q/2}(\eta)=r^q
\right\},
\\
B&=\sup\left\{
\Sigma(\mu):
\mu\in\mathcal M_1(\mathbb R),\
m_p(\mu)=1, m_q(\mu)=r^q
\right\}.
\end{align*}
We claim that $A=2B$. First let $\eta$ be
admissible for $A$, and let $\mu_\eta$ be the law of
$\sigma_{\mathrm{sgn}}\sqrt X$, where $X$ has law $\eta$ and
$\sigma_{\mathrm{sgn}}$ is independent of $X$ with
$\mathbb P(\sigma_{\mathrm{sgn}}=1)
=\mathbb P(\sigma_{\mathrm{sgn}}=-1)=1/2$.
For every $s>0$, $m_s(\mu_\eta)=m_{s/2}(\eta)$, so $\mu_\eta$ is admissible for $B$. 
Averaging over the four sign combinations gives
\begin{align*}
\Sigma(\mu_\eta)
=\frac12\iint_{[0,\infty)^2}
\bigl(\log|\sqrt{x}-\sqrt{y}|+\log(\sqrt{x}+\sqrt{y})\bigr)
 \eta(\dd x)\eta(\dd y)
=\frac12\Sigma(\eta).
\end{align*}
Taking the supremum over all admissible $\eta$ gives $A\leq2B$.

For the reverse inequality, let $\mu$ be admissible for $B$, and
set $\check\mu=(x\mapsto-x)_\#\mu$ and
$\overline\mu=(\mu+\check\mu)/2$. Symmetrization preserves all
absolute moments, while the concavity of $\Sigma$ yields
$\Sigma(\overline\mu)\geq\Sigma(\mu)$. Define
$\eta_\mu=(x\mapsto x^2)_\#\overline\mu$. Then
$m_{s/2}(\eta_\mu)=m_s(\overline\mu)=m_s(\mu)$ for every $s>0$.
The symmetry of $\overline\mu$ also gives
\begin{align*}
\Sigma(\eta_\mu)
&=\iint_{\mathbb R^2}\log|x^2-y^2| 
\overline\mu(\dd x)\overline\mu(\dd y)\\
&=\iint_{\mathbb R^2}
\bigl(\log|x-y|+\log|x+y|\bigr)
 \overline\mu(\dd x)\overline\mu(\dd y)
=2\Sigma(\overline\mu)
\geq2\Sigma(\mu).
\end{align*}
Thus $\eta_\mu$ is admissible for $A$, and
$A\geq\Sigma(\eta_\mu)\geq2\Sigma(\mu)$. Taking the supremum over all admissible $\mu$ implies that $A\geq2B$, and hence $A=2B$. With the corresponding support constraints, the same argument gives for
$p=\infty$,
\begin{align*}
&
\sup\left\{
\Sigma(\eta):
\eta\in\mathcal M_1([0,\infty)),\
\operatorname{supp}(\eta)\subseteq[0,1],\
m_{q/2}(\eta)=r^q
\right\}
\\
={}&2\sup\left\{
\Sigma(\mu):
\mu\in\mathcal M_1(\mathbb R),\
\operatorname{supp}(\mu)\subseteq[-1,1],\
m_q(\mu)=r^q
\right\}.
\end{align*}
Together with \eqref{eq:gcsb}, these variational identities and the definition \eqref{eq:phieq} lead to $\inf_{\eta\in\widehat G_q(r)}\mathcal Y(\eta) =2\Phi_{p,q}(r).$ 
Combining this identity with \eqref{eq:svinfs}, \eqref{eq:svtldp} and \eqref{eq:svvolt}, we have proved part (ii).

Finally, let us study the properties of the rate function $\Phi_{p,q}$. For $0<q<p\leq\infty$, define
$$
H(s)=\Phi_{p,q}(s^{1/q})=\inf_{\substack{\mu\in\mathcal D_p\\m_q(\mu)\geq s}}
\mathcal K(\mu),\qquad r_{p,q}^q\leq s<1.
$$
Since $m_q(\nu^{(p)})=r_{p,q}^q$ and $\mathcal K(\nu^{(p)})=0$, we have $H(r_{p,q}^q)=0$.
For $r_{p,q}^q<s<1$, the measure $U_{s^{1/q}}$ in \eqref{eq:Ur} belongs to $\mathcal D_p$ and satisfies $ m_q(U_{s^{1/q}})>s$.  
Since $U_{s^{1/q}}$ has bounded density with compact support, $\Sigma(U_{s^{1/q}})$ is finite. It follows that
$\mathcal K(U_{s^{1/q}})<\infty$.
As a good rate function, $\mathcal K$ attains the minimum on the closed set $\{\mu\in\mathcal D_p:m_q(\mu)\geq s\}$. But the unique zero of $\kx$ is $\nu^{(p)}$ and $m_q(\nu^{(p)})=r_{p,q}^q<s$. We find
$$
0<H(s)\leq\mathcal K(U_{s^{1/q}})<\infty.
$$
Fix $0<t<1$.
For $\mu_i\in\mathcal D_p$ with $m_q(\mu_i)\geq s_i \in[r_{p,q}^q,1) $, $i=1,2$, let $\mu_t=(1-t)\mu_1+t\mu_2$.
By convexity, $\mu_t\in\mathcal D_p$.
Using the linearity of $m_q$ and the concavity of $\Sigma$, we have
\begin{align*}
m_q(\mu_t)&=(1-t)m_q(\mu_1)+tm_q(\mu_2)
\geq(1-t)s_1+ts_2,\\
\mathcal K(\mu_t)&\leq(1-t)\mathcal K(\mu_1)+t\mathcal K(\mu_2).
\end{align*}
Taking the infimum on both sides, we obtain
$$
H((1-t)s_1+ts_2)\leq(1-t)H(s_1)+tH(s_2).
$$
It follows that $H$ is finite and convex, thus continuous on $(r_{p,q}^q,1)$.
Since $H(r_{p,q}^q)=0$ and $H(s)>0$ for $r_{p,q}^q<s<1$,
convexity also implies that $H(s)$ is strictly increasing.
The same properties hold for $\Phi_{p,q}(r)=H(r^q)$ on
$(r_{p,q},1)$. 
\end{proof}

\begin{remark}\label{re:speedd}
One may wonder what would happen if we applied \eqref{eq:savolt} to the case $p<q$. Since $\mu\mapsto m_q(\mu)$ is lower semicontinuous, we still know that $\widetilde F_q(r)$ is open and the lower bound of \eqref{eq:satldp} is valid. However, it is not hard to construct a family $\nu_n\in\widetilde F_q(r)$ such that $\lim_{n\to\8}\kx(\nu_n)= 0$ for $r>r_{p,q}$. This means that the correct speed should be slower than $n^2$, which is consistent with the speed $n^{1+p/q}$ established in Proposition \ref{thm:diffball}. A similar discussion applies to the general case.
\end{remark}

\subsection{Proofs of technical results}

\begin{proof}[Proof of Lemma \ref{le:pc-Freud}]
For $0\leq k<n$ and $\bx\in\mathbb R^k$, write
$$
I_k^{(n)}(\bx)=\int_{\mathbb R}e^{-\beta nV(t)}
\prod_{i=1}^k|t-x_i|^\beta\dd t.
$$
Restricting the integral defining $I_k^{(n)}$ to $[-1,1]$ and
applying Jensen's inequality with respect to the probability measure
$\frac12\mathbf1_{[-1,1]}(t)\dd t$, we obtain
$$
I_k^{(n)}(\bx)
\geq2\exp\left\{-\beta n\sup_{[-1,1]}V-\beta k\right\}
\geq e^{-C\beta n},
$$
where we used $\inf_{x\in\mathbb R}\frac12\int_{-1}^1\log|t-x|\dd t=-1$.
Consequently,
\begin{equation}\label{eq:Fs}
Z_{k+1}^{(n)}
=\int_{\mathbb R^k}e^{-\beta n\sum_{i=1}^kV(x_i)}
\prod_{1\leq i<j\leq k}|x_i-x_j|^\beta
I_k^{(n)}(\bx)\dd\bx
\geq e^{-C\beta n}Z_k^{(n)}.
\end{equation}
Iterating \eqref{eq:Fs} proves \eqref{eq:Zratio}.

To prove \eqref{eq:widezz}, we use \eqref{e:pairbd} and the power growth
of $V$ to obtain, for some $c,C>0$,
\begin{align*}
I_{n-1}^{(n)}(\bx)
&\leq\left(\int_{\mathbb R}e^{-\beta nV(t)}
(1+|t|)^{\beta(n-1)}\dd t\right)
\prod_{i=1}^{n-1}(1+|x_i|)^\beta\\
&\leq e^{C\beta n}
\left(\int_{\mathbb R}e^{-c\beta n|t|^p}\dd t\right)
\prod_{i=1}^{n-1}(1+|x_i|)^\beta
\leq e^{C'\beta n}\prod_{i=1}^{n-1}(1+|x_i|)^\beta,
\end{align*}
where in the last inequality we used
$\int_{\mathbb R}e^{-c\beta n|t|^p}\dd t
=2\Gamma(1+1/p)(c\beta n)^{-1/p}$ and the fact that $\beta>0$ is fixed.
Substituting into \eqref{eq:Fs} and using
$\log(1+|x|)\leq V(x)+C$ give
\begin{align*}
Z_n^{(n)}
&\leq e^{C\beta n}\int_{\mathbb R^{n-1}}
e^{-\beta n\sum_{i=1}^{n-1}V(x_i)}
\prod_{1\leq i<j\leq n-1}|x_i-x_j|^\beta
\prod_{i=1}^{n-1}(1+|x_i|)^\beta\dd\bx\\
&\leq e^{C'\beta n}\int_{\mathbb R^{n-1}}
e^{-\beta(n-1)\sum_{i=1}^{n-1}V(x_i)}
\prod_{1\leq i<j\leq n-1}|x_i-x_j|^\beta\dd\bx
=e^{C'\beta n}Z_{n-1}^{(n-1)}.
\end{align*}
Next, we compare $Z_{n-1}^{(n-1)}$ with $Z_{n-1}^{(n)}$.
Iterating \eqref{eq:Fs} from $Z_0^{(n)}=1$ yields
$$
Z_{n-1}^{(n)}\geq e^{-C\beta n(n-1)},
\qquad \log Z_{n-1}^{(n)}\geq-C\beta n^2.
$$
By \eqref{e:pairbd} and the power growth of $V$, we obtain
\begin{align*}
Z_{n-1}^{(n/2)}
&\leq\left(\int_{\mathbb R}
e^{-\beta nV(x)/2}(1+|x|)^{\beta(n-2)}\dd x\right)^{n-1}
\leq\left(e^{C\beta n}
\int_{\mathbb R}e^{-c\beta n|x|^p}\dd x\right)^{n-1}
\leq e^{C'\beta n^2}.
\end{align*}
For $n>2$, the identity
$
n-1=\left(1-\frac2n\right)n+\frac2n\frac n2
$
and H\"older's inequality imply
$$
Z_{n-1}^{(n-1)}
\leq\left(Z_{n-1}^{(n)}\right)^{1-2/n}
\left(Z_{n-1}^{(n/2)}\right)^{2/n}.
$$
Consequently,
\begin{align*}
\frac{Z_n^{(n)}}{Z_{n-1}^{(n)}}
\leq e^{C\beta n}
\frac{Z_{n-1}^{(n-1)}}{Z_{n-1}^{(n)}}
\leq e^{C\beta n}
\left(\frac{Z_{n-1}^{(n/2)}}{Z_{n-1}^{(n)}}\right)^{2/n}
\leq e^{C'\beta n},
\end{align*}
which proves \eqref{eq:widezz} for large $n$.
\end{proof}

\begin{proof}[Proof of Theorem \ref{th:momldp}]
For \eqref{pdf} with $D=\mathbb R$ and $V_n=V$,
Assumption A in \cite{WZbeta} follows from the continuity and power growth of $V$, and \cite{WZbeta}*{Lemma 2.3} verifies Assumption B.
Thus \cite{WZbeta}*{Theorem 1.1} gives the LDP of empirical measure at speed $\beta n^2$, with unique minimizer $\mu_V$.

Applying \cite{WZbeta}*{Eq.~A.17)} to both the particles and their
reflections, we obtain
\begin{equation}\label{eq:abstail}
\mathbb P\left(\max_{i\leq n}|\lambda_i|>R\right)
\leq ne^{C\beta n}\int_{|x|>R}
e^{-\beta nV(x)}(1+|x|)^{\beta(n-1)}\dd x
\leq e^{-d_R\beta n}
\end{equation}
for every sufficiently large fixed $R$ and all sufficiently large $n$, where $d_R\to\infty$ as $R\to\infty$. Choose such an $R>\max_{x\in\supp\mu_V}|x|$
and let  $f_R(x)=|x|^q\wedge R^q$. Write $M_{q,n}=\frac1n\sum_{i=1}^n |\la_i|^q$ in this proof.
Since the rate function $J_q(t)$ is clearly continuous and strictly increasing and for $t<m_q$, the LDP of empirical measure and
$M_{q,n}\geq n^{-1}\sum_i f_R(\lambda_i)$ imply
$$
\limsup_{n\to\infty}\frac1{\beta n^{1+p/q}}
\log\mathbb P(M_{q,n}\leq t)=-\infty,
$$
It suffices to show that 
\begin{equation}\label{eq:upperlim}
\lim_{n\to\infty}\frac1{\beta n^{1+p/q}}
\log\mathbb P(M_{q,n}\geq t)
= -b(t-m_q)^{p/q},\qquad t>m_q.
\end{equation}

We first prove the local lower bound
\begin{align}
\lim_{\delta\downarrow0}\liminf_{n\to\infty}
\frac1{\beta n^{1+p/q}}
\log\mathbb P(|M_{q,n}-t|<\delta)
&\geq-J_q(t), \qquad t  > m_q . \label{eq:mqlb}
\end{align}
Once \eqref{eq:mqlb} is established, fix $\eta>0$
and take $0<\delta<\eta$. The inclusion
$\{|M_{q,n}-(t+\eta)|<\delta\}\subset\{M_{q,n}\geq t\}$ 
gives
\begin{align*}
\liminf_{n\to\infty}\frac1{\beta n^{1+p/q}}
\log\mathbb P(M_{q,n}\geq t)
\geq\lim_{\delta\downarrow0}\liminf_{n\to\infty}
\frac1{\beta n^{1+p/q}}
\log\mathbb P(|M_{q,n}-(t+\eta)|<\delta)
\geq-J_q(t+\eta).
\end{align*}
Letting $\eta\downarrow0$, the lower bound in \eqref{eq:upperlim} follows from the continuity of $J_q$.
To show \eqref{eq:mqlb}, we fix  $\delta>0$,
and a sufficiently large $R$ so that
$\supp\mu_V\subset (-R,R)$.
Consider the region
$$
E_n=\left\{
\by=(y_1,\ldots,y_{n-1})\in[-R,R]^{n-1}:
\left|\frac1{n-1}\sum_{i=1}^{n-1}|y_i|^q-m_q\right|<\frac\delta4
\right\}.
$$
We claim $
\mathbb P\bigl((\lambda_1,\ldots,\lambda_{n-1})\in E_n\bigr)
\to1.
$
Indeed, since $f_R(x)$ is bounded and continuous,  the LDP of empirical measure implies that
$\int f_R\dd L_n\xrightarrow{\mathbb P}\int f_R\dd\mu_V=m_q.$
Note that
$$
\left|\frac1{n-1}\sum_{i<n}f_R(\lambda_i)-\int f_R\dd L_n\right|
\leq\frac{2R^q}{n},
$$
and on $\{\max_{i<n}|\lambda_i|\leq R\}$,
$$
\frac1{n-1}\sum_{i<n}f_R(\lambda_i)
=\frac1{n-1}\sum_{i<n}|\lambda_i|^q.
$$
Choosing $n$ sufficiently large so that $2R^q/n<\delta/8$
and using the union bound and \eqref{eq:abstail}, we deduce that
\begin{align*}
\mathbb P\bigl((\lambda_1,\ldots,\lambda_{n-1})\notin E_n\bigr)
\leq\mathbb P\left(\max_{i<n}|\lambda_i|>R\right)
+\mathbb P\left(\left|\int f_R\dd L_n-m_q\right|\geq\frac\delta8\right)
\to0.
\end{align*}
For $\by\in E_n$, by the power growth of $V$, 
\begin{align*}
& \int_{\mathbb R}e^{-\beta nV(z)}
\prod_{i<n}|z-y_i|^\beta\dd z
\leq\int_{\mathbb R}e^{-\beta nV(z)}(|z|+R)^{\beta(n-1)}\dd z\\
\leq{} &(1+R)^{\beta(n-1)}\int_{\mathbb R}
e^{-\beta n[V(z)-\log(1+|z|)]}\dd z
\leq e^{\beta n[C+\log(1+R)]}
\int_{\mathbb R}e^{-c\beta n|z|^p}\dd z
\leq e^{C_R\beta n}.
\end{align*}
Restrict the last particle to
$I_n=\left[(n(t-m_q))^{1/q},\ (n(t-m_q))^{1/q}+1\right].$
Uniformly for $\by\in E_n$ and $z\in I_n$, we have for large $n$,
$$
|M_{q,n}-t|<\delta,\qquad
|z-y_i|\geq1,\qquad
\beta nV(z)=\beta bn^{1+p/q}(t-m_q)^{p/q}+o(\beta n^{1+p/q}),
$$
and thus using conditional density $p(z|\by)=p(\by,z)/p(\by)$ yields
\begin{align*}
\mathbb P(|M_{q,n}-t|<\delta\mid \by )
&\geq e^{-C_R\beta n}\int_{I_n}e^{-\beta nV(z)} \prod_{i<n}|z-y_i|^\beta \dd z\\
&\geq\exp\left\{-\beta bn^{1+p/q}(t-m_q)^{p/q}
+o(\beta n^{1+p/q})\right\}.
\end{align*}
Integrating this uniform bound over $E_n$ with respect to the marginal law of $(\lambda_1,\ldots,\lambda_{n-1})$, we find
\begin{align*}
\mathbb P(|M_{q,n}-t|<\delta)
&\geq\mathbb P\bigl((\lambda_1,\ldots,\lambda_{n-1})\in E_n\bigr)\exp\left\{-\beta bn^{1+p/q}(t-m_q)^{p/q}
+o(\beta n^{1+p/q})\right\}.
\end{align*}
Since the first factor tends to one,  \eqref{eq:mqlb}  follows from taking logarithms and dividing by $\beta n^{1+p/q}$.

Next, we prove the upper bound of \eqref{eq:upperlim}.
Fix $t>m_q$ and $0<\varepsilon<t-m_q$.
For $R>\max_{x\in\supp\mu_V}|x|$, applying the LDP of empirical measure to the function $f_R$, we know 
\begin{equation}\label{eq:frlatail}
    \limsup_{n\to\infty}\frac1{\beta n^{1+p/q}}
\log\mathbb P\left(\frac1n\sum_{i=1}^nf_R(\lambda_i)>m_q+\varepsilon\right)
=-\infty.
\end{equation}
It remains to bound the probability of
$$
\mathcal A_{n,R,\varepsilon}
=\left\{M_{q,n}\geq t,\quad
\frac1n\sum_{i=1}^nf_R(\lambda_i)\leq m_q+\varepsilon\right\}.
$$
We decompose this event according to the number $m$ of coordinates with modulus
greater than $R$. Their values belong to
$$
\Gamma_{n,m}=\left\{\bz\in\mathbb R^m:
|z_i|>R,\quad
\sum_{i=1}^m(|z_i|^q-R^q)\geq n(t-m_q-\varepsilon)\right\}.
$$
By the joint density,
\begin{align}\label{eq:anrv}
\mathbb P(\mathcal A_{n,R,\varepsilon})
\leq\sum_{m=1}^n\binom nm
\int_{\Gamma_{n,m}}e^{-\beta n\sum_{i=1}^mV(z_i)}
\prod_{1\leq i<j\leq m}|z_i-z_j|^\beta
\mathcal B_{n,m}(\bz)\dd\bz,
\end{align}
where
$$
\mathcal B_{n,m}(\bz)=\frac1{Z_n^{(n)}}
\int_{[-R,R]^{n-m}}e^{-\beta n\sum_{j=1}^{n-m}V(u_j)}
\prod_{1\leq j<k\leq n-m}|u_j-u_k|^\beta
\prod_{i=1}^m\prod_{j=1}^{n-m}|z_i-u_j|^\beta\dd\bu.
$$
The idea is to treat the interactions between the large coordinates and between large and small coordinates separately.  For $R\geq1$, let
$$
h_R=\sup_{u\geq R}\frac{\log(2u)}{u^p},
\qquad
b_R=\inf_{|x|\geq R}\frac{V(x)}{|x|^p}.
$$
For $|z|,|z'|>R$ and $|u|\leq R$,
$$
\log|z-u|\leq h_R|z|^p,\qquad
\log|z-z'|\leq h_R(|z|^p+|z'|^p).
$$
Together with Lemma \ref{le:pc-Freud}, we derive
\begin{align}\label{eq:sabn}
\mathcal B_{n,m}(\bz)
\leq\frac{Z_{n-m}^{(n)}}{Z_n^{(n)}}
\exp\left\{\beta(n-m)h_R\sum_{i=1}^m|z_i|^p\right\}\leq\exp\left\{C\beta nm+\beta(n-m)h_R\sum_{i=1}^m|z_i|^p\right\}.
\end{align}
The interaction among the $m$ large coordinates is at most
$\exp\{\beta(m-1)h_R\sum_i|z_i|^p\}$.
Also, on $\Gamma_{n,m}$, subadditivity of $u^{p/q}$ gives
\begin{equation}\label{eq:I-p-cost}
\begin{aligned}
\sum_{i=1}^m|z_i|^p
\geq\left(\sum_{i=1}^m(|z_i|^q-R^q)\right)^{p/q}
\geq n^{p/q}(t-m_q-\varepsilon)^{p/q}, \qquad
\sum_{i=1}^m|z_i|^p \geq mR^p.
\end{aligned}
\end{equation}
Using this bound and writing $c_R=b_R-h_R-C/R^p$, we find
\begin{align*}
\mathbb P(\mathcal A_{n,R,\varepsilon})
\leq\sum_{m=1}^n\binom nm
\int_{\Gamma_{n,m}}e^{-\beta nc_R\sum_{i=1}^m|z_i|^p}\dd\bz
\leq 2^n
e^{-\beta(c_R-\kappa)n^{1+p/q}(t-m_q-\varepsilon)^{p/q}}
\end{align*}
for fixed $0<\kappa<c_R$ and all sufficiently large $n$, where we used the fact that
$$
\int_{\mathbb R^m}e^{-\beta n\kappa\sum_{i=1}^m|z_i|^p}\dd\bz
=\left[\frac2p\Gamma(1/p)(\beta n\kappa)^{-1/p}\right]^m\leq1.
$$
Letting $n\to\infty$ and combining \eqref{eq:frlatail}, we arrive at
\[
\limsup_{n\to\8} \frac1{\beta n^{1+p/q}}
\log\mathbb P(M_{q,n}\geq t)\le -  (c_R-\kappa) (t-m_q-\eps)^{p/q}.
\]
Noting that $h_R\to0$ and $b_R\to b$ as $R\to\infty$, we have $c_R\to b$. Therefore, for fixed $0<\kappa<b$, the condition $\kappa<c_R$ holds for all sufficiently large $R$.
Sending first $R\to\infty$, and then $\kappa\downarrow0$ and $\varepsilon\downarrow0$, we obtain the upper bound of \eqref{eq:upperlim}.
\end{proof}

\begin{proof}[Proof of Lemma \ref{le:snoutlier}]
Write $V_p(x)=|x|^p/(2p)$. First suppose that $p<q<\infty$. Then Theorem \ref{th:momldp} holds with $b=1/(2p)$, and Remark \ref{rem:ullman} gives the equilibrium measure $\mu_{V_p}=\nu^{(p)}$ and $m_q=r_{p,q}^q$.
Since $B_{q,n}=M_{q,n}^{1/q}$ and
$J_q(s^q)=\frac1{2p}(s^q-r_{p,q}^q)^{p/q}=\Phi_{p,q}(s)$
by \eqref{eq:phipq}, \eqref{eq:upperlim} and the continuity of $J_q$
yield \eqref{eq:oneoutq} for $p<q<\infty$.

It remains to prove \eqref{eq:oneoutq} for $q=\infty$. By Remark \ref{rem:ullman},   $\max\supp\nu^{(p)}=b_p$.
For $r>b_p$, define
$$J_{V_p}^H(r)
=
-\int \log|r-y| \nu^{(p)}(\dd y)
+
V_p(r)
-
\ell_{V_p},$$
where
$
\ell_{V_p}
=
-\int \log(b_p-y) \nu^{(p)}(\dd y)
+
\frac{b_p^p}{2p}.
$
For $r>b_p$, the logarithmic potential formula in \cite{ST97}*{Section IV.5, Eq.~(5.2)} gives 
\begin{align}\label{e:powlog}
& -\int\log(r-y) \nu^{(p)}(\dd y)
=-\log b_p+\log2+\frac1p
-\log\frac{r+\sqrt{r^2-b_p^2}}{b_p}
-\frac{r}{b_p^p}\int_0^{b_p}
\frac{x^{p-1}}{\sqrt{r^2-x^2}}\dd x\notag\\
={} & -\log b_p+\log2+\frac1p
-\log\frac{r+\sqrt{r^2-b_p^2}}{b_p}
-\frac1p {}_2F_1\left(
\frac12,\frac p2;\frac p2+1;\frac{b_p^2}{r^2}\right),
\end{align}
where the last equality follows from \eqref{e:hypint} after the substitution $z=x^2/b_p^2$.
By continuity, taking $r=b_p$ in the same formula yields
$
\ell_{V_p}
=
-\log b_p+\log 2+\frac1p.
$
Therefore, 
$$J_{V_p}^H(r)=
\frac{r^p}{2p}
-
\frac{r}{b_p^p}
\int_0^{b_p}
\frac{x^{p-1}}{\sqrt{r^2-x^2}} \dd x
-
\log\frac{r+\sqrt{r^2-b_p^2}}{b_p}=\Phi_{p,\infty}(r),$$
where $\Phi_{p,\infty}$ is defined in \eqref{eq:phipq} after a change of variable. 
Moreover, from \cite{ST97}*{Section IV.5}, we know that the function $J_{V_p}^H$ is continuous and strictly increasing on $(b_p,\infty)$. Hence
$$
\inf_{t\ge r}J_{V_p}^H(t)
=
J_{V_p}^H(r),
\qquad r>b_p.
$$
For every $r>b_p$,  by \cite{BG13}*{Proposition 2.1}, we have
$$
\lim_{n\to\infty}
\frac1{\beta n}
\log \mathbb P(\lambda_{\max}\ge r)
=-J_{V_p}^H(r)
=
-\Phi_{p,\infty}(r).
$$
The same limit holds with $>$ in place of $\ge$ by continuity. 
Finally, since $V_p$ is even, the joint density is invariant under $(\lambda_1,\ldots,\lambda_n)\mapsto (-\lambda_1,\ldots,-\lambda_n).$
Hence $-\lambda_{\min}$ has the same distribution as $\lambda_{\max}$. It follows that $
\{B_{\infty,n}>r\}
=
\{\lambda_{\max}>r\}\cup\{-\lambda_{\min}>r\},
$
and
$$
\mathbb P(\lambda_{\max}>r)
\le
\mathbb P(B_{\infty,n}>r)
\le
2\mathbb P(\lambda_{\max}>r).
$$
Therefore,
$\lim_{n\to\infty}
\frac1{\beta n}
\log
\mathbb P(B_{\infty,n}>r)
=
-\Phi_{p,\infty}(r)$ for $r>b_p.$
This completes the proof.
\end{proof}

\begin{proof}[Proof of Lemma \ref{le:sv-pc}]
The argument is similar to that of Lemma \ref{le:pc-Freud}. For $0\le k\le n-1$, we rewrite
$$    
\bar Z_{k+1,p,\beta}^{(n)}
=
\int_{[0,\infty)^k}
\prod_{1\le i<j\le k}|x_i-x_j|^\beta
\prod_{i=1}^k x_i^{\frac{\beta}{2}-1}
\exp\left\{
-\frac{\beta n}{p}\sum_{i=1}^k x_i^{p/2}
\right\}
J_k^{(n)}(\bx) \dd\bx,
$$
where
$
J_k^{(n)}(\bx)
=
\int_0^\infty
\exp\left\{-\frac{\beta n}{p}t^{p/2}\right\}
t^{\frac{\beta}{2}-1}
\prod_{i=1}^k |t-x_i|^\beta \dd t
$
for $\bx=(x_1,\ldots,x_k)\in[0,\infty)^k$. Since the function $x\mapsto \int_1^2\log|t-x| \dd t$
is bounded from below on $[0,\infty)$, restricting $J_k^{(n)}$ to $t\in[1,2]$ and applying Jensen's inequality yield
$$J_k^{(n)}(\bx)
\ge
e^{-C\beta n}
e^{\beta\sum_{i=1}^k\int_1^2\log|t-x_i| \dd t}
\ge e^{-C\beta n-C\beta k}.$$
It follows that $
\bar Z_{k+1,p,\beta}^{(n)}
\ge
\exp\{-C\beta n-C\beta k\}
\bar Z_{k,p,\beta}^{(n)}.
$
Iterating this estimate gives the first assertion.

To prove \eqref{eq:barzz}, we apply \eqref{e:pairbd} and obtain
\begin{align*}
\bar Z_{n,p,\beta}^{(n)}
\le{}&
J_n
\int_{[0,\infty)^{n-1}}
\prod_{1\le i<j\le n-1}|y_i-y_j|^\beta
\prod_{i=1}^{n-1}y_i^{\frac{\beta}{2}-1}
\exp\left\{
-\frac{\beta n}{p}\sum_{i=1}^{n-1}y_i^{p/2}
\right\}
\prod_{i=1}^{n-1}(1+y_i)^\beta
 \dd\by,
\end{align*}
where $J_n=\int_0^\infty
z^{\frac{\beta}{2}-1}
\exp\left\{
-\frac{\beta n}{p}z^{p/2}
\right\}
(1+z)^{\beta(n-1)}
 \dd z.
$
 Applying \eqref{eq:Frelb} with $\alpha=p/2$ and $\varepsilon=1/(2p)$, we obtain
$\log(1+x)
\le
\frac{1}{2p}x^{p/2}+C_p.$ Thus,
$$
\prod_{i=1}^{n-1}(1+y_i)^\beta
\exp\left\{
-\frac{\beta n}{p}\sum_{i=1}^{n-1}y_i^{p/2}
\right\}
\le
\exp\{C_{p,\beta}n\}
\exp\left\{
-\frac{\beta(n-\frac12)}{p}
\sum_{i=1}^{n-1}y_i^{p/2}
\right\}.
$$
The same bound gives $J_n\le e^{C_{p,\beta}n}$.
It follows that
$
\bar Z_{n,p,\beta}^{(n)}
\le e^{C_{p,\beta}n} \bar Z_{n-1,p,\beta}^{(n-\frac12)}.
$
Under the change of variable  $y_i=\left(\frac{n}{n-1/2}\right)^{2 / p} x_i,$ we have
\begin{align*}
\bar Z_{n-1,p,\beta}^{(n-1/2)}
&=
\left(
\frac{n}{n-1/2}
\right)^{\beta(n-1)^2/p}
\bar Z_{n-1,p,\beta}^{(n)}\le
\exp\{C_{p,\beta}n\}
\bar Z_{n-1,p,\beta}^{(n)}.
\end{align*}
It follows that $\bar Z_{n,p,\beta}^{(n)}
\le e^{C_{p,\beta}n}
\bar Z_{n-1,p,\beta}^{(n)}
$
and the proof is complete.
\end{proof}

\begin{proof}[Proof of Lemma \ref{le:sv-o}] 
The current setting has a new subtlety that the exponent $\frac\bt2 -1$ may be negative and is analogous to Theorem \ref{th:momldp} and Lemma \ref{le:snoutlier} otherwise.

\noindent\emph{(i) Case $q=\8$.} 
The vector $\bY$ has the distribution in Lemma \ref{x:powemp}. 
As in the proof of Lemma \ref{le:snoutlier} for $q=\infty$, the LDP in Lemma  \ref{x:powmax} together with the continuity and strict monotonicity of $\Phi_{p,\infty}$ yields
$$\lim_{n\to\infty}\frac1{\beta n}
\log\mathbb P(\tilde B_{\infty,n}>r^2)
=-2\Phi_{p,\infty}(r),\qquad r>b_p,$$
which proves the claim thanks to the limits in \eqref{eq:brslim}. 

\noindent\emph{(ii) Case $q<\infty$.} 
We start with the lower bound. Fix $\varepsilon>0$, and let 
$$
I_{n,\varepsilon}
=
\left[
n^{2/q}(s^q-r_{p,q}^q+\varepsilon)^{2/q},
n^{2/q}(s^q-r_{p,q}^q+\varepsilon)^{2/q}+1
\right].
$$
Let $
f_R(x)=x^{q/2}\wedge R, x\ge0. $ Choose $R>b_p^q$. 
Remark \ref{rem:squlm} gives
$\int f_R \dd\eta^{(p)}= \int x^{q/2} \eta^{(p)}(\dd x) =r_{p,q}^q.$
Choose $\eta_R>0$ such that $\nu\in B(\eta^{(p)},\eta_R)$  implies $
|\int f_R \dd\nu-\int f_R \dd\eta^{(p)}|<\frac{\varepsilon}{4}.
$
Let $\bar{\mathbb P}_{n-1}^{(n)}$ denote the $(n-1)$-particle probability measure with density proportional to
$$
\prod_{1\le i<j\le n-1}|y_i-y_j|^\beta
 \prod_{i=1}^{n-1} y_i^{\frac{\beta}{2}-1} \nc 
\exp\left\{
-\frac{\beta n}{p}\sum_{i=1}^{n-1}y_i^{p/2}
\right\}
\mathbf 1_{\{y_i\ge0\}} \dd\by.
$$
Let $E_{n-1}= \{L_{n-1}\in B(\eta^{(p)},\eta_R),\
\max_{1\le i\le n-1}Y_i\le R^{2/q} \}.$ 
By Lemmas \ref{x:powemp} and \ref{x:powmax}, $\bar{\mathbb P}_{n-1}^{(n)}(E_{n-1})\to1.$
For every configuration in $E_{n-1}$, we have 
$
\frac1{n-1}\sum_{i=1}^{n-1}Y_i^{q/2}
\ge
r_{p,q}^q-\frac{\varepsilon}{4}.
$
Hence, on $E_{n-1}$ and for $z\in I_{n,\varepsilon}$,
$$
\frac1n
\left(
\sum_{i=1}^{n-1}Y_i^{q/2}+z^{q/2}
\right)
\ge
\frac{n-1}{n}\left(r_{p,q}^q-\frac{\varepsilon}{4}\right)
+s^q-r_{p,q}^q+\varepsilon
>
s^q
$$
for all sufficiently large $n$. 
Separating the last particle in the joint density, we find
\begin{align*}
\mathbb P(\tilde B_{q,n}>s^2)
\ge{}&
n \frac{\bar Z_{n-1,p,\beta}^{(n)}}{\bar Z_{n,p,\beta}^{(n)}}
\int_{I_{n,\varepsilon}}
\exp\left\{
-\frac{\beta n}{p}z^{p/2}
\right\}
z^{\frac{\beta}{2}-1} 
\bar{\mathbb E}_{n-1}^{(n)}
\left[
\mathbf 1_{E_{n-1}}
\prod_{i=1}^{n-1}|z-Y_i|^\beta
\right]
 \dd z,
\end{align*}
where $\bar Z_{k,p,\beta}^{(n)}$ is defined in \eqref{eq:barzk}. 
On $E_{n-1}$, for $z\in I_{n,\varepsilon}$ and all large $n$, $|z-Y_i|\ge z-R^{2/q}.$
Therefore,
$$
\prod_{i=1}^{n-1}|z-Y_i|^\beta
\ge
\left(z-R^{2/q}\right)^{\beta(n-1)}
=
\exp\{\Theta(n\log n)\}.
$$
Moreover, \eqref{eq:barzz} implies that $
\liminf_{n\to\infty}
\frac{1}{\beta n^{1+p/q}}
\log
\frac{\bar Z_{n-1,p,\beta}^{(n)}}
     {\bar Z_{n,p,\beta}^{(n)}}
\ge0. 
$
Finally, for
$z\in I_{n,\varepsilon}$, we have $z^{p/2}
=
n^{p/q}(s^q-r_{p,q}^q+\varepsilon)^{p/q}
+
o(n^{p/q}).$
Collecting the above estimates, we arrive at
$$
\liminf_{n\to\infty}
\frac1{\beta n^{1+p/q}}
\log
\mathbb P(\tilde B_{q,n}>s^2)
\ge
-\frac1p
\left(s^q-r_{p,q}^q+\varepsilon\right)^{p/q}.
$$
Sending $\varepsilon\downarrow0$ yields the desired lower bound.

Let us turn to the upper bound. Fix $0<\varepsilon<s^q-r_{p,q}^q.$ Choose $R>b_p^q$ so large that
$
\int f_R \dd\eta^{(p)}
\le
r_{p,q}^q+\frac{\varepsilon}{2}.
$
Since $f_R$ is bounded and continuous, the LDP for empirical measure $L_n=\frac{1}{n} \sum_{i=1}^n \delta_{Y_i}$ gives
$$
\limsup_{n\to\infty}
\frac1{\beta n^{1+p/q}}
\log
\mathbb P\left(
\frac1n\sum_{i=1}^n f_R(Y_i)>r_{p,q}^q+\varepsilon
\right)
=
-\infty.
$$
It remains to estimate
$$
\mathcal A_{n,R,\varepsilon}
:=
\left\{
\tilde B_{q,n}>s^2,\
\frac1n\sum_{i=1}^n f_R(Y_i)\le r_{p,q}^q+\varepsilon
\right\}.
$$
For $I\subset\{1,\ldots,n\}$, we write $m=|I|$ and $I^c=\{1,\dots,n\} \setminus I$. For $I=\{i_1,...,i_m\}$, let 
$$\Gamma_I=\{ (Y_{i_1},\dots,Y_{i_m}): \bY=(Y_1,\dots,Y_n)\in \ax_{n,R,\varepsilon}, |Y_{i}|> R^{1/q}, \forall i\in I \mbox{ and } |Y_j|\le R^{1/q},\forall j\in I^c\}$$ 
be the projection of $\ax_{n,R,\varepsilon}$ onto the coordinates from $I$. We decompose $\mathcal A_{n,R,\varepsilon}$ according to the large coordinates of $\bY$ by requiring $i\in I$ exactly when  $Y_i>R^{2/q}$ 
 and deduce
\begin{align*}
\mathbb P(\mathcal A_{n,R,\varepsilon})
\le
\sum_{I\subset\{1,\ldots,n\}}
\int_{\Gamma_I}
&
\exp\left\{
-\frac{\beta n}{p}\sum_{i\in I}z_i^{p/2}
\right\}
\prod_{i\in I}z_i^{\frac{\beta}{2}-1}
\prod_{\substack{i<j\\ i,j\in I}}|z_i-z_j|^\beta
\mathcal B_{n,I}(\bz)
 \dd\bz
\end{align*}
where
\begin{align*}
\mathcal B_{n,I}(\bz)
=
\frac1{\bar Z_{n,p,\beta}^{(n)}}
\int_{[0,R^{2/q}]^{n-m}}
&
\prod_{\substack{i<j\\ i,j\in I^c}}|u_i-u_j|^\beta
\prod_{j\in I^c}u_j^{\frac{\beta}{2}-1}
\exp\left\{
-\frac{\beta n}{p}\sum_{j\in I^c}u_j^{p/2}
\right\} 
\prod_{\substack{i\in I\\ j\in I^c}}|z_i-u_j|^\beta
 \dd\bu.
\end{align*}
Using \eqref{e:pairbd}, \eqref{eq:Frelb}, and Lemma \ref{le:sv-pc}, we estimate the contribution from small coordinates and obtain that for every $\eta>0$ and $R\geq1$,
$$
\mathcal B_{n,I}(\bz)
\le
\exp\left\{
\eta\beta n\sum_{i\in I}z_i^{p/2}
+
C_{\eta} \log(1 + R^{2/q}) nm
\right\}.
$$
For the large coordinates, applying \eqref{e:logpair} with $\alpha=p/2$ and $\varepsilon=\eta$ to the pairs inside $I$, we find
$$
\prod_{\substack{i<j\\ i,j\in I}}|z_i-z_j|^\beta
\le
\exp\left\{
\eta\beta m\sum_{i\in I}z_i^{p/2}
+
C_\eta m^2
\right\}
\le
\exp\left\{
\eta\beta n\sum_{i\in I}z_i^{p/2}
+
C_\eta nm
\right\}.
$$
Moreover, since $z_i>R^{2/q}>0$, we have
$
|\frac{\beta}{2}-1 | |\log z_i|
\le
\eta\beta n z_i^{p/2}+C_\eta n
$
after adjusting $C_\eta$ if necessary. It follows that
$\prod_{i\in I}z_i^{\frac{\beta}{2}-1}
\le
\exp\left\{
\eta\beta n\sum_{i\in I}z_i^{p/2}
+
C_\eta nm
\right\}.$
Therefore,
$$
\mathbb P(\mathcal A_{n,R,\varepsilon})
\le
\sum_{I\subset\{1,\ldots,n\}}
\exp\{C_{\eta,R}nm\}
\int_{\Gamma_I}
\exp\left\{
-\beta n\left(\frac1p-3\eta\right)
\sum_{i\in I}z_i^{p/2}
\right\}
 \dd\bz.
$$
On the event $\mathcal A_{n,R,\varepsilon}\cap\{|Y_i|>R^{2/q},\forall i\in I \mbox{ and } |Y_j|\le R^{2/q}, \forall j\in I^c\}$, we have
$$
\frac1n\sum_{i\in I}\left(z_i^{q/2}-R\right)
=
\frac1n\sum_{a=1}^n y_a^{q/2}
-
\frac1n\sum_{a=1}^n f_R(y_a)
>
s^q-r_{p,q}^q-\varepsilon.
$$
Observing that $\sum_{i\in I}z_i^{p/2}\ge |I| R^{p/q},$
we find $\sum_{i\in I}z_i^{p/2}
\ge
\max\{n^{p/q} \left(s^q-r_{p,q}^q-\varepsilon\right)^{p/q}, |I| R^{p/q}\}.$
Choose $\eta>0$ and $\kappa>0$ so small that $1/p-3\eta-\kappa >0$. Then
\begin{align*}
&\int_{\Gamma_I}
\exp\left\{
-\beta n\left(\frac1p-3\eta\right)
\sum_{i\in I}z_i^{p/2}
\right\}
 \dd\bz \\
\le{}&
\exp\left\{
-\beta n\left(\frac1p-3\eta-\kappa\right)
\max\{n^{p/q}\left(s^q-r_{p,q}^q-\varepsilon\right)^{p/q}, mR^{p/q}\}+ C_\kappa m\log n
\right\},
\end{align*}
since  $\int_{[0,\infty)^m}
\exp\left\{
-\beta n\kappa\sum_{i=1}^m z_i^{p/2}
\right\}
 \dd\bz\le e^{C_\kappa m\log n}$.
Letting $
\delta_R=\frac{(C_\eta+C_\kappa)\log(1+R^{2/q})}{\beta R^{p/q}}
$, we deduce
\begin{align*}
\mathbb P(\mathcal A_{n,R,\varepsilon})
&\le
\sum_{m=1}^n \binom{n}{m}
\exp\Bigg\{
-\beta n\left(\frac1p-3\eta-\kappa\right)
\max\{n^{p/q}\left(s^q-r_{p,q}^q-\varepsilon\right)^{p/q}, mR^{p/q}\} \\
&\hspace{0.3\textwidth}
+C_{\eta} \log(1 + R^{2/q}) nm
+C_\kappa m\log n
\Bigg\}\\
&\le \sum_{m=1}^n \binom{n}{m}
\exp\Big\{
-\beta n\Big(\frac1p-3\eta-\kappa-\de_R\Big)
\max\{n^{p/q}\left(s^q-r_{p,q}^q-\varepsilon\right)^{p/q}, mR^{p/q}\}\Big\}
\end{align*}
Taking limits  in the order $n\to\8, R\to\8,\eta\to0,\kappa\to0$, the rest of proof follows that of Theorem \ref{th:momldp} verbatim.
\end{proof}

\begin{proof}[Proof of Lemma \ref{x:powmax}]
Take $D=\mathbb R_+$ with the decomposition \eqref{e:powwgt}. The proof of Lemma \ref{x:powemp} establishes Assumptions A and B in \cite{WZbeta} for every $\mu\in\mathcal M_1(\mathbb R_+)$ with $I_{W_p}(\mu)<\infty$. Therefore, we only need to verify Assumption C in \cite{WZbeta} to apply \cite{WZbeta}*{Proposition 2.8}. Condition (C0) holds by assumption.

By Remark \ref{rem:squlm}, the equilibrium measure for $W_p$
on $\mathbb R_+$ is $\mu_{W_p}=\eta^{(p)}$
with support $[0,b_p^2]$. Hence $b_{W_p}=b_p^2\in\mathcal O$.
By \eqref{eq:arcavg}, we have
\begin{align*}
F_{W_p}
=\frac1{\pi p}\int_0^{b_p^2}
\frac{x^{p/2}}{\sqrt{x(b_p^2-x)}}\dd x
-\log\frac{b_p^2}{4}
=\frac2p-2\log\frac{b_p}{2}.
\end{align*}
Since $\eta^{(p)}=(t\mapsto t^2)_\#\nu^{(p)}$, by change of variable in the definition \eqref{eq:phipq} we have for $r>b_p^2$
\begin{align*}
J_{W_p}(r)
&=\frac{r^{p/2}}{p}
-\int_{\mathbb R_+}\log|r-x| \eta^{(p)}(\dd x)-F_{W_p}\\
&=\frac{r^{p/2}}{p}
-2\int_{\mathbb R}\log|\sqrt r-t| \nu^{(p)}(\dd t)-F_{W_p}
=2\Phi_{p,\infty}(\sqrt r);
\end{align*}
see \eqref{e:powlog} below for details of the last equality. 
Both sides extend continuously to zero at $r=b_p^2$.
Differentiating yields
\begin{align*}
J_{W_p}'(r)
=\frac12r^{p/2-1}-\frac{p}{b_p^p\sqrt r}
\int_0^{b_p}\frac{t^{p-1}}{\sqrt{r-t^2}}\dd t
=\frac{p}{b_p^p\sqrt r}
\int_{b_p}^{\sqrt r}\frac{t^{p-1}}{\sqrt{r-t^2}}\dd t>0.
\end{align*}
Thus \eqref{e:powmax} holds and (C1) is verified.

 For (C2), the only finite endpoint is $0$.  
For every $\delta>0$, we have
\begin{align*}
\int_0^\delta e^{-\beta nW_p(x)}\omega_n(x)\dd x
\leq\int_0^\delta x^{\beta/2-1}\dd x
=\frac{2\delta^{\beta/2}}{\beta}.
\end{align*}
Since $\beta$ is fixed and $W_p(0)=0$, this verifies \cite{WZbeta}*{Eq.~(2.56)} and thus Remark 2.7 of that paper gives (C2). 
Noting that $(\beta n)^{-1}\log\omega_n(x) \leq\log(1+x)$ for $x\geq1$, \eqref{e:powgro} yields
\begin{align*}
W_p(x)-\frac1{\beta n}\log\omega_n(x)
&\geq W_p(x)-\log(1+x)
\geq\frac{x^{p/2}}{2p}-C,\qquad x\geq1.
\end{align*}
By \cite{WZbeta}*{Remark 2.9}, This uniform bound implies (C3).
The asserted LDP now follows from \cite{WZbeta}*{Proposition 2.8}.
\end{proof}

\bibliography{references}

\end{document}